\hoffset=0.2in
\magnification=1150
\input amstex
\documentstyle{amsppt}
\nologo

\define\occ{\operatorname {CC}}
\define\ct{T^*X}
\define\db#1{\operatorname {D}^b(#1)}
\define\oh{{\operatorname H}}
\define\ohf#1{{\oh^{inf}_{#1}}}
\define\oc#1{{\operatorname C_{#1}^{inf}}}
\define\oz#1{{\operatorname Z_{#1}^{inf}}}
\define\e{{\epsilon}}
\define\fg{{\frak g}}
\define\fa{{\frak a}}
\define\fk{{\frak k}}
\define\fb{{\frak b}}
\define\fm{{\frak m}}
\define\fn{{\frak n}}
\define\fh{{\frak h}}
\define\fc{{\frak c}}
\define\ft{{\frak t}}
\define\fp{{\frak p}}

\define\fv{{\frak v}}
\define\CB{{\Cal B}}

\define\br{{\Bbb R}}
\define\bc{{\Bbb C}}
\define\AR{{A_\Bbb R}}
\define\CR{{C_\Bbb R}}
\define\TR{{T_\Bbb R}}
\define\ttr{{\tilde T_\Bbb R}}
\define\GR{{G_\Bbb R}}
\define\UR{{U_\Bbb R}}
\define\KR{{K_\Bbb R}}
\define\gr{{\fg_{\Bbb R}}}
\define\kr{{\fk_{\Bbb R}}}

\define\ug{{\Cal U(\fg)}}
\define\zg{{\Cal Z(\fg)}}
\define\cf{{\Cal F}}
\define\cg{{\Cal G}}
\define\h{\hat}
\define\TP{{\Theta_\pi}}
\define\tp{{\theta_\pi}}

\define\pro{G_{\Bbb R}\times H \times \hat X}

\define\Hom{\operatorname{H}}  

\topmatter
\title Two geometric character formulas \\ for reductive Lie groups
\endtitle
\author Wilfried Schmid and Kari Vilonen
\endauthor
\address {Department of mathematics, Harvard University, Cambridge, MA
02138, USA}
\endaddress
\email schmid\@math.harvard.edu
\endemail
\address{Department of mathematics, Brandeis University, Waltham, MA 02254,
USA}
\endaddress
\email vilonen\@math.brandeis.edu
\endemail
\thanks W.Schmid was partially supported by NSF
\endgraf K.Vilonen was partially supported by NSA and NSF
\endthanks
\endtopmatter

\document

\subheading{\bf 1. Introduction}

\vskip 5\jot

In this paper we prove two formulas for the characters of 
representations of reductive groups. Both express the character of a 
representation $\pi$ in terms of the same geometric data attached to $\pi$. 
When specialized to the case of a compact Lie group, one of them 
reduces to Kirillov's character formula in the compact case, and the 
other, to an application of the Atiyah-Bott fixed point formula to 
the Borel-Weil realization of the representation $\pi$.

To set the stage, let us first recall the Borel-Weil theorem for a 
connected, compact Lie group $\GR$. For simplicity, we assume that $\GR$ 
is simply connected. We let $\fg_\br$ denote its Lie algebra and 
$\fg=\Bbb C\otimes_\Bbb R \fg_\br$ the complexified Lie algebra. Via the 
adjoint action, $\GR$ operates on $i\gr^*$, the space of all $\br$-linear 
functions $\lambda : \fg_\br \to i\br$. Every $\GR$-orbit $\Omega$ in 
$i\gr^*$ -- \lq\lq coadjoint orbit" for short -- carries a canonical 
$\GR$-invariant symplectic structure. The orbit $\Omega$ is said to be 
integral if some, or equi\-valently any, $\lambda\in\Omega$ exponentiates 
to a character $e^\lambda$ of the isotropy subgroup $(\GR)_\lambda$. In that 
case, the one dimensional representation spaces of the characters 
$e^\lambda$, $\lambda\in\Omega$, fit together into a $\GR$-equivariant, 
real algebraic, Hermitian line bundle $\bold L_\Omega \to \Omega$. The 
pair $(\Omega,\bold L_\Omega \to \Omega)$ carries a unique \lq\lq positive 
polarization": a $\GR$-invariant complex structure on the manifold 
$\Omega$ and the structure of $\GR$-equivariant holomorphic line bundle on 
$\bold L_\Omega$, positive in the sense of algebraic geometry. The hypothesis 
of simple connectivity ensures that  the square root $\sqrt{\bold
K_\Omega}$ of the canonical bundle 
$\bold K_\Omega$ exists as $\GR$-equivariant holomorphic line bundle. 
According  to the Borel-Weil theorem, the natural action of $\GR$ on the space
of  holomorphic $\bold L_\Omega$-valued half forms 
$\operatorname H^0(\Omega, \Cal O(\bold K_\Omega^{1/2}\otimes \bold L_\Omega))$
is  irreducible; moreover, the association
$$
\Omega \ \rightsquigarrow \ 
\operatorname H^0(\Omega, \Cal O(\bold K_\Omega^{1/2}\otimes \bold L_\Omega))
\tag1.1
$$
sets up a bijection between the integral regular -- i.e. maximal 
dimensional -- coadjoint orbits $\Omega \subset i\fg_\br$ and the 
irreducible representations of $\GR$. We should remark that   the 
passage from sections to half forms corresponds to the $\rho$-shift 
familiar in representation theory.

Still in the setting of a compact group, we write $\Theta_\Omega$ for 
the character of the representation associated to the $\Omega$, and 
$\theta_\Omega$ for the pullback, as half form, of $\Theta_\Omega$ to 
the Lie algebra $\fg_\br$; concretely,
$$
\theta_\Omega \ = \ (\operatorname{det}\operatorname{exp}_*)^{1/2}
\,\operatorname{exp}^*\Theta_\Omega \,.
\tag1.2
$$
We define the Fourier transform $\hat \phi$ of a test function $\phi \in 
C^{\infty}_c(\fg_\br)$ without the usual factor $i$ in the exponent,
$$
\hat \phi(\zeta) \ = \ \int_{\gr}\, \phi(x)\,e^{\langle \zeta, 
x \rangle}\,dx \qquad\qquad(\,\zeta \in i\fg_\br^*\,)\,,
\tag1.3
$$
as function on $i\fg_\br^*$. Let $\sigma_\Omega$ denote the canonical
symplectic form and $n$ the complex dimension of $\Omega$. With these
conventions, Kirillov's formula can be stated as follows:
$$
\int_\gr \theta_\Omega\,\phi\,dx \ = \
\frac 1{(2\pi )^n\, n!}\,\int_\Omega\hat\phi\,\sigma_\Omega^n\,.
\tag1.4
$$
In effect, the formula expresses the Fourier transform $\hat \theta_\Omega$ as
integration over $\Omega$. 

Translation $\ell_g:\Omega\to\Omega$ by a generic $g\in\GR$ has finitely many
fixed points. The Atiyah-Bott fixed point formula gives the Lefschetz number 
$$
L(g) \ = \ \sum_p \ (-1)^p \ \operatorname{tr}\{(\ell_g)^*: \operatorname
H^p(\Omega, \Cal O(\bold K_\Omega^{1/2}\otimes \bold L_\Omega)) \to
\operatorname H^p(\Omega, \Cal O(\bold K_\Omega^{1/2}\otimes \bold L_\Omega))\}
$$
as a sum over the fixed point set $\Omega^g$, 
$$
L(g) \ = \ \sum_{\zeta\in\Omega^g} \  \frac{\operatorname{tr}\{(\ell_g)^*:(\bold
K_\Omega^{1/2}\otimes \bold L_\Omega)_\zeta \to (\bold K_\Omega^{1/2}\otimes
\bold L_\Omega)_\zeta\}}{\operatorname{det}\{(\operatorname{Id} - \ell_{g*}):
T_\zeta\Omega\to T_\zeta\Omega\}}\ .
$$
The positivity of the polarization implies the vanishing of the higher
cohomology groups $\,\operatorname
H^p(\Omega, \Cal O(\bold K_\Omega^{1/2}\otimes \bold L_\Omega))$\,, \ $p>0$.
Hence 
$$
\Theta_\Omega(g) \ = \ \sum_{\zeta\in\Omega^g} \ 
\frac{\operatorname{tr}\{(\ell_{g^{-1}})^*:(\bold K_\Omega^{1/2}\otimes \bold
L_\Omega)_\zeta \to (\bold K_\Omega^{1/2}\otimes
\bold L_\Omega)_\zeta\}}{\operatorname{det}\{(\operatorname{Id} -
(\ell_{g^{-1}})_*): T_\zeta\Omega\to T_\zeta\Omega\}}\ ,
\tag1.5
$$
for every generic $g\in\GR$. This conclusion translates readily into the Weyl
character formula -- the Atiyah-Bott fixed point formula provides a geometric
interpretation of Weyl's formula. 

Our objective in this paper is the generalization of the two formulas (1.4) and
(1.5) to the case of a reductive Lie group. When the group is noncompact, no
simple, completely explicit character formula is known in general, certainly no
formula as simple as the Weyl character formula in the compact case. Our two
formulas express the character of any irreducible unitary representation $\pi$
-- more generally,  of any admissible  $\pi$ having an infinitesimal character
-- explicitly in geometric terms. The formulas, though simple in appearance, do
not seem to lead to simple numerical expressions for the values of characters:
in effect, the combinatorial complexity of the character values is reflected by
the complexity of the geometry. 

In the following, we consider a connected, linear, reductive Lie group $\GR$,
with Lie algebra $\gr$. Such a group has a complexification $G$, whose Lie
algebra $\fg$ is naturally isomorphic to the complexification of $\gr$. The
complex algebraic group $G$ acts algebraically and transitively on the flag
variety $X$, i.e., the variety of Borel subalgebras $\fb\subset\fg$. The real
group $\GR$ acts real algebraically on $X$, with finitely many orbits
\footnote{When $\GR$ happens to be compact, it acts transitively, and
$X\cong \Omega$ as $\GR$-homogenous complex manifolds, for every regular
coadjoint orbit $\Omega\subset i\gr^*$.}. Representations of $\GR$ can be
associated to twisted $\GR$-equivariant sheaves $\cf$ on $X$ \cite{KSd}, by an
analytic process resembling  Beilinson-Bernstein localization \cite{BB1,2}.
This analytic construction, which is recalled in detail in
\S 2, attaches cohomology groups  to $\cf$, on which $\GR$ acts by translation,
continuously with respect to a certain canonical topology. The resulting
representations are admissible, and have an infinitesimal character which is
determined by the twisting parameter. Every admissible representation with
infinitesimal character can be realized in this way, up to infinitesimal
equivalence; the realization can be chosen in degree zero, with all other
cohomology groups vanishing. Formally, this is analogous to the Borel-Weil
realization (1.1): in the compact case,
$\cf$ reduces to the constant sheaf, and the twisting parameter specifies the
line bundle. We write
$\Theta(\cf)$ for the alternating sum of the Harish-Chandra characters of the
cohomology groups attached to $\cf$.  These are conjugation invariant functions
on the regular set in $\GR$. The $\Theta(\cf)$ exhaust the set of characters of
admissible representations with infinitesimal character, as follows from what
was just said. As in (1.2), the character $\Theta(\cf)$ has a Lie algebra
analogue $\theta(\cf)$, which determines $\Theta(\cf)$ at least near the
identity. 

Our first formula expresses the Fourier transform of $\theta(\cf)$ as an
integral over a certain cycle in the cotangent bundle $\ct$. The idea that this
can be done in principle is due to Rossmann, who also worked out a certain
special case \cite{R3}. A precise description of the relation between his work
and ours, as well as some other references, can be found at the beginning of
section \S 3 below. Our formula involves the negative of the twisting parameter
of $\cf$, which we denote by $\lambda$. Rossmann has defined a twisted moment
map
$$
\mu_\lambda\  :\ \ct \ \longrightarrow \ \fg^* \,,
\tag1.6
$$
depending on $\lambda$. When $\lambda$ happens to be regular, $\mu_\lambda$
maps $\ct$ real algebraically and isomorphically on $\Omega_\lambda=G\cdot
\lambda$, the coadjoint orbit through $\lambda$ for the complex group $G$; at
the opposite extreme, for $\lambda = 0$, the map $\mu_\lambda$ reduces to the
ordinary moment map, and takes values in the nilpotent cone $\Cal N
\subset \fg\cong\fg^*$. Let $\phi\in C_c^\infty(\gr)$ be a test function, and
$\h \phi$ its Fourier transform, normalized as in (1.3), but viewed as a
holomorphic function on $\fg^*$. We can now state our first formula:
$$
\int_{\gr}\theta(\cf)\,\phi\,dx \ = \ \frac 1 {(2\pi i)^n n!}
\int_{\occ(\cf)} \mu^*_\lambda \hat \phi \ (-\sigma + \pi^* \tau_\lambda)^n\,.
\tag1.7
$$
Here $\sigma$ denotes the canonical (holomorphic) symplectic form on $\ct$,
$\pi^*\tau_\lambda$ is the pullback from $X$ to $\ct$ of a two form
representing the  class determined by $2\pi i\lambda$ in
$\operatorname{H}^2(X,\Bbb C)$, and $n$ the complex dimension of $X$. Most
importantly, $\occ(\cf)$ refers to the characteristic cycle of the sheaf $\cf$
as defined by Kashiwara \cite{Ka2}. It is a top dimensional cycle, of possibly
infinite support, carried by $T^*_\GR X$, the union of conormal bundles of
all the $\GR$-orbits on $X$. 

For $\cf = \Bbb C_X=$ constant sheaf on $X$,
the characteristic cycle $\occ(\Bbb C_X)$ coincides with the zero section $X$ in
$\ct$; in that case, the formula (1.7) reduces directly to Kirillov's formula
(1.4). More generally, for $\lambda$ regular but $\cf$ arbitrary, the
integration on the right in (1.7) can be performed  over the middle dimensional 
cycle $\mu_\lambda(\occ(\cf))$ in the complex coadjoint
orbit $\Omega_\lambda$, 
$$
\int_{\gr}\theta(\cf)\,\phi\,dx \ = \ \frac 1 {(2\pi i)^n n!}
\int_{\mu_\lambda(\occ(\cf))} \ \hat \phi \ \sigma_\lambda^n\,,
\tag1.8
$$
with $\sigma_\lambda$ denoting the natural (holomorphic) symplectic form on
$\Omega_\lambda$. When $\cf$ corresponds to a tempered irreducible
representation, $\mu_\lambda(\occ(\cf))$ turns out to be homologous to a
coadjoint orbit of $\GR$ in $i\gr^*$, and (1.8) becomes equivalent to
Kirillov's \lq\lq universal character formula" \cite{Ki3} in this special case,
which was established by Rossmann
\cite{R1}. For non-tempered representations of the reductive group $\GR$, the
\lq\lq universal character formula" fails, since there are not enough coadjoint
orbits. It was Rossmann's idea to express irreducible characters
not necessarily   as integrals  over coadjoint orbits, but more generally
as integrals over cycles in complexified coadjoint orbits. In this sense, the
formulas (1.7-8) can be regarded as remedies for the failure of the  \lq\lq
universal character formula" in the reductive case.

Our second formula was conjectured by Kashiwara \cite{Ka4}. It comes in several
flavors: global, pointwise, both as formulas on the group and on the Lie
algebra. Here, in this introduction, we describe two of the four versions. Let
$\tilde G$ denote the set of pairs $(g,x)$ with $g\in G$ and $x$ in the fixed
point set $X^g$. This is naturally a smooth complex algebraic variety.
Projection to the first factor defines a map $q:\tilde G \to G$, which is a
covering over the regular set in $G$, with covering group equal to the Weyl
group $W$. The exponentials $e^\alpha$ of the roots $\alpha$ exist as
multiple valued holomorphic functions on $G$; when pulled back to $\tilde G$,
they become globally defined algebraic functions. By the same process, the
exponential  of the twisting parameter becomes a well defined
algebraic function on $\tilde G$ when $\lambda$ is integral, and a multiple
valued function in general. In particular, this exponential generates a rank one
local system $\Bbb C_\lambda\subset \Cal O_{\tilde G}$. Kashiwara, in the case
of trivial twisting parameter \footnote{That is, $\lambda=\rho$ in
representation theoretic notation. This corresponds to the infinitesimal
character of the trivial representation, and makes the  twisted sheaves $\cf$
into ordinary ones.}, associates to the sheaf $\cf$ a top dimensional cycle
$c(\cf)$ in $\tilde G_\Bbb R = q^{-1}(\GR)$. The crux of the matter is a fixed
point formalism, which we extend to the general, twisted case in \S 5  below.
The result is a top dimensional cycle $c(\cf)$ in  $\tilde G_\Bbb R$, whose
coefficients are not integers, but sections of the local system $\Bbb
C_\lambda$.  This is the global version, on the group, of our second
formula:
$$
\int_\GR \Theta(\cf) \phi\, dg \ = \ \int_{c(\cf)}(q^*\phi)\tilde
\omega\qquad\text{for $\phi\in C^\infty_c(\GR)$}\,; 
\tag1.9
$$  
here $\tilde \omega$ denotes the holomorphic form on $\tilde G$ obtained by
complexifying  the Haar measure $dg$ of $\GR$ to a holomorphic form $\omega$
of top degree on
$G$, and dividing the pull back $q^*\omega$  by the function $\prod_{\alpha >0}
(1-e^{-\alpha})$, which vanishes to first order on the singular locus of $q$. 

Letting $\phi$ run through a sequence converging (weakly) to the delta 
function
$\delta_g$ at any particular regular $g\in \GR$, we can evaluate both sides of
the identity (1.9) on $\phi=\delta_g$. The result is a cohomological expression
for $\Theta(\cf)(g)$ in terms of the action of $g$ at its various fixed points
-- an analogue, in the noncompact case, of (1.5). This local formula for
$\Theta(\cf)$ has a counterpart on the Lie algebra, as does the global formula
(1.9). For simplicity, we state only the Lie algebra version of the local
formula here, in the introduction -- both local versions can be found in \S 5.

We consider a regular $\zeta\in\gr$,  and write $E$ for the connected 
component of
$\zeta$ in the regular set. Each fixed point $x\in X^\zeta$ of the
infinitesimal action of
$\zeta$ -- in other words, a zero of $\zeta$ considered as  vector field on
$X$ -- determines a pair $(\zeta, x)\in \tilde\fg_\Bbb R$, the Lie algebra
analogue of the space
$\tilde G_\Bbb R$. This pair, in turn, gives meaning to the value 
$\alpha_x(\zeta)$, for every positive root $\alpha$, and similarly to
$\lambda_x(\zeta)$.  On general grounds, there exist constants 
$d_{E,x}$, such that
$$
\theta(\cf)(\zeta) \ = \ \sum_{x\in X^{\zeta}}\ \frac { d_{E,x}\,
e^{\lambda_x}(\zeta)} {\tsize \prod_{\alpha \in {\Phi}^{ +}}
\, \alpha_x(\zeta)} \ .
\tag1.10a
$$
This is Harish-Chandra's local expression for invariant eigendistributions on
the Lie algebra. The Lie algebra version of (10.9), evaluated on the delta
function at $\zeta$, computes the $d_{E,x}$, and thus makes (1.10a) explicit:
$$   
d_{E,x}\  = \ \sum_p \ (-1)^p \dim \,\operatorname {H}_c^p(N'(\zeta,x)\cap
D_\epsilon,\cf)  \,, \qquad 0<\epsilon \ll 1 \,.
\tag1.10b
$$
In this formula, $\operatorname {H}_c$ refers to cohomology with compact
support, $D_\epsilon$ is the $\epsilon$-ball at $x$, and $N'(\zeta,x)$ denotes a
subspace of $X$ obtained by extending the space of shrinking directions of
$\zeta$ at $x$ by means of a suitable unipotent subgroup of $\GR$.

We should point out that the group version of the fixed point formula (1.10)
implies, and is implied by, a statement known as Osborne's conjecture. The
argument for the equivalence can be found in the note \cite{SV1}, which
announces our proof of Kashiwara's conjecture. The integral formula (1.7-8)
and its proof were announced in \cite{SV2}; it is a major ingredient of our proof
of a conjecture of Barbasch and Vogan \cite{SV3}. The  formulas (1.7) and
(1.10) provide two radically different expressions for the same quantity
$\theta(\cf)$. One may wonder if it is possible to go directly between them,
without the \lq\lq detour" via representation theory. In the compact case, the
equivalence of (1.4) and (1.5) was established by Berline and Vergne,
who established their localization formula for this purpose \cite{BV}. The
note \cite{S5} speculates on the possibility of deducing (1.10) from (1.7) by
a hypothetical localization formula for non-elliptic  fixed points.

The proofs of the two formulas follow the same pattern, though the details are
very different. We develop them side-by-side; however, they can be read
separately. Section 2, which recalls the description of representations in
terms of $\GR$-equivariant sheaves, is equally relevant for both formulas.
Section 3 develops the statement of the integral formula, and the two
subsequent sections, that of the fixed point formula.  Both the statement and
the proof of the fixed point formula would simplify greatly if we limited
ourselves to the case  $\lambda=\rho$, as in Kashiwara's statement of his
conjecture \cite{Ka4}; his note provides a good introduction to the fixed point
formula. The actual proofs occupy sections 7-10. Except for \S 8, which is
short, each of the last four sections begins with an introduction common to
both formulas, then treats the integral formula, and finally, the fixed point
formula.

\vskip 1cm

\subheading{\bf 2. Geometric data}
\vskip .5cm

Let us begin with the notation and hypotheses that will be in force for the
rest of the paper. We fix a connected, complex algebraic, reductive group
$G$ which is defined over $\br$. The representations we consider will be
representations of a real form $\GR$  of $G$ -- in other words, $\GR$ is a
subgroup of
$G$ lying between the group of real points $G(\br)$ and the identity
component
$G(\br)^0$. We regard $\GR$ as a reductive Lie group. All of our results
remain valid in the larger class of reductive Lie groups considered by
Harish-Chandra
\cite{HC7, \S 3}; we shall  comment later on certain   modifications 
necessary to cover this larger class of groups. We pick a maximal compact
subgroup $\KR$ of $\GR$ and recall that this is not an essential choice: any
two maximal compact subgroups are conjugate. The complexification $K$ of
$\KR$ lies naturally as a subgroup in $G$. We denote the Lie algebras of the
four groups by the corresponding lower case German symbols $\fg$, $\gr$,
$\kr$, and $\fk$; the latter three are subalgebras of $\fg$.

By a representation of $\GR$, we shall always mean a continuous
representation on a complete, locally convex, Hausdorff topological vector
space. Such a representation is called admissible if its restriction to $\KR$
involves any irreducible representation of $\KR$ only finitely often. A
representation of
$\GR$ is said to have finite length if every increasing chain of closed,
invariant subspaces breaks off after finitely many steps. The universal
enveloping algebra $\ug$, and hence its center $\zg$, operates on the dense
subspace of all $\KR$-finite vectors of any particular admissible
representation $\pi$ of finite length -- the so-called Harish-Chandra module
of $\pi$. When $\zg$ acts on the Harish-Chandra module via a character, one
says that $\pi$ has an infinitesimal character. Every irreducible, admissible
representation $\pi$  does have an infinitesimal character; this follows from
the irreducibility of the Harish-Chandra module of $\pi$ \cite{HC2}.

A construction of Harish-Chandra attaches a (global) character $\TP$ to
every admissible representation $\pi$ of finite length -- a conjugation
invariant distribution
\footnote{In this paper, \lq\lq distribution" will be synonymous to \lq\lq
generalized function", i.e., continuous linear functional on the space of smooth,
compactly supported measures. Thus, under coordinate changes,
distributions transform like functions.} on
$\GR$. The characters of any set $\Pi$ of irreducible, admissible
representations are linearly independent, provided no two representations
in
$\Pi$ are infinitesimally equivalent. Infinitesimal equivalence, we recall, is
the equivalence relation defined by isomorphism of the underlying
Harish-Chandra modules. The character is an additive invariant in
short exact sequences. It follows that the composition factors of any
admissible representation $\pi$ of finite length are determined, up to
infinitesimal equivalence,  by the character $\TP$ 
\cite{HC3}.  

We identify
$\ug$  with the algebra of left invariant linear differential operators on
$\GR$ via infinitesimal right translation. Under this identification, $\zg$
corresponds to the algebra of all bi-invariant linear differential operators. 
As such,
$\zg$ acts on functions and distributions on the group $\GR$. When a
representation $\pi$ has an infinitesimal character $\chi$, the algebra $\zg$
operates on the distribution $\TP$ via this same character $\chi$. In
Harish-Chandra's terminology, $\TP$ is then an invariant eigendistribution: a
conjugation invariant distribution which is an eigenvector for the algebra
$\zg$.

Harish-Chandra's regularity theorem \cite{HC4} asserts that every invariant
eigendistribution $\Theta$ is a locally $L^1$ function. This locally $L^1$
function is real analytic on $G'_\Bbb R$, the set of regular semisimple
elements in $\GR$. The complement of $G'_\Bbb R$ has measure zero, so the
restriction of
$\Theta$ to the open subset $G'_\Bbb R\subset \GR$  completely determines
$\Theta$. Consequently, all finite linear combinations of invariant
eigendistributions -- and in particular all characters of admissible
representations of finite length -- may be thought  of as real analytic
functions on $\GR$ with potential singularities along the complement of 
$G'_\Bbb R$.

The exponential map from $\gr$ to $\GR$ is real analytic. Near 0, it maps
$\gr'$, the set of regular semisimple elements in $\gr$, into $\GR'$.  More
precisely, 
$$
\gathered
\exp(\gr'\cap U) \subset \GR'\,, \ \ \ \text{where}   
\\ 
U=\{ x\in \gr \mid
\text{all eigenvalues $\eta$ of $ \operatorname{ad}(x)$ satisfy
$|\operatorname{Im}\eta|<2\pi$ }\}\,.
\endgathered
$$
Thus, for any character $\Theta_\pi$, we can define
$$
\tp =
\sqrt{\operatorname{det}(\operatorname{exp}_*)}\,\operatorname{exp}^*\TP\,,
\tag2.1
$$ 
at least as  real analytic function on $\gr'\cap U$. In fact, the proof of the
regularity theorem shows that $\tp$ extends real analytically to all of $\gr'$.
The resulting function is  locally $L^1$  as  function on $\gr$,  hence a
globally defined distribution -- a conjugation invariant distribution, since the
exponential map commutes with conjugation. When the representation $\pi$
has an infinitesimal character, $\tp$ is an eigendistribution for $\zg\cong
S(\fg)^G$, the algebra of conjugation invariant, constant coefficient
differential operators on $\gr$ \cite{HC4}. This is the reason for the factor
$\sqrt{\operatorname{det}(\operatorname{exp}_*)}$ in (2.1):  without it, the
preceding statement would not be correct. We shall call $\tp$ the character
of $\pi$ on the Lie algebra. It determines $\TP$ near the identity, at least, since
the exponential map is a local diffeomorphism near the origin. 

Our character formulas express $\TP$ and $\tp$ in terms of  certain
geometric data attached to the representation $\pi$. Let us describe these
geometric data. The flag variety $X$ of $\fg$ carries a tautological bundle
$\CB \to X$, whose fiber at any $x\in X$ is the Borel subalgebra $\fb_x
\subset \fg$ which fixes $x$.  The various quotients $\fb_x/[\fb_x,\fb_x]$ are
canonically isomorphic, and hence $\CB/[\CB,\CB]$ has a canonical
trivialization; 
$$
\fh \, =_{\text{def}} \, \text{fiber of the trivialized bundle $\CB/[\CB,\CB]$}
\tag2.2
$$ is called the universal Cartan algebra. It comes equipped with a root
system
$\Phi\subset\fh^*$ and an action of the Weyl group $W$. Every concrete
Cartan subalgebra $\fc\subset\fg$ has as many fixed points  $x\in X$ as the
order of
$W$; the choice of a particular fixed point $x$ determines  concrete
isomorphisms
$$
\fc \,\cong \,\fb_x/[\fb_x,\fb_x] \,\cong \,\fh\,.
\tag2.3
$$ We specify a notion of positivity in the universal root system
$\Phi$ so that the weights of $\fg/\fb_x$ constitute the set positive roots
$\Phi^+$. In analogy to the universal Cartan algebra (2.2), we define the
universal Cartan group $H$ for $G$,
$$ H\,\cong\, B_x/[B_x,B_x]\,,
\tag2.4
$$ where $B_x$ denotes the Borel subgroup of $G$ stabilizing $x$. This group
has Lie algebra $\fh$ and contains 
$$ Z \ = \ \text{center of $G$}
\tag2.5
$$ canonically, i.e., independently of the choice of $x$.

Beilinson-Bernstein \cite{BB3} introduce the enhanced flag variety of $\fg$ as
the fiber bundle  $\hat X_{\fg} \to X$ whose fiber at $x\in X$ is a collection of
generators of the $\alpha$-root spaces in $\fg/\fb_x$ for all simple roots
$\alpha \in \Phi^+$. It is visibly a principal bundle over $X$ with structure
group equal to the product of copies of $\bc^*$ indexed by the simple roots.
The universal Cartan group $H$ acts on the simple root spaces, hence on $\hat
X$, and $Z$ is the kernel of this action. One can therefore identify 
$$ H/Z\ = \ \text{structure group of $\hat X_{\fg} \to X$}\,.
\tag2.6
$$ The group $G$ acts equivariantly on $\hat X_{\fg} \to X$, transitively on
both $\hat X_{\fg}$ and
$X$, and the action on $\hat X_{\fg}$ commutes with the action of $H$. To see
all of this more concretely (but less invariantly), we visualize $X$ as the
quotient
$G/B_x$  and
$$
\hat X_{\fg} \ \cong \ G/(ZN_x)\  \cong \ Z\backslash G/N_x\,,
\tag2.7a
$$  with $N_x= \text{unipotent radical of $B_x$}$. Here $H$ acts on $\hat
X_{\fg}$ by
$$ h: g(ZN_x) \mapsto gh^{-1}(ZN_x)\,.
\tag2.7b
$$  The surjection $G \to Z\backslash G$ determines a principal bundle 
$$
\hat X \to X\,, \  \text{with structure group $H$}\,,
\tag2.8
$$ which lies over $\hat X_{\fg} \to X$ with fiber $Z$. We shall call $\hat X$ the
enhanced flag variety of $G$. Via the identification $X \cong G/B_x$, we get
the description
$$
\hat X \ \cong \ G/N_x\,, \ \ \ \text{with $H$-action} \ \ \ h: gN_x \mapsto
gh^{-1}N_x\,,
\tag2.9
$$ analogous to (2.7a,b). 

In the definition of the twisted equivariant derived category which is about
to follow, we shall work with $\hat X$, unlike Beilinson-Bernstein \cite{BB3},
who use $\hat X_{\fg}$ instead. One can show that both
choices result in the same equivariant derived category. When $\GR$ is
non-linear and semisimple, as in \cite{BB3}, $\hat X_{\fg}$ is the better
choice, since it is associated to the Lie algebra, rather than a specific
complex group with Lie algebra $\fg$. On the other hand, if $\GR$ is linear but
with  center of positive dimension,  $\hat X$ is preferable. In the general
case, with $\GR$ reductive and not necessarily linear, neither $\hat X_{\fg}$
nor $\hat X$ is completely satisfactory; instead  one may want to work with
yet another bundle over $X$, with structure group $\fh$.

In the discussion below, we keep fixed a particular \lq\lq localization
parameter" $\lambda \in \fh^*=\text{dual space of $\fh$}$. As is customary,
we set 
$$\tsize
\rho \ = \ \frac 1 2 \sum_{\alpha \in \Phi^+} \alpha\,.
\tag2.10
$$ Let us define the
\lq\lq $\GR$-equivariant derived category on $X$ with twist
$(\lambda-\rho)$", to be denoted by $\operatorname D_\GR(X)_\lambda$.
Technically speaking, it is not a derived category of sheaves on $X$, but
rather a pre-stack -- concretely, a subcategory of the $(\GR \times
\fh)$-equivariant derived category on $\hat X$. For $\lambda=\rho$,
$\operatorname D_\GR(X)_\rho = \operatorname D_\GR(X)$ is the bounded
$\GR$-equivariant derived category in the sense of Bernstein-Lunts
\cite{BL}; a short summary of their construction can be found in
\cite{MV}. Our shift by $\rho$ serves the purpose of making the labeling
compatible with Harish-Chandra's description of the characters of $\zg$
\cite{HC1}.

We let $\GR\times \fh$ act on $\hat X$ via the $(G\times H)$-action, the
inclusion $\GR \hookrightarrow G$, and the exponential homomorphism
$\exp: \fh
\to H$. A $\bc$-sheaf $\cf$ on $\hat X$ is said to be
$(\lambda-\rho)$-monodromic with respect to the $H$-action on $\hat X \to
X$ if it is locally constant on every fiber, with monodromy
$e^{\lambda-\rho}$. To clarify what we mean, note that each fiber can be
identified with $H$ (canonically up to translation), so the restriction of $\cf$ to
any fiber pulls back to a locally constant sheaf on
$H$. The $(\lambda-\rho)$-monodromicity condition on $\cf$ requires that the
locally constant sheaves on $H$, corresponding to the various fibers, have
the same monodromy as the rank one local system generated by the 
multiple valued function  $e^{\lambda-\rho}$. 

The $(\lambda-\rho)$-monodromic sheaves constitute an abelian category
$Sh_{X,\lambda}$. One can thus form the bounded derived category $\db
{Sh_{X,\lambda}}$ and the full subcategory $\operatorname
D^b_c({Sh_{X,\lambda}})$ of complexes whose cohomology is constructible
with respect to a subanalytic stratification. We then pass from
$\operatorname D^b_c({Sh_{X,\lambda}})$ to the bounded
$\GR$-equivariant derived category $\operatorname D_\GR(X)_\lambda$ as
described in \cite{BL}. A seemingly different description of the twisted
equivariant derived category appears in \cite{B,MUV}; one can  show
that the two definitions agree.

Often we shall view objects in  $\operatorname D_\GR(X)_\lambda$ as
complexes of sheaves on $\hat X$, disregarding the $\GR$-equivariance.
Since $\GR\times H$ acts real algebraically on $\hat X$, with finitely many
orbits \cite{Ma,W}, the orbits define a semi-algebraic (Whitney) stratification
of $\hat X$. Any
$\cf\in\operatorname D_\GR(X)_\lambda$, viewed as an element in $\db
{Sh_{X,\lambda}}$, has locally constant cohomology along the orbits, i.e.,
$$
\gathered
\text{the image of $\cf\in\operatorname D_\GR(X)_\lambda$ in  $\db
{Sh_{X,\lambda}}$ under the forgetful functor}
\\
\text{is constructible with respect to the $(\GR\times H)$-orbit
stratification.}
\endgathered
\tag2.11
$$ In particular, the construction of $\operatorname D_\GR(X)_\lambda$ can
be carried out in the semi-algebraic context. The twisting disappears when
$\lambda = \rho$, hence
$$
\operatorname D_\GR(X)_\rho \ = \  \operatorname D_\GR(X)\,.
\tag2.12
$$ The construction of $\operatorname D_\GR(X)_\lambda$ involves only the
monodromy of the multiple valued function
$e^{\lambda-\rho}$ on
$H$, not the parameter $\lambda$ itself. This implies the periodicity
condition
$$
\operatorname D_\GR(X)_\lambda \ = \ \operatorname
D_\GR(X)_{\lambda+\mu} \ \
\ \text{if $\mu\in\fh^*$ is $H$-integral}.
\tag2.13
$$ Bernstein-Lunts \cite{BL} extend the standard operations on derived
categories to the equivariant setting, among them Verdier duality operator
$\Bbb D$. Duality reverses the twisting from $\lambda-\rho$ to
$-(\lambda-\rho) = (-\lambda+2\rho)-\rho$, so 
$$
\Bbb D : \operatorname D_\GR(X)_\lambda @>{\ \sim\ }>> \operatorname
D_\GR(X)_{-\lambda+2\rho} \ = \ \operatorname D_\GR(X)_{-\lambda}
\tag2.14
$$ because $2\rho$ is integral.

The geometric description of representations we shall use involves an
additional ingredient: $\Cal O_X(\lambda)$, the twisted  sheaf of holomorphic
functions, with twist $(\lambda-\rho)$. Concretely, this is a subsheaf of the
sheaf of holomorphic functions  on $\hat X$, the subsheaf consisting of all
germs $f$ whose restriction to the fiber of $\hat X
\to X$ is a constant multiple of $e^{\lambda-\rho}$. Here, as before, we
identify the fiber with $H$ via the action in (2.9). This definition manifestly
depends on the particular value of $\lambda-\rho$, so $\Cal O_X(\lambda)$
does not satisfy the periodicity  analogous to (2.13). When
$\lambda-\rho$ happens to be $H$-integral, the character $e^{\lambda-\rho}$
of
$H$ pulls back to the isotropy groups $B_x$, and thus determines a
$G$-equivariant holomorphic line bundle 
$$
\bold L_{\lambda-\rho} \to X\,.
\tag2.15
$$ Equivalently, $\bold L_{\lambda-\rho}$ can be described as the line bundle
associated to the principal bundle $q:\hat X \to X$ by the character 
$e^{\lambda-\rho}$ of the structure group $H$. In this case the twisted
sheaf    $\Cal O_X(\lambda)$  becomes an actual sheaf on $X$ with an action of
$G$, and as such coincides with the sheaf of holomorphic sections of the
equivariant line bundle $\bold L_{\lambda-\rho}$\,. We should remark that
these statements depend on the presence of $h^{-1}$, rather than $h$, in (2.9).

By construction, the sheaves  $\Cal O_X(\lambda)$ on $\hat X$ are
$(\lambda-\rho)$-monodromic. Hence, for $\cf\in\operatorname
D_\GR(X)_\lambda$\,, we can define the groups
$\operatorname{Ext}^*(\cf,\Cal O_X(\lambda))$ by deriving the functor
$\operatorname{Hom}$ on the category
$Sh_{X,\lambda}$, which has enough injectives. Equivalently, one can interpret
$\Cal H\italic{om}(\cf,\Cal O_X(\lambda))$ as a sheaf on $X$, and define
$\operatorname{Ext}^*(\cf,\Cal O_X(\lambda))$ as the cohomology of
$$ R\operatorname{Hom}(\cf,\Cal O_X(\lambda))\  =\ R\Gamma(X,R\Cal
H\italic{om}(\cf,\Cal O_X(\lambda)))\,. 
$$ The paper \cite{KSd} defines a natural, functorial Fr\'echet topology and a
continuous, functorial $\GR$-action on $\operatorname{Ext}^*(\cf,\Cal
O_X(\lambda))$, with the following property:
$$
\operatorname{Ext}^p(\cf,\Cal O_X(\lambda)) \text{ is an admissible
$\GR$-module of finite length},
\tag2.16
$$ for all $p \in \Bbb Z$ and $\cf\in\operatorname D_\GR(X)_\lambda$\,.  In
particular, this representation has a $\GR$-character
$\Theta(\operatorname{Ext}^p(\cf,\Cal O_X(\lambda)))$ and a $\gr$-character
$\theta(\operatorname{Ext}^p(\cf,\Cal O_X(\lambda)))$.

The correspondence between $\cf$ and the $\GR$-module (2.16) is
contravariant. We make it covariant by inserting the Verdier duality
operator\footnote{In defining the Verdier duality (2.14), we think of twisted
sheaves as objects on
$X$; alternatively, we may think of them as sheaves on $\hat X$ and apply
Verdier duality there. The two operations coincide except for a shift in
degree by the real dimension of the fiber $H$, which is even. Thus, in the
definition (2.17), the two interpretations of $\Bbb D$ have the same effect.}
(2.14). Taking alternating sums, we define the virtual characters
$$
\aligned
\Theta(\cf) \ &= \ \tsize \sum_p (-1)^p
\Theta(\operatorname{Ext}^p(\Bbb D\cf,\Cal O_X(\lambda)))
\\
\theta(\cf) \ &= \ \tsize \sum_p (-1)^p
\theta(\operatorname{Ext}^p(\Bbb D\cf,\Cal O_X(\lambda)))\,,
\endaligned
\tag2.17
$$ which depend covariantly on $\cf\in \operatorname D_\GR(X)_{-\lambda}$.
Our geometric character formulas, which will be stated in sections 3 and 4,
express $\Theta(\cf)$ and $\theta(\cf)$ in terms of $\cf$. We should remark
that the assignments $\cf \mapsto \Theta(\cf)$ and $\cf \mapsto \theta(\cf)$
descend to the K-group of $\operatorname D_\GR(X)_{-\lambda}$\,, which is
generated by \lq\lq standard sheaves" -- i.e., direct images of equivariant,
twisted local systems on orbits. Hence, to define $\Theta(\cf)$ and
$\theta(\cf)$, it is not absolutely necessary to appeal to the results of
\cite{KSd}; instead, one may appeal to the less functorial version of (2.16)
proved in \cite{SW}.

Let $\Cal R(\GR)$ denote the category of admissible representations of
$\GR$ of finite length, and $\Cal R(\GR)_\lambda$ the full subcategory of
representations with infinitesimal character $\chi_\lambda$, in
Harish-Chandra's notation. It is important to note that
$$
\aligned &\text{every  $\pi\in\Cal R(\GR)_\lambda$ is infinitesimally
equivalent to 
$\operatorname{Ext}^0(\Bbb D\cf,\Cal O_X(\lambda))$}\\ &\text{for some
$\cf\in \operatorname D_\GR(X)_{-\lambda}$\,, such that
$\operatorname{Ext}^p(\Bbb D\cf,\Cal O_X(\lambda))=0$ if $p\neq 0$.}
\endaligned
\tag2.18
$$ In particular, then,  our character formulas  apply to all admissible
representations of finite length. The assertion (2.18) follows from
\cite{KSd, (1.1f)}, combined with the Beilinson-Bernstein's localization 
functor \cite{BB1} and the Riemann-Hilbert correspondence \cite{Ka1,Me}. This
process produces a specific $\cf\in 
\operatorname D_\GR(X)_{-\lambda}$ for any $\pi\in \Cal R(\GR)_\lambda$\,:
let
$V_\pi$ be the Harish-Chandra module of $K_\Bbb R$-finite vectors in the
representation space of $\pi$, and $\Cal M$ the (derived) Beilinson-Bernstein
localization of $V_\pi$ at $\lambda$\,; we apply the deRham functor to
$\Cal M$ (this has the effect of switching the sign of twisting parameter) and
then the equivalence of categories $\operatorname D_{K_\Bbb
C}(X)_{-\lambda}
\cong \operatorname D_\GR(X)_{-\lambda}$\,\cite{MUV}; this, up to a shift in
degree, gives the sheaf $\cf$.

The K-group of $\Cal R(\GR)_\lambda$ is generated by standard
representations, i.e., representations parabolically induced from discrete
series representations. A more explicit description of the correspondence
(2.18), in the case of standard representations, is crucial both for
applications of our formulas and for their proof. We give such an explicit
description in the sections below.
\vskip 1cm

\subheading{\bf 3. The  integral formula}
\vskip .5cm

Let $\GR$ be a unimodular Lie group of type I, $\pi$ an irreducible unitary
representation of $\GR$, and $\tp$ the character of $\pi$ on the Lie algebra. 
Kirillov's \lq\lq universal character formula" \cite{Ki3} attempts to express
the Fourier transform $\hat\theta_\pi$ as an integral over a coadjoint orbit.
Beyond the cases of compact and nilpotent groups, which were established by
Kirillov himself \cite{Ki1,Ki2}, the formula has been proved, in full generality,
for type I solvable groups  \cite{C}. The reductive case is more subtle:
Rossmann has shown that the formula applies in the case of tempered
representations
\cite {R1}, but the formula fails for non-tempered representations -- there
are simply not enough coadjoint orbits. Duflo's construction \cite{D2}
produces irreducible unitary  representations of algebraic groups by induction
from representations of reductive groups. When this is combined with Rossmann's
formula, it implies the validity of the \lq\lq universal character formula" for
{\it generic} irreducible unitary representations of type I algebraic groups.

Rossmann \cite{R3} has proposed a remedy for the failure of the  \lq\lq
universal  formula" in the reductive case: the Fourier transform $\hat\theta$
of any invariant eigendistribution $\theta$ on $\gr$, with regular
infinitesimal character, can be expressed as an integral over an appropriate
cycle in a coadjoint orbit for the complexification $G$. Our result is an
explicit Rossmann type formula for the character of any admissible
representation of finite length, unitary or not, with regular or singular
infinitesimal character.

Let us recall the definition of Rossmann's twisted moment map \cite{R3}. For
details we refer the reader to \cite{SV4, section 8}, which was written with
the present application in mind. We fix a compact real form $\UR \subset G$
which contains $\KR$. The twisted moment map corresponding to any particular
$\lambda\in\fh^*$\,,
$$
\mu_\lambda\ : \ \ct \longrightarrow \fg^*\,,
\tag3.1a
$$
is defined as follows. Each $x\in X$ is fixed by a unique maximal torus $T_\Bbb
R\subset \UR$\,; its complexified Lie algebra $\ft$ becomes canonically
isomorphic to the universal Cartan algebra $\fh$ via $\ft \cong
\fb_x/[\fb_x,\fb_x] \cong \fh$. This identification makes $\lambda$ correspond
to a $\lambda_x \in \ft^*$, which we extend to a linear functional on $\fg$
by means of the canonical splitting $\fg = \ft \oplus [\ft,\fg]$\,. The map $x
\mapsto \lambda_x \in \fg^*$ is real algebraic and $\UR$-equivariant, and the
ordinary moment map $\mu:\ct \to \fg^*$ is even complex algebraic and
$G$-equivariant, so 
$$
\mu_\lambda(x,\xi) \ = \ \lambda_x \, + \, \mu(x,\xi) \qquad (\,\xi \in
T^*_xX\,)
\tag3.1b
$$
is a $\UR$-equivariant, real algebraic map. The $G$-orbit through $\lambda_x$
depends only on $\lambda$, not on $x$, so we denote it by $\Omega_\lambda$.
The correspondence $\lambda \mapsto \Omega_\lambda$, we note, induces a
bijection between $W\backslash \fh^*$ and the set of semisimple coadjoint
orbits for $G$. When $\lambda$ is regular, $\mu_\lambda$ takes values in
$\Omega_\lambda$, and 
$$
\mu_\lambda \ : \ \ct \, @>{\ \sim \ }>> \, \Omega_\lambda
\tag3.1c
$$
is a real algebraic isomorphism. At the opposite extreme, for $\lambda  = 0$\,,
the twisted moment map reduces to the ordinary moment map, of course.

Both the cotangent bundle $\ct$ and the  coadjoint orbit $\omega_\lambda$ come
equipped with canonical, complex algebraic, $G$-invariant symplectic forms,
$\sigma$ and $\sigma_\lambda$\,, respectively. We define a $\UR$-invariant, real
algebraic two form
$\tau_\lambda$ on $X$ by the formula
$$
\tau_\lambda(u_x,v_x) = \lambda_x([u,v])\,;
\tag3.2
$$  
here $u_x,v_x \in T_xX$ are the tangent vectors at $x$ induced by $u,v \in \frak
u_\Bbb R$ via differentiation of the $\UR$-action. When
$\lambda \in \frak h^*$ happens to be integral, $-\tau_\lambda$ is  the
curvature form of the (essentially unique) $ U_\Bbb R$-invariant metric on the
$G$-invariant algebraic line bundle $\bold L_\lambda \to X$ parametrized by
$\lambda \in \frak h^*$ [GS]. Moreover, in this situation, ${(2\pi i)}^{-1}
{\tau_\lambda}$ represents the Chern class of $\bold L_\lambda$\,. 
Whenever $\lambda$ is regular, the three differential forms $\sigma$,
$\sigma_\lambda$, $\tau_\lambda$ are related. Let $\pi:\ct \to X$ denote the
natural projection; then:

\proclaim{3.3 Proposition} For $\lambda$ regular,
$
\mu_\lambda^*\sigma_\lambda  =   -\sigma + \pi^*\tau_\lambda\,.
$
\endproclaim

\demo{Proof} The inverse image $\mu^{-1}(\Omega)$ of the regular nilpotent
coadjoint orbit $\Omega$ is dense in $\ct$, so it suffices to verify the
identity on $\mu^{-1}(\Omega)$. Let $\sigma_\Omega$ be the canonical
$G$-invariant, complex algebraic symplectic form on $\Omega$. Rossmann
\cite{R3, 7.2, p.172} proves the identity
$$
\mu_\lambda^*\sigma_\lambda  = \mu^*\sigma_\Omega + 
\pi^*\tau_\lambda \ \ \ \ \text{on  $\mu^{-1}(\Omega)$}\,;
$$  
on the other hand, by \cite{SV4, 8.19},
$$
\mu^*\sigma_\Omega \ = \ - \sigma \ \ \ \ \text{on  $\mu^{-1}(\Omega)$}\,.
$$
We should remark that both  identities hold for any nilpotent coadjoint orbit
$\Omega$, at the smooth points of  $\mu^{-1}(\Omega)$. The proofs of this more
general version of the two identities in \cite{R3,SV4} simplify considerably in
our more special situation. The two identities together imply the proposition.
\enddemo

The $\GR$-action on the flag variety $X$ is real algebraic and has finitely many
orbits. The resulting orbit stratification of  $X$ is semi-algebraic and
satisfies the Whitney conditions. Hence
$$
T^*_\GR X \ = _{\text{def}} \ \text{union of the conormal bundles of the
$\GR$-orbits}
\tag3.4
$$
is a closed, semi-algebraic subset of $\ct$. Like any finite union of conormal
bundles, $T^*_\GR X$ is Lagrangian (at its smooth point) with respect to the
canonical symplectic structure on $\ct$, viewed as the $C^\infty$ (not
holomorphic!) cotangent bundle of $X$. Equivalently this is the symplectic
structure defined by the real two form\footnote{There are at least two
different but equally natural identifications between the holomorphic and the
$C^\infty$ cotangent bundles; we use the convention of \cite{KSa, Chapter XI}.}
$2\operatorname{Re}\sigma$. In particular, if $n$ denotes the complex
dimension of $X$, $T^*_\GR X$ has dimension $2n$. We write $\ohf *(T^*_\GR X,
\Bbb Z)$ for the homology of $T^*_\GR X$ with infinite support -- i.e., for the
Borel-Moore homology. Then $\ohf {2n}(T^*_\GR X,\Bbb Z)$ is the group of top
dimensional, possibly infinite cycles on $T^*_\GR X$.

Kashiwara's characteristic cycle construction \cite{Ka2,KSa,SV4} defines a
$\Bbb Z$-linear map
$$
\occ \ : \ K(\operatorname D_\GR(X)_{\lambda}) \longrightarrow  \ohf
{2n}(T^*_\GR X,\Bbb Z)\,,
\tag3.5
$$
from the $K$-group of the twisted equivariant derived category 
$\operatorname D_\GR(X)_{\lambda}$. Strictly speaking, the cited references
deal only with the untwisted case. However, the characteristic cycle
construction is local with respect to the base manifold $X$, and locally twisted
sheaves  \lq\lq can be untwisted" -- they can identified with sheaves -- so
$\occ$ makes sense also for twisted sheaves. Alternatively but equivalently, we
may view any object 
$\cf \in \operatorname D_\GR(X)_{\lambda}$ as a  monodromic sheaf on $\hat X$
and take its characteristic cycle $\occ_{\hat X}(\cf)$ in $T^*\hat X$; because of
the monodromicity condition, $\occ_{\hat X}(\cf)$ is the pullback of a Lagrangian
cycle in $\ct$, namely $\occ(\cf)$. The map $\occ$ becomes intrinsic only when
one fixes an orientation of the base manifold $X$. We use the complex structure
on $X$ to do so. 

Like Rossmann \cite{R3}, we define the Fourier transform $\hat\phi$ of a test
function $\phi$ in $C^\infty_c(\gr)$ without choosing a square root of -1\,, as
a holomorphic function on $\fg^*$\,,
$$
\hat \phi (\zeta) \ = \ \int_\gr e^{\zeta(x)}\phi(x) dx \qquad
(\zeta\in\fg^*)\,.
\tag3.6
$$
Then $\hat \phi$ decays rapidly in the imaginary directions. The particular
normalization of the Euclidean measure $dx$ on $\gr$ will not matter in the end.
One may choose to think of the Fourier transform more invariantly as attached to
a smooth compactly supported measure; we leave it to the reader to reinterpret
our character formula in these terms.

\proclaim{3.7 Proposition} For all $\phi \in C^\infty_c(\gr)$, $C \in \ohf
{2n}(T^*_\GR X,\Bbb Z)$, and $\lambda \in \fh^*$, the integral
$$
\int_C \mu^*_\lambda \hat \phi \ (-\sigma + \pi^* \tau_\lambda)^n
$$
converges absolutely. The value of the integral depends holomorphically on
$\lambda$\,.
\endproclaim

Rossmann \cite{R3}, in the case of a regular $\lambda \in \fh^*$, states and
proves the convergence and the holomorphic behavior of the corresponding 
(via 3.3) integral on
$\mu_\lambda(C) \subset \Omega_\lambda$.  Whether or not $\lambda$ is regular,
this comes down to the rapid decay of the holomorphic function $\hat \phi$ in the
imaginary directions, as  will become clear at the end of this section where
we establish a more general convergence criterion.

Recall the definition (2.17) of the virtual character $\theta(\cf)$. It is a
generalized function, and as such can be integrated against any smooth, compactly
supported test function $\phi$.

\proclaim{3.8 Theorem}  Let $\lambda \in \fh^*$ and $\cf \in \operatorname
D_\GR(X)_{-\lambda}$ be given.  Then
$$
\int_{\gr}\theta(\cf)\,\phi\,dx \ = \ \frac 1 {(2\pi i)^n n!}
\int_{\occ(\cf)} \mu^*_\lambda \hat \phi \ (-\sigma + \pi^* \tau_\lambda)^n 
\qquad (\,\phi \in C^{\infty}_c(\gr)\,)\,.
$$
\endproclaim

Appearances notwithstanding, the formula has canonical meaning, independent of
the particular choice of $i=\sqrt{-1}$. We mentioned already that the
orientation of $\occ(\cf)$ depends on the choice of an orientation of $X$. When
we orient the complex manifold $X$, we use the complex structure operator $J$,
which reflects the choice of $i=\sqrt{-1}$. This same choice affects also the
sign on the right hand side of our formula, as it must.

In effect, the theorem provides several integral formulas for any virtual
character, corresponding to the various choices of the localization parameter
$\lambda$ in any particular $W$-orbit. When 
$\lambda$ is regular, we can use the isomorphism (3.1a) and proposition 3.3 to
rewrite our formula on the complex, regular, elliptic coadjoint orbit 
$\Omega_\lambda$\,, 
$$
\int_{\gr}\theta(\cf)\,\phi\,dx \ = \ \frac 1 {(2\pi i)^n n!}
\int_{\mu_\lambda(\occ(\cf))}  \hat \phi \ \sigma_\lambda ^n 
\,.
\tag3.9
$$
Rossmann \cite{R3} shows that any invariant eigendistribution $\theta$ with
regular infinitesimal character can be written as  an integral of this type, over
an unspecified cycle $C\in \ohf {2n}(T^*_\GR X,\Bbb C)$. He identifies the
cycle explicitly under the following circumstances: $\GR$ is a complex Lie
group, and $\lambda$ is integral, regular, and dominant (with Rossmann's
conventions, anti-dominant). 

The integrand on the right in (3.9) is a holomorphic $2n$-form, hence is
closed, which gets integrated over an (infinite) $2n$-cycle in the
$4n$-manifold $\Omega_\lambda$\,. As Rossmann \cite{R3} points out, the
integral remains unchanged when one replaces  the cycle of integration $C_0$ by
another cycle  $C_1$, homologous to $C_0$ under an appropriately restricted
notion of homology which takes into account the growth of cycles and chains at
infinity. This observation plays a crucial role in our proof of theorem 3.8. When
$\theta(\cf)$ is a tempered irreducible character, with $\lambda$ regular and
appropriately chosen within its $W$-orbit,  the
cycle $\occ(\cf)$ turns out to be homologous in the restricted sense to the
$\mu_{\lambda}^{-1}$-image of a
$\GR$-orbit in $i\frak g^*_\Bbb R$. In this particular case, then, the integral
formula (3.9) is equivalent to the \lq\lq universal formula" as proved by
Rossmann
\cite{R1} -- our justification for viewing (3.8) as an appropriate
analogue of Kirillov's formula. 

We need to be more precise about the restricted notion of homology. Initially,
we consider a real algebraic manifold $M$ of dimension $m$ and a $C^\infty$
differential form $\omega$ on $M$, of degree $d$; $\omega$ maybe real or
complex valued. Recall the notion of a semi-algebraic, locally finite $d$-chain
$C$
\cite{SV4}.  We want to define the
$\omega$-norm
$\|C\| = \|C\|_\omega$ in case $|C|=\,$ support of $C$ is  compact  --
which makes
$C$ finite. Thus we can express
$C$ as a finite integral linear combination 
$$
C\ = \ \sum_{j=1}^N \, n_j \,S_j
\tag3.10
$$
of pairwise disjoint, connected, smooth, oriented, relatively compact,
$d$-dimensional semi-algebraic sets $S_j$. These hypotheses ensure that 
$$
\|C\| \ = \ \sum_{j=1}^N \, |n_j| \,\int_{S_j}\,|\omega|
\tag3.11
$$
is well defined and independent of the particular choice of the expression
(3.10). In this defining formula, the absolute value $|\omega|$ of the form
$\omega$ is viewed as a measure on each of the $S_j$. 

Let us apply the preceding discussion to the case of $M=\ct$, with $\omega =
(-\sigma + \pi^* \tau_\lambda)^n$ and $d=2n$. We equip $X$ with a $U_\Bbb
R$-invariant hermitian metric. Any such metric is real algebraic, since $U_\Bbb
R$ acts transitively on $X$. We define 
$$
D(r) \ = \ \{\,(x,\xi)\in\ct \ | \ \| \xi\| \leq r\,\}\,,
\tag3.12
$$
the disc bundle of radius $r$ in $\ct$; here $\| \xi\|$, the norm of a cotangent
vector $\xi$\,, is measured with respect to the Hermitian metric. The
particular choice of metric will not matter, since any two are mutually
bounded. For convenience, we choose it so that 
$$
\mu \ : \ T^*_xX \ \longrightarrow \ \fg^* \ \ \  \text{is an isometry for each
$x\in X$}
\tag3.13
$$
with respect to some $U_\Bbb R$-invariant inner product on $\fg^*$. We consider
a locally finite, semi-algebraic $2n$-chain $C$. Since $D(r)$ is compact
semi-algebraic, the intersection $C\cap D(r)$ can be viewed as a $2n$-chain
\footnote{Stratify
 $|C|=$ support of $C$ compatibly with its intersection with $D(r)$; then
each
$2n$-stratum in $|C|\cap D(r)$ inherits a multiplicity from $C$.} with
compact support. Then $C$  has polynomial growth with respect to $\omega =
(-\sigma + \pi^*
\tau_\lambda)^n$, in the sense that
$$
\|C\cap D(r)\| \ \ \ \text{grows at most polynomially with $r$}\,.
\tag3.14
$$
To see this, we compactify $\ct$ by regarding it as a real algebraic $\Bbb R$
-vector bundle with distinguished Riemannian metric, hence with structure group
$O(2n, \Bbb R)$. The linear action of $O(2n, \Bbb R)$ on $\Bbb R^{2n}$ extends to
an algebraic action on the one point compactification $S^{2n}$. Thus, replacing
the typical fiber $\Bbb R^{2n}$ by $S^{2n}$, we obtain a real algebraic
compactification $\overline {\ct}$ of $\ct$. The function $r^{-1}$ extends real
algebraically the the complement of the 0-section in  $\overline {\ct}$, and the
extended function vanishes to first order on the locus at infinity. Since
$\omega$ is real algebraic on $\ct$, $r^{-k}\omega$  extends real algebraically
across the locus at infinity for $k$ sufficiently large -- in fact, for $k\geq
n$, though the precise value will not matter to us. The chain $C$, being
semi-algebraic, extends to a finite chain in  $\overline {\ct}$, and this now
implies (3.14).  

Now let $\phi\in C^\infty_c(\gr)$ be a test function. The Fourier transform
$\hat\phi$ as defined in (3.6) is holomorphic on $\fg^*$ and decays rapidly in
imaginary directions; the rate of decay can be uniformly bounded as long as
the real part of the argument is restricted to a compact set (Paley-Wiener
theorem). It is clear from the definition of the twisted moment map that
$\mu_\lambda$ differs from the ordinary moment map $\mu$ by a  term  whose norm
can be bounded by a constant multiple of the norm of $\lambda$\,. 
Hence, for every $N\in \Bbb N\,$ and $R>0\,$,  there exists a constant
$A=A(\phi,N,R)$, such that
$$
|\mu_\lambda^*\hat\phi(x,\xi)| \ \leq \ \frac A {1+\|\xi\|^N}\ \ \ \ \text{if
$\ \ \|\operatorname{Re}\mu(x,\xi)\|\,+\,\|\lambda\|\,\leq \,R\,$}.
\tag3.15
$$
The values of the differential form $\mu_\lambda^*\hat\phi\, (-\sigma +
\pi^* \tau_\lambda)^n$ at the various points in $\ct$ depend holomorphically
on $\lambda$, and the polynomial growth condition (3.14), for any particular
semi-algebraic $2n$-chain $C$, holds locally uniformly in $\lambda$. We conclude:

\proclaim{3.16 Lemma} If the
real part of $\mu(x,\xi)$ is bounded on the support $|C|$, the integral
$$
\int_C \, \mu_\lambda^*\hat\phi\, (-\sigma + \pi^* \tau_\lambda)^n
$$
converges absolutely. The value of this integral depends holomorphically on the
parameter $\lambda$.
\endproclaim

The moment map $\mu$ takes purely imaginary values on any chain supported on
$T^*_\GR X$ -- in fact, 
$$
T^*_\GR X \ = \ \mu^{-1}(i\fg^*_\Bbb R)\,,
\tag3.17
$$
as is easy to see; e.g., in \cite{R3}. Thus 3.16 implies proposition 3.7.

We now consider two semi-algebraic $2n$-cycles $C_1,C_2$ in
$\ct$, both of which satisfy the boundedness hypothesis of lemma 3.16. We suppose
that
$C_1$ and $C_2$ are homologous, in the sense that 
$$
C_1 - C_2 \ = \ \partial \tilde C
\tag3.18a
$$
for some semi-algebraic $(2n+1)$-chain $\tilde C$. We require that
$$
\gathered
\text{the image of the support  $|\tilde C|$ }
\\
\text{under the map
$(x,\xi)\mapsto
\operatorname{Re}\mu(x,\xi)$ is bounded}\,.
\endgathered
\tag3.18b
$$
\proclaim{3.19 Lemma} Under the hypotheses (3.18),
$$
\int_{C_1} \mu_\lambda^*\hat\phi\, (-\sigma + \pi^*
\tau_\lambda)^n \ = \ \int_{C_2} \mu_\lambda^*\hat\phi\, (-\sigma + \pi^*
\tau_\lambda)^n\,.
$$
\endproclaim
\demo{Proof}
Let us argue first of all that the $2n$-form $\mu_\lambda^*\hat\phi\, (-\sigma + \pi^*
\tau_\lambda)^n$ on $\ct$ is closed. It suffices to show this when $\lambda$ is regular, since the form
depends real analytically on $\lambda$. In the regular case, 
$$
\mu_\lambda^*\hat\phi\, (-\sigma + \pi^* \tau_\lambda)^n \ = \ \mu_\lambda^*
(\hat\phi\,\sigma_\lambda^n)
$$
is the pullback to $\ct$ of a holomorphic form of top degree,
$\hat\phi\,\sigma_\lambda^n$, on the complex manifold $\Omega_\lambda$. As
such, it is closed. We restrict the homology relation (3.18a) to the disc bundle
$D(r)$: 
$$
C_1\cap D(r) - C_2\cap D(r) \ = \ \partial(\tilde C\cap D(r)) -C_3(r),
\tag3.20
$$
with $C_3(r)$ supported on the boundary $\partial D(r)$. Thus Stokes' theorem
implies
$$
\gathered
\int_{C_1\cap D(r)} \mu_\lambda^*\hat\phi\, (-\sigma + \pi^*
\tau_\lambda)^n \ -\ \int_{C_2\cap D(r)} \mu_\lambda^*\hat\phi\, (-\sigma +
\pi^*\tau_\lambda)^n \ = 
\\
-\int_{C_3(r)} \mu_\lambda^*\hat\phi\, (-\sigma + \pi^*
\tau_\lambda)^n \,,
\endgathered
\tag3.21
$$
so  the lemma will follow if the term on the right can be made arbitrarily small
as $r$ tends to infinity. Arguing as in the proof of lemma 3.16, we use (3.15)
and (3.18b) to remove $ \mu_\lambda^*\hat\phi$ from the picture: we only need to
show that the $\omega$-norm 
$$
\|\tilde C\cap \partial D(r)\| \ \ \ \text{grows at most polynomially with
$r$ as $r\to\infty$}\,.
\tag3.22
$$
In interpreting this statement, we need to eliminate the -- finitely many
-- values of
$r$ for which the intersection $\tilde C\cap \partial D(r)$ fails to be
transverse. Note that the complement of the 0-section in $\overline{\ct}$ is
naturally  isomorphic to $\ct$; the isomorphism transforms the semi-algebraic
chain $\tilde C$ into another semi-algebraic chain, the function $r$ to
$r^{-1}$, and the real algebraic form $\omega$ to a form with a pole of finite
order  along the (new) zero-section. We remove the pole by multiplying the form
with an appropriate power of the (new) function $r$.  Thus, changing notation, we
only need to show that 
$$
\|\tilde C\cap \partial D(r)\| \ \ \ \text{remains bounded as $r\to 0$}
\tag3.23
$$
when we measure the size of the chain $\tilde C\cap \partial D(r)$ with respect
to a form which is algebraic  around the 0-section. 

To establish (3.23), one can argue as follows. The problem is local around the
0-section, which is compact. The boundedness with respect to some algebraic form
is therefore implied by the boundedness of the volume of $\tilde C\cap \partial
D(r)$ near $r=0$ when the volume is measured with respect to any particular
Riemannian metric on the ambient manifold $\ct$. Since we only need to establish
boundedness of the volume, we can make several simplifying assumptions. First of
all, the assumption that $\tilde C$ consists of a single  component of
multiplicity one. Secondly, we can uniformize $\tilde C$ \cite{Hn1,2;\, BM,
theorem 0.1} by a  real analytic map $f: N \to \ct$, with $N$ a compact real
analytic manifold of the same dimension as $\tilde C$, and $f(N)=\tilde C$. The
pullback to $N$ of the metric on $\ct$ is bounded by a metric on $N$, and
$f^{-1}(\partial D(r))$ covers $\tilde C\cap \partial D(r)$ generically
finitely. This leaves us with the following problem: a compact, connected, real
analytic manifold $N$, a non-constant real analytic function $r$, and we must
show that the level sets $r=c$ have volume depending continuously on $c$. At
non-critical values of $r$, that is obvious. To deal with critical values, one
can put $r$ into normal form by a further step of uniformization \cite{Hn1,2;\,
BM, corollary 4.9}. 

Lemmas 3.16 and 3.19 give us the right notion of restricted homology: we
consider semi-algebraic cycles, on whose support $\operatorname{Re}\mu$ remains
bounded, modulo boundaries of chains satisfying the same conditions. According
to the two lemmas, the integral of  $ \mu_\lambda^*\hat\phi\, (-\sigma +
\pi^*
\tau_\lambda)^n$ over a restricted homology class becomes a well defined quantity.
\enddemo
 \vskip 1cm

\subheading{\bf 4. Character cycles}
\vskip .5cm

Invariant eigendistributions, either on $\GR$ or on $\gr$, are locally $L^1$
functions, real analytic on, respectively, $\GR'$ and $\gr'$, the sets of regular
semisimple elements in $\GR$ and $\gr$ -- that is the assertion of
Harish-Chandra's regularity theorem \cite{HC4}. The complements of $\GR'$ and
$\gr'$ have measure zero, so an invariant eigendistribution is completely
determined by its  restriction to the set of regular semisimple elements. 

For definiteness, let us talk about the group case; at various points we shall
indicate how statements need to be modified for the case of the Lie algebra. We
choose a set of representatives $\{T_{1,\Bbb R}, \dots ,T_{m,\Bbb R}\}$ of the
finitely many conjugacy classes of Cartan subgroups. Then the conjugates of the
$\GR'\cap T_{i,\Bbb R}$, $1\leq i \leq m$, cover $\GR'$, so any invariant
eigendistribution $\Theta$ is completely determined by its restriction to the
various $\GR'\cap T_{i,\Bbb R}$. To simplify the notation, we pick out one of
the $ T_{i,\Bbb R}$, drop the subscript $i$, and set $T_\Bbb R '= \GR'\cap
T_{\Bbb R}$. The invariant eigendistribution $\Theta$ corresponds to an
infinitesimal character $\chi_\lambda: \Cal Z(\fg) \to \Bbb C$\,, with
$\lambda\in \fh^*$. The parameter $\lambda$ is only determined up to Weyl
conjugation, but we pick a specific representative. As we have done before, we
index the various identifications between the universal Cartan $\fh$ and the
complexified Lie algebra $\ft$ of $T_\Bbb R$ by the fixed points of $T_\Bbb R$
on the flag manifold $X$:
$$
\ft \ \cong \ \fb_x/[\fb_x,\fb_x] \ \cong \ \fh \qquad (\,x\in X^{\TR} \,=
\,\text{fixed point set of $\TR$\,})\,.
\tag4.1
$$
These isomorphisms induce isomorphisms $\fh^* \cong \ft^*$, $\mu\mapsto
\mu_x$\,. The universal Weyl group $W$ acts canonically and simply transitively on
the fixed point set of $T_\Bbb R$. Initially we choose a particular fixed point
$x$. Then, for any given $g\in T_\Bbb R'$, there exist polynomial functions
$p_{g,w}$ on $\ft$, indexed by $g$ and by $w\in W/W_\lambda$\,, the quotient of
$W$ by the isotropy subgroup at $\lambda$, such that 
$$
\Theta(g \exp \zeta) \ = \ \frac {\tsize\sum_{w\in W/W_\lambda} p_{g,w}(\zeta)
e^{(w\lambda)_x(\zeta) - \rho_x(\zeta)}} {\tsize \prod_{\alpha \in
{\Phi}^{ +}} (1-e^{-\alpha_x})(g \exp(\zeta))}
\tag4.2
$$
for all small $\zeta\in\ft_\Bbb R$. Here $\Phi^+$ and $\rho$ have the same
meaning as in (2.10), and $e^{-\alpha_x}$ denotes the character of $T_\Bbb R$
corresponding to the root $-\alpha_x$, as is customary. These local
expressions, even for only a single $g$ in every connected component of
$\GR'\cap T_{i,\Bbb R}$ and every $i$, completely determine $\Theta$.
Conversely $\Theta$ determines the coefficient polynomials $p_{g,w}$, which have
degree not exceeding  the order of $W_\lambda$ -- in particular, the $p_{g,w}$
are constant whenever $\lambda$ is regular. 

The preceding statements constitute the easy part of Harish-Chandra's proof of
the regularity theorem. One may ask whether a collection of local expressions
as in (4.2) does define an invariant eigendistribution. For this it is
certainly necessary that the local expressions fit together and give global,
well defined functions on the various $\GR'\cap T_{i,\Bbb R}$\,. Beyond this
obvious necessary condition, there are necessary and sufficient conditions due
to Harish-Chandra \cite{HC4}, usually referred to as the Harish-Chandra
matching conditions. These involve the continuations of the numerators of
the local expressions (4.2) to the subregular semisimple points; such points lie
in (conjugates of) at most two of the $T_{i,\Bbb R}$, and the matching conditions  relate
the numerators on those Cartans.

The coefficient polynomials are constant not only when $\lambda$ is regular, as
we had remarked earlier, but also when $\Theta$ is a virtual character
corresponding to a possibly singular infinitesimal character; this observation
was made by Fomin-Shapovalov \cite{FS}. Conversely, an invariant
eigendistribution $\Theta$ whose local expressions involve only constant
coefficients is a $\Bbb C$-linear combination of characters. From now on we
shall restrict attention to this particular situation: we assume that the
$p_{g,w}$ are constants, as will be the case for all virtual characters.

Let us rewrite the local expression (4.2) in slightly more invariant terms. As
before, $g\in T_\Bbb R$ will be regular, and $\zeta\in\ft_\Bbb R$ sufficiently
small.  We note that  $(w\lambda)_{wx} = \lambda_x$\,, hence
$$
\Theta(g \exp \zeta) \ = \ \sum_{x\in X^{T_\Bbb R}}\ \frac { c_{g,x}\,
e^{\lambda_x(\zeta) - \rho_x(\zeta)}} {\tsize \prod_{\alpha \in {\Phi}^{ +}}
(1-e^{-\alpha_x})(g \exp(\zeta))} 
\tag4.3
$$
summed over the  fixed point set $ X^{T_\Bbb R}$ of $T_\Bbb R$. Because of
our assumption, the coefficients $c_{g,x}$ are constants. They are parametrized
by $x\in  X^{T_\Bbb R}$, rather than $w \in W/W_\lambda$ as in (4.2), so there
will be repetitions in the sum when $\lambda$ is singular.  Now, to make the
$c_{g,x}$ unique, we require that $c_{g,vx} = \operatorname{sgn}(v)
c_{g,x}$ for $v\in W_\lambda$ -- the sign accounts for the fact that the fixed
point $x$ has become variable whereas we had chosen a particular fixed point in
(4.2).

The formula (4.3) has an analogue on the Lie algebra, which follows from
(4.3) and the universal expression 
$$
\exp_*|_\zeta\  = \ \ell(\exp(\zeta))_* \circ\frac
{1-e^{-\operatorname{ad}\zeta}}{\operatorname{ad}\zeta}
$$
for the differential of the exponential map. The analogue on the Lie algebra has
a slightly simpler appearance, because roots and weights are globally well
defined functions on the Cartan subalgebra
$\ft_\Bbb R$\,: there exist constants
$d_{E,x}$, indexed by the various connected components $E\subset \ft_\Bbb R'=
\gr'\cap \ft_\Bbb R$ and $x \in X^{T_{\Bbb R}}$, such that
$$
\theta(\zeta) \ = \ \sum_{x\in X^{T_\Bbb R}}\ \frac { d_{E,x}\,
e^{\lambda_x }(\zeta)} {\tsize \prod_{\alpha \in {\Phi}^{ +}}
\, \alpha_x(\zeta)} \qquad \text{for all $\zeta \in E \subset \ft_{\Bbb
R}'$\,}\,;
\tag4.4
$$
here $\theta$ and $\Theta$ are related by the formula (2.1). As in (4.3), the
coefficients become uniquely determined when one requires $d_{E,vx}=
\operatorname{sgn}(v) d_{E,x}$ for all $v\in W_\lambda$. 

We return to the group case (4.3). Hotta-Kashiwara \cite{HK,Ka5} have given 
geometric meaning to the coefficients $c_{g,x}$\,. Let $\tilde G$ denote the
variety of pairs $(g,x) \in G\times X$, such that $x$ is a fixed point of $g$.
When one makes $G$ act on itself by conjugation and on $X$ by translation, the
two projections
$$
\pi \ :\  \tilde G \longrightarrow X \,, \qquad q \ :\  \tilde G
\longrightarrow G\,.
\tag4.5
$$
become $G$-equivariant.  The former exhibits $\tilde G$ as the tautological
bundle of Borel subgroups; in particular, $\tilde G$ is smooth. On the other
hand, $q$ induces  a (Galois) covering map over the regular set,
$$
 \tilde G'\ = \ q^{-1}(G')\  @>{ \ \ q \ \ }>>\  G' \,, \qquad \text{with
covering group
$W$}\,,
\tag4.6
$$
since $W$ acts simply transitively on the fixed point set $X^g$ of any $g\in
G'$.  Note that
$$
\tilde \GR \ = \ \{\,(g,x) \in \tilde G \ | \ g \in \GR \,\} \ = \ q^{-1}\GR
\tag4.7
$$
is a -- usually singular -- real algebraic subvariety of $\tilde G$.  Recall
(2.4): since $\pi^{-1}(x) = B_x$ for $x \in X$, there exists a natural map
$$
\tilde G \longrightarrow H
\tag4.8
$$
from $\tilde G$ to the universal Cartan group $H$, which becomes
$G$-equivariant when one thinks of $G$ as acting trivially on $H$.  For
$\lambda \in \fh$\,, the multiple valued function $e^{\lambda - \rho}$ generates
a rank one local system on $H$, which entered the definition of the twisted
$\GR$-equivariant derived category $\operatorname{D}_{\GR}(X)_\lambda$ in
section 2.  The pullback, via (4.8),  of this local system from $H$ to
$\tilde G$ defines a $G$-equivariant rank one local system $\Bbb C_\lambda$ and
a $G$-equivariant inclusion
$$
\Bbb C_\lambda \, \hookrightarrow \, \Cal O_{\tilde G}
\tag4.9
$$
into the sheaf of holomorphic functions.  The local systems $\Bbb C_\lambda$\,,
$\lambda \in \fh$, satisfy the periodicity condition analogous to (2.13): there
exists a natural isomorphism
$$
\Bbb C_\lambda \ \cong \ \Bbb C_{\lambda + \mu} \qquad \text{if $\mu \in \fh^*$
is $H$-integral}\,,
\tag4.10
$$
given by multiplication with the globally defined holomorphic function
$e^\mu$ on $\tilde G$.   As in (2.12),
$\lambda =
\rho$ corresponds to the \lq\lq untwisted" case $\Bbb C_\rho = \Bbb C_{\tilde
G}$\,.

Appearances to the contrary, the local formula (4.3) has global meaning: each
regular $g\in\GR'$ lies in a unique Cartan subgroup $T_\Bbb R \subset \GR$\,,
and $\Theta$ is conjugation invariant and real analytic on $\GR'$. It follows
that the coefficients $c_{g,x}$, which  are indexed by pairs $(g,x)\in \tilde\GR'
= q^{-1}(\GR')$, are $\GR$-invariant, real analytic   functions on
$\tilde\GR'$. Setting $\zeta=0$ in (4.3) gives the global expression 
$$
\Theta(g)  \ = \ \sum_{x\in X^{g}}\ \frac { c_{g,x}} {\tsize \prod_{\alpha \in {\Phi}^{ +}} (1-e^{-\alpha_x})(g
)} \qquad (\,g\in \GR'\,)\,.
\tag4.11
$$
Again by (4.3), $c_{g,x}e^{(\lambda-\rho)_x(\zeta)} = c_{g\exp\zeta,x}$ whenever
$\zeta\in \gr$ is small and centralizes $g$. It follows that
$$
c\,:\, (g,x) \mapsto c_{g,x} \  \ \text{is a section of $\,\Bbb C_\lambda \subset
\Cal O_{\tilde G}\,$ over $\tilde\GR'$}\,.
\tag4.12
$$
We fix an orientation of $\GR$\,; this induces an orientation also on
$\tilde\GR'$ via the covering map (4.6). Then, by Poincar\'e duality, 
$$
c\in\oh^0(\tilde\GR',\Bbb C_\lambda) = \ohf
d(\tilde\GR',\Bbb C_{-\lambda+2\rho}) \cong \ohf d(\tilde\GR',\Bbb
C_{-\lambda})\,;
\tag4.13
$$
here $d = \operatorname{dim}_\Bbb R \GR$ is the dimension of $\tilde \GR$ and
$\ohf *(\  \cdot \ , \Bbb C_{-\lambda})$ denotes homology with locally finite
support and values in the local system $\Bbb C_{-\lambda}$\,. The passage from
$\Bbb C_\lambda$ to the dual local system $\Bbb C_{-\lambda+2\rho}$ is forced
by the formalism of Poincar\'e duality, since cohomology behaves covariantly with
respect to the coefficient system, whereas homology behaves
contravariantly\footnote{more precisely, one can take cohomology with values
in a sheaf, but homology with values in a cosheaf; local systems can be
regarded as either, but the identification between the two is contravariant.}. 

Homology with locally finite supports can be restricted to open
subsets. Since $d$ is the dimension of $\tilde\GR$\,, 
$$
\ohf d(\tilde\GR,\Bbb C_{-\lambda}) \  \hookrightarrow \ \ohf
d(\tilde\GR',\Bbb C_{-\lambda})
\tag4.14
$$
is an inclusion.  Hotta-Kashiwara \cite{HK,Ka5} show that every $c=c(\Theta)$,
coming from an invariant eigendistribution $\Theta$ via (4.11-14), lies in $\ohf
d(\tilde\GR,\Bbb C_{-\lambda})$; conversely, every $c\in \ohf d(\tilde\GR,\Bbb
C_{-\lambda})$ determines an invariant eigendistribution. In effect, this is a
reinterpretation of the Harish-Chandra matching conditions. Put differently,
there exists a natural surjective map 
$$
\ohf d(\tilde\GR,\Bbb C_{-\lambda}) \ \longrightarrow\  \Cal{IE}_{cc}(\lambda)
\tag4.15
$$
onto $\Cal{IE}_{cc}(\lambda)$, the space of invariant eigendistributions with
constant coefficients $c_{g,w}$ and infinitesimal character $\chi_\lambda$. For
regular $\lambda\in\fh^*$, this map is also injective. We had required earlier
that the coefficients  $ c_{g,x}$ of any \lq\lq character cycle" $c(\Theta)$
satisfy the symmetry condition
$$
c_{g,vx} = \operatorname{sgn}(v) c_{g,x} \qquad 
\text{for $v\in W_\lambda$}\,.
\tag4.16
$$
The map (4.15) becomes an isomorphism, even for singular $\lambda$, when
restricted to the subspace of cycles $c$ subject to the symmetry condition
(4.16).

We recall that we had fixed an orientation on $\GR$ which affects the sign of
the Poincar\'e duality map (4.13). Let us also fix a Haar measure $dg$ on
$\GR$\,. Together, these two choices determine a non-zero,
$\GR$-invariant form of top degree. Its complexification defines a
$G$-invariant, holomorphic form $\omega$ on $G$, of top degree. The quotient
map (4.8) allows us to  view the $e^\alpha$, with $\alpha\in\Phi$, as
holomorphic functions on $\tilde G$. The form $q^*\omega$ is divisible by the
product of the $(1-e^{-\alpha})$, $\alpha\in\Phi^+$ -- see \cite{A}, for
example -- so 
$$
\tilde \omega \ = \ \frac {q^*\omega}
{\tsize\prod_{\alpha\in\Phi^+}(1-e^{-\alpha})}
\tag4.17
$$
exists as  a well defined, $G$-invariant, holomorphic $d$-form  on
$\tilde G$. As such, it is locally integrable over any real algebraic
$d$-chain in $\tilde G$. We make sense of the integral
$$
\int_c (q^*\phi) \ \tilde\omega \qquad \qquad (\,c\in \ohf d(\tilde\GR,\Bbb
C_{-\lambda})\,)\,,
$$ 
for any test function $\phi\in C_c^\infty(\GR)$, by viewing $c$ as a collection
of sections of $\Bbb C_\lambda\subset \Cal O_{\tilde G}$ over the various
connected components of $\tilde \GR$\,. These sections are locally bounded on
$\tilde G$. Thus we can multiply these sections and the integrand, and
integrate the product over the various components of $\tilde\GR'$, using the
orientation induced from $\GR$\,. When $c$ corresponds to $\Theta\in
\Cal{IE}_{cc}(\lambda)$ via the homomorphism (4.15), 
$$
\int_\GR \Theta\phi\,dg \ = \ \int_c (q^*\phi)\ \tilde\omega\,.
\tag4.18
$$ 
This identity amounts to a slight rephrasing of the correspondence (4.15)
between cycles and invariant eigendistributions -- see the proof of lemma 6.2
below. Note that the form
$\tilde
\omega$ was normalized by the choice of the Haar measure $dg$.

The passage from the local formula (4.3) to the character cycle $c(\Theta)$ has
a counterpart on the Lie algebra. In analogy to (4.5-8), we define $\tilde\fg$
as the set of pairs $(\zeta,x)\in\fg\times X$ such that $\zeta \in \fb_x$,
projections $\pi:\tilde\fg\to X$, $q:\tilde\fg\to \fg$, $\tilde\fg\to \fh$, as
well as $\tilde\gr = q^{-1}\gr$ and $\tilde\gr' = q^{-1}(\gr')$. Note that
$\tilde \fg$ and the quotient map $\tilde\fg \to \fh$ have come up already, in
section 2, where $\tilde\fg$ was denoted by $\Cal B$; we are switching to
different notation now to emphasize the analogy with the group case. The
function $e^{\lambda-\rho}$, pulled back from $\fh$ to $\tilde\fg$, spans a
rank one $\Bbb C$-subsheaf $\Bbb C_\lambda\subset\Cal O_{\tilde\fg}$. Unlike in
the group case, this sheaf is trivial as a sheaf of $\Bbb C$-vector spaces. The
numerators $d_{E,x}e^{\lambda_x-\rho_x}$ in (4.4) determine a chain
$c=c(\theta)\in \ohf d(\tilde\gr',\Bbb C_{-\lambda})$, and the Harish-Chandra
matching conditions on the Lie algebra imply that $c(\theta)$ is a cycle,
i.e., $c(\theta)\in \ohf d(\tilde\gr,\Bbb C_{-\lambda})$. As before, we get a
natural surjective map from $\ohf d(\tilde\gr,\Bbb C_{-\lambda})$ to the
vector space of invariant eigendistributions on $\gr$, with constant local
coefficients and with infinitesimal character $\chi_\lambda$. This map relates
a cycle $c$ to the invariant eigendistribution $\theta$ via the formula
$$
\int_\gr \theta\phi \,d\zeta \ = \ \int_c (q^*\phi)\ \tilde\omega \qquad\qquad
(\,\phi\in C_c^\infty(\gr)\,)\,.
\tag4.19
$$
Here $d\zeta$ is the Euclidean measure on $\gr$, normalized so that $d\zeta$ and
$\exp^*dg$ coincide at the identity. Also, $\tilde\omega$ denotes
the holomorphic $d$-form on $\tilde\fg$ obtained by complexi\-fying the
Euclidean measure $d\zeta$, pulling it back to $\tilde\fg$, and dividing by
the product of the positive roots.
\vskip 1cm

\subheading{\bf 5. The fixed point formalism and Kashiwara's conjecture}
\vskip .5cm

The discussion in section 2 attaches invariant eigendistributions
$\Theta(\cf)$, $\theta(\cf)$ to objects $\cf\in \operatorname
 D_\GR(X)_{-\lambda}$\,. These, in turn, correspond to character cycles, in
$\ohf d (\tilde\GR,\Bbb C_{-\lambda})$ and $\ohf d (\tilde\gr,\Bbb
C_{-\lambda})$, respectively. Kashiwara \cite{Ka4} has conjectured a direct
geometric relationship between the sheaf $\cf$ and the character cycle
corres\-ponding to
$\Theta(\cf)$. This conjecture is equivalent, via (4.18), to an integral
formula for the virtual character $\Theta(\cf)$, and also provides geometric
expressions for the coefficients $c_{g,x}$ in the local formula (4.3) for
$\Theta(\cf)$. When $\GR$ happens to be compact, the local expressions for
the  $c_{g,x}$ reduce to the fixed point formula (1.5). There is one other
situation in which characters had been calculated by means of a fixed point
formula, prior to Kashiwara's conjecture: Hecht's formula for the
characters of holomorphic discrete series representations  \cite{Hec}, which
also follows easily from the local fixed point formula presented here. In this
section we recall Kashiwara's  conjecture, which will be proved -- along with its
counterpart on the Lie algebra -- in subsequent sections. We treat the case of a
general infinitesimal character, which makes the discussion a bit heavy. 
Kashiwara  only discusses the  much more transparent  case
$\lambda=\rho$, so the reader may want to consult \cite{Ka4} first.

The commuting actions of $G$  and the universal Cartan $H$  -- see (2.9) -- 
on the enhanced flag variety  $\hat X$ induce  an action map
$$ a : \pro \ \to \ \hat X\,. 
\tag5.1a
$$  The projection onto the factor $\hat X$ defines a second morphism 
$$ p: \pro \to\hat X\,.
\tag5.1b
$$  In analogy to the spaces $\tilde G$, $\tilde\GR$, which were introduced in
section 4, we set
$$
\aligned &\hat G \ = \ \{\,(g,h,\hat x)\in G\times H \times\hat X \ | \
a(g,h,\hat x) =
\hat x = p(g,h,\hat x)\,\}\,, \\ &\hat\GR \ = \  \hat G\  \cap \ (\pro)\,.
\endaligned
\tag5.2
$$ Note that there is a natural $H$-fibration $\hat G \to \tilde G$, obtained by
eliminating the $H$-factor in the product $G\times H \times\hat X$ and
mapping
$\hat X$ to $X$. Composing this fibration with the natural map (4.8) results in
a morphism $\tau:\hat G \to H$. We claim: the diagram
$$
\CD
\ \ \ \ \hat G \ \ \ \ @>{\subset}>>  G\times H \times\hat X
\\ @V{\tau}VV       @VV{p_H}V
\\
\ \ \ \ H \ \ \ \ @=  H 
\endCD
\tag5.3
$$ commutes. It suffices to check this over any given point $x\in X$. We can
identify $\hat X$ with  $G/N_x$ and $H$ with $B_x/N_x$; with these
identifications, the verification is straightforward. 

The rank one local system $\Bbb C_\lambda$ on $\tilde G$ was defined as the
pullback of a local system on $H$ via the projection (4.8) -- the local system
generated by the multiple valued function $e^{\lambda-\rho}$. The
commutativity of the diagram (5.3) allows us to think of $\Bbb C_\lambda$
as a \lq\lq universal"  rank one local system also on the spaces $\hat G$ and
$G\times H
\times\hat X$. We shall use the same symbol $\Bbb C_\lambda$ in all cases,
to avoid complicated notation. Note that $\Bbb C_\lambda$  is canonically a
subsheaf of the sheaf of holomorphic functions on the various spaces, in
analogy to (4.9), compatibly with the natural morphisms between the spaces.

We recall the definition, in section 2, of the twisted $\GR$-equivariant
derived category $ \operatorname D_{G_{\Bbb R}}(X)_{\lambda}$.
Disregarding part of the structure, we get a \lq\lq forgetful functor"
$$
 \operatorname D_{G_{\Bbb R}}(X)_{\lambda} \ \longrightarrow \  \db
{Sh_{X,\lambda}}
\tag5.4a
$$ into the bounded derived category of $(\lambda -\rho)$-monodromic
sheaves on
$X$. The latter derived category is a full subcategory of $\db {\hat X}$,
$$
\db {Sh_{X,\lambda}} \ \hookrightarrow \  \db {\hat X}\,.
\tag5.4b
$$ The composition of (5.1a,b) gives a functor
$$
 \operatorname D_{G_{\Bbb R}}(X)_{\lambda} \ \longrightarrow \  \db {\hat
X}\,,
\qquad
\cg\mapsto \hat\cg\,.
\tag5.4c
$$ Recall the definitions (5.1). We claim: for each  $\cg \in \operatorname
D_{G_{\Bbb R}}(X)_{\lambda}$\,, there exists a cano\-nical morphism
$$
\phi : a^*\hat\cg\to p^*\hat\cg\otimes\Bbb C_\lambda 
\tag5.5
$$  in $\db\pro$, which will play a crucial role in the fixed point formalism.

To construct $\phi$, we consider the diagram
$$  G_{\Bbb R}\times \frak h \times \hat X @>{\epsilon}>> \pro @>{a,p}>> \hat
X,
$$  where $\epsilon(g,\zeta,\hat x) = (g,\text{exp}(\zeta),\hat x)$. Note that
the compositions
$\e\, a$, $\e\, p$ are the action and projection morphisms of the
$(\GR\times \fh)$-action on $\hat X$. By definition, 
$\operatorname D_{G_{\Bbb R}}(X)_
\lambda$ is a full subcategory of the $(\GR\times\fh)$-equivariant derived
category on
$\hat X$. The formalism of this latter equivariant derived category implies
the existence of a distinguished isomorphism
$$
\phi' : \epsilon^*a^*\hat\cg \to \e^*p^*\hat\cg\,.
\tag5.6a
$$ The single isomorphism $\phi'$ encodes a  family of isomorphisms
$$
\aligned &\phi'_{g,\zeta} \ : \ \ell_{g,\zeta}^*\hat \cg \ @>{\ \sim \ }>> \hat
\cg \,,
\qquad \text{(\, $\ell_{g,\zeta}=$\ translation by
$(g,\zeta)\in\GR\times\fh$\,)\, },
\\ &\phi'_{g,\zeta} \ = \ \text{restriction of $\phi'$ to
$\{g\}\times\{\zeta\}\times\hat X$}\,,
\endaligned
\tag5.6b
$$ which depend multiplicatively and continuously on the parameters. Since
$\frak h$ is simply connected, we may view $e^{\lambda-\rho}$ as a
well-defined  section of $\e^*\Bbb C_\lambda$\,. We now define 
$$
\phi'' : \e^*a^*\hat\cg \to \e^*p^*\hat\cg\otimes \e^*\Bbb C_\lambda
\tag5.7
$$ by the formula $\phi'' = e^{\lambda-\rho} \phi'$. Any
$\zeta\in\fh_\Bbb Z =_{\text{def}}\operatorname{Ker}(\fh \to H)$ acts on
$\hat X$ as translation by
$\exp(\zeta) = 1$. The monodromicity condition in the definition of
$Sh_{X,\lambda}$ implies:
$$
\phi'_{e,\zeta} \ : \ \ell_{e,\zeta}^* \hat\cg = \hat \cg \ \to \ \hat \cg \ \
\text{is multiplication by $e^{-(\lambda-\rho)(\zeta)}$}\,,
\tag5.8
$$ for $\zeta\in\fh_\Bbb Z$, as before. It follows that
$\phi''$ is invariant under translation by the lattice $\fh_\Bbb Z$, and
therefore drops to an isomorphism $a^*\hat\cg @>{\sim}>>
p^*\hat\cg\otimes\Bbb C_\lambda$. That is the isomorphism $\phi$ whose
existence we postulated earlier. 

Two comments are in order. The statement (5.8) amounts to a
reinterpretation of the monodromicity condition in the language of the
equivariant derived category. Secondly,  the category $\operatorname
D_{G_{\Bbb R}}(X)_\lambda$ depends only on the image of $\lambda$ in
$\frak h^*/\Lambda$\,, where $\Lambda$ denotes the weight lattice, i.e.,
the lattice dual to $(2\pi i)^{-1}\fh_\Bbb Z$\,. The definition of
$\phi$, on the other hand, involves the section $e^{\lambda-\rho}$ of
$\e^*\Bbb C_\lambda$, and thus depends on $\lambda$ itself.

To each $\cg\in \operatorname D_{G_{\Bbb R}}(X)_\lambda$, we associate a
cycle in
$\Hom^{inf}_d( {\tilde G}_{\Bbb R} , \Bbb C_{-\lambda})$, of degree
$d=\operatorname{dim}(\GR)$, as follows. We consider the Cartesian
diagram
$$
\CD
\pro     @>s=(p,a)>>           \hat X \times \hat X  \\ @AAA                                
@AA{\Delta}A   \\ {\hat G}_{\Bbb R}      @>>>      \hat X
\endCD
\tag5.9
$$  where $\Delta$ is the diagonal map. Consider the following sequence of
maps
$$
\aligned &\operatorname{Hom}({\hat\cg},{\hat\cg})  \ @>{\ \sim \ }>>\ 
\Hom^0_{\Delta
\hat X}(\hat X\times\hat X ,
\Bbb D{\hat\cg} \boxtimes {\hat\cg}) 
\\ &\ \ \ \ \ @>{\ s^* \ }>>\  \Hom^0_{{\hat G}_{\Bbb R}}(\pro , s^*(\Bbb
D{\hat\cg}
\boxtimes {\hat\cg})) 
\\ &\ \ \ \   \ @>{\ \, \sim \ }>> \ \Hom^0_{{\hat G}_{\Bbb R}}(\pro , p^*\Bbb
D{\hat\cg} 
\otimes a^*{\hat\cg}) 
\\ &\ \ \ \  \ @>{\ \, \phi\ \, }>> \ \Hom^0_{{\hat G}_{\Bbb R}}(\pro , p^*\Bbb
D{\hat\cg}
\otimes p^*{\hat\cg} \otimes \Bbb C_\lambda) 
\\  &\ \ \ \ \ @>{\ \ \eta \, \ }>> \  \Hom^0_{{\hat G}_{\Bbb R}}(\pro , p^*\Bbb
D_{\hat X}
\otimes
\Bbb C_\lambda) 
\\ &\ \ \ \ \ \ \ \ \ \ \ @>{\ \sim\ }>> \ \Hom^{inf}_{d+r}( {\hat G}_{\Bbb R} ,
\Bbb C_{-\lambda}) @>{\ \sim \ }>> \ 
\Hom^{inf}_d( {\tilde G}_{\Bbb R} , \Bbb C_{-\lambda})\,,
\endaligned
\tag5.10
$$  with $r = \operatorname{dim}_{\Bbb R}(H)$. We use our previous
notational conventions: $\Hom^p_Z(\dots)$ denotes
local cohomology along
$Z$,   $\Hom_p^{inf}(\dots)$ is homology with infinite supports, and $d =
\operatorname{dim}(\GR)$. However, we define the Verdier dual  $\Bbb
D{\hat\cg}$ of  ${\hat\cg}$ using the dualizing sheaf
$\Bbb D_{\hat X}$  on $\hat X$; this differs formally from our definition in
section 2, but does not effect the results to be proved, as pointed out  in
footnote (2) in section 2.

Let us explain the various steps in (5.10). The isomorphism on the first line 
follows from the constructibility of ${\hat\cg}$ -- see, for example,
propositions 3.1.14 and 3.4.4 in \cite{KSa}. The morphism
$s^*$  is the pullback of local cohomology, and the next isomorphism follows
from the interpretation of the external tensor product  $\boxtimes$ in
terms of the usual (internal) tensor product of sheaves. The morphism
$\eta$ is induced by the duality pairing $\Bbb D{\hat\cg}\otimes{\hat\cg}
\to
\Bbb D$. On the most naive level, the next to last isomorphism can be
understood by identifying
$p^*\Bbb D_{\hat X}$ with the constant sheaf, appropriately shifted, and then
applying Poincar\'e duality; recall the reason for the passage from $\Bbb
C_\lambda$ to
$\Bbb C_{-\lambda}$ as explained in section 4. More formally, 
$$ p^*\Bbb D_{\hat X} \simeq p^! \Bbb D_{\hat X} [-(d+r)] \simeq \Bbb
D_{\pro}[-(d+r)]\,,
$$ hence
$$
\aligned &p^*\Bbb D_{\hat X}\otimes \Bbb C_\lambda \simeq (\Bbb
D_{\pro}\otimes \Bbb C_\lambda)[-(d+r)]
\\ &\ \ \ \ \simeq R\,\Cal H om (\Bbb C_{-\lambda}, \Bbb D_{\pro})[-(d+r)]
\qquad
\text{(\,$\Bbb C_{-\lambda}$ is locally free!\,)}
\\ &\ \ \ \  =_{\text{def}}\, \Bbb D (\Bbb C_{-\lambda})[-(d+r)]\,.
\endaligned
$$ Let $i$ denote the inclusion $\hat\GR \subset \pro$\,; then
$$
\aligned &\Hom^0_{{\hat G}_{\Bbb R}}(\pro , p^*\Bbb D_{\hat X}\otimes \Bbb
C_\lambda)
\simeq \Hom^{-d-r}_{{\hat G}_{\Bbb R}}(\pro ,\, \Bbb D (\Bbb C_{-\lambda}))
\\ &\ \ \ \ \simeq  \Hom^{-d-r}(\hat\GR\,, i^!\,\Bbb D (\Bbb C_{-\lambda}))
\qquad
\text{(\,by definition of local cohomology)}
\\ &\ \ \ \ \simeq  \Hom^{-d-r}(\hat\GR\,, \,\Bbb D (i^*\Bbb C_{-\lambda})) 
\\ &\ \ \ \ = 
\Hom^{-d-r}(\hat\GR\,, \,\Bbb D (\Bbb C_{-\lambda}))
\qquad\text{(\,according to our notational conventions) }
\\ &\ \ \ \ \simeq \Hom^{inf}_{d+r}( {\hat G}_{\Bbb R} ,\Bbb C_{-\lambda})
\qquad
\text{(\,by definition of homology)}\,.
\endaligned
$$ 
The last isomorphism in (5.10) follows from two facts: we are dealing with
homology (with locally finite support) in top degree, and $\hat\GR \to
\tilde\GR$ is a fibration with $r$-dimensional fibers. The symbol $\Bbb
C_\lambda$, we recall, refers to sheaves on the various spaces, all coming 
by pull back from the same local system on $H$, so $\Bbb C_\lambda$ on
$\hat\GR$ is the pullback of $\Bbb C_\lambda$ on $\tilde\GR$.

We can assign a cycle in $\Hom^{inf}_d({\tilde G}_{\Bbb R} , \Bbb
C_{-\lambda})$ to any given $\cg\in  \operatorname D_{G_{\Bbb
R}}(X)_\lambda$ by taking the image of
$1_{\hat\cg}\in\operatorname{Hom}({\hat\cg},{\hat\cg})$ under the
sequence of maps (5.10). As a matter of notation, we write
$$ c(\cf) \,\in\, \Hom^{inf}_d({\tilde G}_{\Bbb R} , \Bbb C_{-\lambda})
\qquad\qquad (\,\cf\in \operatorname D_{G_{\Bbb R}}(X)_{-\lambda}\,)
\tag5.11
$$ 
for the cycle attached to $\cg = \Bbb D\cf$\,. Recall (2.14): since the
definition of $D_{G_{\Bbb R}}(X)_{\lambda}$ involves the monodromy
behavior of $e^{(\lambda-\rho)}$, Verdier duality maps $D_{G_{\Bbb
R}}(X)_{-\lambda}$ to $D_{G_{\Bbb R}}(X)_{\lambda +2\rho}$; on the other
hand, $D_{G_{\Bbb R}}(X)_{\lambda +2\rho}=D_{G_{\Bbb R}}(X)_{\lambda}$
since $2\rho$ is integral, so $\Bbb D:D_{G_{\Bbb R}}(X)_{-\lambda} \to
D_{G_{\Bbb R}}(X)_{\lambda }$\,.

\proclaim{5.12 Theorem} The cycle $c(\cf)$, for $\cf\in \operatorname
D_{G_{\Bbb R}}(X)_{-\lambda}$,\,
\ 
 is the character cycle of the virtual character $\Theta(\cf)$. In particular,
$$
\int_\GR \Theta(\cf) \phi\, dg \ = \ \int_{c(\cf)}(q^*\phi)\tilde \omega\,,
$$  
for any test function $\phi\in C^\infty_c(\GR)$.
\endproclaim

This statement was conjectured by Kashiwara \cite{Ka4}. We shall prove it in
the later sections. We remind the reader that the homomorphism (4.15) is
not injective, in general. However, it has a distinguished left inverse,
namely the assignment $\Theta \mapsto$\ character cycle of $\Theta$\ \  -- 
\, the character cycle, by definition, satisfies the symmetry condition
(4.16).  Our theorem asserts, in particular, that $c(\cf)$ satisfies this
symmetry condition.

The theorem -- more precisely, the construction preceding the theorem -- 
can be simplified whenever $\lambda$ is integral. In this situation, the
twisted equivariant derived category $\operatorname D_{G_{\Bbb
R}}(X)_{-\lambda}$ agrees with the usual equivariant derived category 
$\operatorname D_{G_{\Bbb R}}(X)$ (recall 2.12-13). Also, the \lq\lq
universal" local system $\Bbb C_\lambda$ has a distinguished, globally
defined generating section $e^{\lambda-\rho}$. The cycle
$c(\cf)\in
\Hom^{inf}_d({\tilde G}_{\Bbb R} , \Bbb C_{-\lambda})$ can now be written as
a product of the generating section $e^{\lambda-\rho}$ with the absolute
cycle
$c(\cf_\rho)\in \Hom^{inf}_d({\tilde G}_{\Bbb R} , \Bbb C)$, where
$\cf_\rho\in
\operatorname D_{G_{\Bbb R}}(X)$ is the image of $\cf$ under the periodicity
isomorphism (2.12-13). This allows us to work on $\GR\times X$, rather than
$\pro$ as before. In particular, the existence of the canonical isomorphism
$\phi: a^*\cg \to p^*\cg$ immediately follows from the definition of the
equivariant derived category, so (5.6-8) become superfluous, and (5.9-10)
simplify correspondingly. Kashiwara's conjecture \cite{Ka4} is phrased in
these terms.

The proof of theorem (5.12) depends on a local fixed point formula for the
coefficients $c_{g,x}$ in the local expression (4.3) for $\Theta =
\Theta(\cf)$. This local formula is of independent interest. Kashiwara
 develops a fixed point formalism in \cite{Ka4} and uses it to calculate the
cycle (5.11) from the geometric datum of the sheaf $\cf$.

To state the local formula, we fix a regular semisimple 
$g \in G'_\Bbb R$  and  a fixed point $x\in X$ of $g$. The centralizer of $g$ in
$\GR$ is a   Cartan subgroup $T_\Bbb R$\,, whose Lie algebra we denote by
$\ft_\Bbb R$\,; $\ft$ is the complexification of  $\ft_\Bbb R$. The
identification 
$\ft\cong
\fb_x/[\fb_x,\fb_x]\cong\fh$ lifts to a natural identification
$$ T\ \cong \ B_x/[B_x,B_x]\ \cong \  H\,.
\tag5.13
$$
 It induces an
identification $\alpha\mapsto \alpha_x$ between the universal root
system $\Phi$ and the concrete root system $\Phi(\fg,\ft)$. The 
complexified Cartan subgroup
$T$ normalizes the Borel subalgebra
$$
\frak b_x \ = \ \frak t \oplus (\oplus_{\alpha \in \Phi^+} \ \frak g^
{-\alpha_x})\,.
$$ 
Here $\fg^{\alpha_x}$, for
$\alpha\in\Phi$, denotes the
$\alpha_x$-root space of
$(\fg,\ft)$.  Let
$$
\fn^+(g,x) \ = \ \oplus_{\alpha \in \Phi^+} \ \frak g^ {\alpha_x}
\tag5.14
$$ denote the nilpotent radical of the opposite Borel, and $\fn'(g,x)$, 
$\fn''(g,x)$ two subalgebras of $\fn^+(g,x)$, of the following type. We
choose  subsets 
$\Psi',\Psi''\subset \Phi^+$, such that
$$
\aligned
\text{a) \ \ }&\text{for\ \ } \alpha\in\Phi^+ \,, \ \ \text{if \ \
}|e^{\alpha_x}(g)|\neq 1
\,, \text{\ \ then}
\\
&\alpha
\in\Psi'\ \Longleftrightarrow  \ |e^{\alpha_x}(g)|<1\ \ \ \text{and}  \ \ \ 
\alpha
\in\Psi''\ \Longleftrightarrow  \ |e^{\alpha_x}(g)|>1\,;
\\
\text{b) \ \ }&\text{for\ \ } \alpha_1\,,\alpha_2\in\Phi^+ \,, \ \ \text{if \ \
} \alpha_1+\alpha_2\in\Phi^+
\,, \text{\ \ then}
\\
&\alpha_1,\alpha_2\in\Psi' \ \Longrightarrow \ 
\alpha_1+\alpha_2\in\Psi'\, ,\ \ \ \alpha_1,\alpha_2\in\Psi'' \ 
\Longrightarrow \ 
\alpha_1+\alpha_2\in\Psi''\,.
\endaligned
\tag5.15
$$ 
Subsets $\Psi',\Psi''$ satisfying these conditions do exist: for example, the
subsets defined by a) without the restriction $|e^{\alpha_x}(g)|\neq 1$
satisfy also b). Because of b), 
$$
\fn'(g,x) \ = \ \oplus_{\alpha \in \Psi'} \ \frak g^ {\alpha_x}\,, \qquad
\fn''(g,x) \ = \ \oplus_{\alpha \in \Psi''} \ \frak g^ {\alpha_x}\,,
\tag5.16
$$ 
are subalgebras of $\fn^+(g,x)$. Our eventual statement, theorem 5.24, will
not depend on the particular choice of
$\Psi$. Indeed, the theorem is valid even with more general choices
of $\fn'(g,x)$ and $\fn''(g,x)$ -- see the comment at the end of this
section. Further notation:
$$
\gathered N^+(g,x) \ = \  \exp(\fn^+(g,x))x\,, 
\\ N'(g,x) \ = \  \exp(\fn'(g,x))x\,, \qquad N''(g,x) \ = \  \exp(\fn''(g,x))x\,.
\endgathered
\tag5.17
$$ Then $N^+(g,x)\subset X$ is an open Schubert cell which contains $N'(g,x)$
and
$N''(g,x)$ as affine linear subspaces, both invariant under the action of $T$. 

The inverse image of  $N^+(g,x)$ in $\hat X$ splits into a product 
$N^+(g,x)\times H$, and the $G\times H$-action on $\hat X$ restricts to a
$T\times H$-action on $N^+(g,x)\times H$\,:
$$ (t,h)\,:\,(\exp(\zeta)x,h_1) \ \mapsto \
(\exp((\operatorname{Ad}t)(\zeta))x,t\,h_1\,h^{-1})\,,\qquad
\zeta\in\fn^+(g,x)\,;
\tag5.18
$$ here we identify $t\in T$ with an element of the universal Cartan $H$ via
the isomorphism (5.13). In particular, there is a natural restriction functor
$$
 \operatorname D_{G_{\Bbb R}}(X)_\lambda \to  \operatorname D_{T_{\Bbb
R}}(N^+(g,x))_\lambda\,, 
\tag5.19a
$$
from the twisted equivariant derived category $\operatorname
D_{G_{\Bbb R}}(X)_\lambda$, viewed as the $\GR$-equi\-variant derived
category of $(-\lambda-\rho)$-monodromic sheaves on $\hat X$, to the 
twisted equivariant derived category $\operatorname D_{T_{\Bbb
R}}(N^+(g,x))_\lambda$, viewed as the $T_\Bbb R$-equi\-variant derived
category of $(-\lambda-\rho)$-monodromic sheaves on $N^+(g,x)\times H$. 
The base of the bundle  $N^+(g,x)\times H \to N^+(g,x)$ is contractible.   It
follows that there is a canonical equivalence  between the category of
$(\lambda -\rho)$-monodromic sheaves on $N^+(g,x)$ on the one hand, and
the category of (ordinary) sheaves on $N^+(g,x)$ on the other,
$$
\operatorname D_{T_{\Bbb R}}(N^+(g,x))_\lambda \  @>{\ \sim \ }>>\
\operatorname D_{T_{\Bbb R}}(N^+(g,x))\,;
\tag5.19b
$$ 
in one direction, 
sheaves on $N^+(g,x)$ can be pulled back
to $N^+(g,x) \times H$ and tensored with the universal local system $\Bbb
C_\lambda$, and in the opposite direction, twisted sheaves on $N^+(g,x)
\times H$ can be restricted to  $N^+(g,x)\times\{e\} \cong N^+(g,x)$. Thus
we get a restriction functor
$$
 \operatorname D_{G_{\Bbb R}}(X)_\lambda \to \operatorname D_{T_{\Bbb
R}}(N^+(g,x))\,, \qquad \cg \mapsto \cg(x)\,,
\tag5.19c
$$
the composition of the functors (5.19a,b).

The action and projection maps $a,p:\GR\times H\times \hat X \to \hat X$
restrict to maps from $T_\Bbb R\times H \times N^+(g,x)
\times H$ to $N^+(g,x)
\times H$, to which we refer by the same letters $a,p$. We then get a
commutative diagram
$$
\CD
T_\Bbb R\times H \times N^+(g,x)\times H @>{\ a,p \ }>>  N^+(g,x)\times H
\\
@A{(t,y)\mapsto (t,t,y,e)}AA @AA{y\mapsto (y,e)}A
\\ 
T_\Bbb R\times N^+(g,x) @>{\ a,p \ }>>  \ \ N^+(g,x)\ \ ;
\endCD
\tag5.20
$$
note that the action map in the top row is the action (5.18), and that  the
universal local system $\Bbb C_\lambda$ \lq\lq lives" on the first of the
two factors $H$ in  $T_\Bbb R\times H \times N^+(g,x)\times H$. The
morphism (5.5), which came from the structure of the twisted equivariant
derived category, restricts to a morphism $\phi: a^*\cg \to p^*\cg \otimes
\Bbb C_\lambda$ in the derived category $\db{T_\Bbb R\times H \times
N^+(g,x)\times H}$, for each $\cg\in \operatorname D_{G_{\Bbb
R}}(X)_\lambda$. In (5.5), we wrote $\hat \cg$ to signify that we were
disregarding the equivariance. To be consistent, we might use the same
convention now; however, to avoid proliferating notation, we shall not
notationally distinguish between $\cg$ and $\hat \cg$ from now on.  Via the
commutative diagram (5.20), our present incarnation of $\phi$ induces   a
morphism
$$
\phi\,:\,a^*\cg(x) \ \longrightarrow \   p^*\cg(x)\otimes\Bbb C_\lambda
\tag5.21
$$ 
in the category $\db{T_\Bbb R \times N^+(g,x)}$. In the statement that
follows, we think of $x$ as fixed and $g$ as variable within $T'_\Bbb R$. To
emphasize the point, we shall denote this variable element as $t$. Note that 
the subspaces  $N'(t,x)$, $N''(t,x)$, defined in analogy to  $N'(g,x)$ and
$N''(g,x)$, depend only on the connected component of $T'_\Bbb R$ in which
$t$ lies, whereas
$N^+(t,x)=N^+(g,x)$ for all $t\in T'_\Bbb R$. 

To keep the notation simple, we denote the restriction of the morphism
(5.21) to the slice $\{t\}\times N^+(t,x)$ by the same letter $\phi$. This
gives us a morphism 
$$
\phi\ : \ t^*\cg(x) \ \longrightarrow \  \cg(x)\otimes\Bbb C_\lambda\,.
\tag5.22
$$ The diffeomorphism $t:N^+(t,x) \to N^+(t,x)$ preserves the subspaces
$N'(t,x)$,
$N''(t,x)$. Hence  $\phi$ induces  morphisms 
$$
\aligned &\phi\ : \ t^*\Cal H_{N'(t,x)}(\cg(x)) \ \longrightarrow \  \Cal
H_{N'(t,x)}(\cg(x))\otimes\Bbb C_\lambda\,,
\\ &\phi\ : \ t^*(\cg(x)|_{N''(t,x)}) \ \longrightarrow \ 
(\cg(x)|_{N''(t,x)})\otimes\Bbb C_\lambda\,.
\endaligned
\tag5.23a
$$ Since $x$ is a fixed point of $t$, (5.23a) induces morphisms of stalks,
respectively costalks -- i.e., local cohomology at a point -- 
$$
\aligned &\phi_t\ : \ \Cal H_{N'(t,x)}(\cg(x))_x \ \longrightarrow \  \Cal
H_{N'(t,x)}(\cg(x))_x\otimes\Bbb C_\lambda\,,
\\ &\phi_t\ : \ \Cal H_{\{x\}}(\cg(x)|_{N''(t,x)}) \ \longrightarrow \  \Cal
H_{\{x\}}(\cg(x)|_{N''(t,x)})\otimes\Bbb C_\lambda\,.
\endaligned
\tag5.23b
$$ We can now state the local fixed point theorem, which gives us two
separate, formally dual expressions for the  section $c_{t,x}$ of $\Bbb
C_\lambda$ appearing in the expression (4.11) of the virtual character
$\Theta = \Theta(\cf)$:

\proclaim{Theorem 5.24} For $\Cal F \in \operatorname D_{G_{\Bbb
R}}(X)_{-\lambda}$ and $t\in T'_\Bbb R$, the  section $c_{t,x}$ of $\Bbb
C_\lambda$  satisfies
$$ 
\aligned c_{t,x}\,  &= \, \tsize\sum (-1)^i \,\operatorname{tr}(\,\phi_t:\Cal
H_{N'(t,x)}^i(\Bbb D
\cf(x))_x
\,
\to
\, 
\Cal H_{N'(t,x)}^i(\Bbb D \cf(x))_x\otimes\Bbb C_\lambda\,)
\\ &= \, \tsize\sum (-1)^i \,\operatorname{tr}(\,\phi_t:\Cal H_{\{x\}}^i((\Bbb D
\cf(x))|_{N''(t,x)})
\to 
\Cal H_{\{x\}}^i((\Bbb D \cf(x))|_{N''(t,x)})\otimes\Bbb C_\lambda\,)\,.
\endaligned
$$
\endproclaim

Note that the presence of the factors $\Bbb C_\lambda$ makes the traces
into sections, as they should be, rather than scalars.

Let us suppose, for the moment, that a specific local section of the sheaf
$\Bbb C_\lambda$ on $T'_\Bbb R$ has been chosen near $t$, and let us denote
this section by the symbol $e^{\lambda-\rho}$, for simplicity. In this
situation
$$
\aligned &e^{-(\lambda-\rho)}\phi_t:\Cal H_{N'(t,x)}^i(\Bbb D \cf(x))_x
\to
\Cal H_{N'(t,x)}^i(\Bbb D \cf(x))_x\,, \\ &e^{-(\lambda-\rho)}\phi_t:\Cal
H_{\{x\}}^i((\Bbb D \cf(x))|_{N''(t,x)})
\to
\Cal H_{\{x\}}^i((\Bbb D \cf(x))|_{N''(t,x)})\,,
\endaligned
$$ are well defined, and 
$$
\gathered c_{t,x}\ = \ d_{t,x} \, e^{\lambda-\rho}\,, \ \ \text{with  }
\\ d_{t,x} = \tsize\sum (-1)^i \,\operatorname{tr}(e^{\rho-\lambda}\phi_t:\Cal
H_{N'(t,x)}^i(\Bbb D
\cf(x))_x
\to 
\Cal H_{N'(t,x)}^i(\Bbb D \cf(x))_x\,) =
\\
\tsize\sum (-1)^i \,\operatorname{tr}(e^{\rho-\lambda}\phi_t:\Cal
H_{\{x\}}^i((\Bbb D
\cf(x))|_{N''(t,x)})
\to 
\Cal H_{\{x\}}^i((\Bbb D \cf(x))|_{N''(t,x)})\,;
\endgathered
\tag5.25a
$$ here $d_{t,x}$ is now a scalar. In particular, this applies for any regular
$t\in T_\Bbb R$ near the identity: we choose the branch of
$e^{\lambda-\rho}$ with value 1 at the identity. For such $t$,
$e^{\rho-\lambda}\phi_t$ is homotopic to the identity map, so 
$$ d_{t,x} \ = \ \chi(\,\Cal H_{N'(t,x)}^*(\Bbb D \cf(x))_x\,) \ = \ \chi(\,\Cal
H_{\{x\}}^*(\Bbb D \cf(x)|_{N''(t,x)})\,)
\tag5.25b
$$ becomes an Euler characteristic. The resulting formulas can be continued
to the connected component of $t$ in $T'_\Bbb R$\,. Every connected
component $C$ of
$T'_\Bbb R\cap T_\Bbb R^0$ contains the identity in its closure, hence
$c_{t,x}$ either vanishes identically on $C$, or equals an integral multiple of
a well defined branch of $e^{\lambda-\rho}$. One can argue similarly for
connected components of $T'_\Bbb R$ outside $T_\Bbb R^0$, where $c_{t,x}$
can be expressed as a cyclotomic integer times a well defined branch  of
$e^{\lambda-\rho}$, unless $c_{t,x}$ vanishes identically.

The analogues on the Lie algebra of theorems 5.12 and 5.24 are formally
simpler, because the local system $\Bbb C_\lambda$ on $\tilde\fg$ has a
canonical global generating section $e^{\lambda-\rho}$. We fix a particular
point $(\zeta,x)\in \tilde \gr'$\,, write $T$ for the centralizer of $\zeta$ in
$G$, and define 
$$
\gathered N^+(\zeta,x) = N^+(\exp(\zeta),x)\,, 
\\ N'(\zeta,x) = N'(\exp(\zeta),x)\,,
\qquad  N''(\zeta,x) = N''(\exp(\zeta),x)\,.
\endgathered
\tag5.26
$$ Let $E$ be the connected component of $\ft'_\Bbb R$ containing $\zeta$\,.  

\proclaim{5.27 Theorem} For $\Cal F \in \operatorname D_{G_{\Bbb
R}}(X)_{-\lambda}$\,,  the  constant $d_{E,x}$  appearing in the local
expression (4.4) of the virtual character $\theta =
\theta(\cf)$ is given by the formula
$$   d_{E,x} = \chi(\,\Cal H_{N'(\zeta,x)}^*(\Bbb D \cf(x))_x\,) = \chi(\,\Cal
H_{\{x\}}^*(\Bbb D \cf(x)|_{N''(t,x)})\,)\,.
$$
\endproclaim

The character cycle $c\in \Hom^{inf}_d({\tilde \fg}_{\Bbb R} , \Bbb
C_{-\lambda})$ of the virtual character $\theta(\cf)$ is the product of the
generating section $e^{\lambda-\rho}$ with an absolute cycle $c_\rho\in
\Hom^{inf}_d({\tilde \fg}_{\Bbb R} , \Bbb C)$, just as in the simplified version
of theorem 5.12 for integral $\lambda$ -- see the discussion below the
statement of that theorem. We shall omit a formal statement of the Lie
algebra version of theorem 5.12, since it is a direct analogue of the
simplified statement -- with $\lambda$ integral -- of 5.12. 

We remarked earlier that the hypotheses of  theorem 5.24 are
unnecessarily restrictive. First of all, the condition on $\Psi'$ and $\Psi''$
embodied by (5.15a) needs to be satisfied only by roots $\alpha\in\Phi^+$
such that $e^{\alpha_x}$ is real valued and positive near $g$. To {\it
prove\/} this more general version would be considerably more involved,
since it does not directly follow from the fixed point formalism of either
Goresky-MacPherson \cite{GM} or Kashiwara \cite{KSa}.  Secondly,
condition (5.15b) serves the purpose of making $\fn'(g,x)$ and
$\fn''(g,x)$ subalgebras, which is  natural from a group
theoretic point of view. This condition is used crucially in section 8, where
we show that theorem 5.24 is compatible with parabolic induction. It is not
used in section 6, however, when we show that theorems 5.12 and 5.24 are
equivalent.
\vskip 1cm

\subheading{\bf 6. Formal aspects of the proofs}
\vskip .5cm

In broad outline, we shall prove our main theorems -- the character formulas
3.8, 5.12, 5.24, 5.27 -- by verifying them for certain \lq\lq standard sheaves".
In this section, we describe the standard sheaves, and we argue that they
generate   the K-groups of the equivariant derived categories
$\operatorname D_\GR(X)_{-\lambda}$\,. In subsequent sections we break
up  the verification for standard sheaves into several steps, which we then
carry out. 

To begin with, the main theorems descend to the level of the K-group: if
three sheaves $\cf_i\in\operatorname D_\GR(X)_{-\lambda}$ fit into a
distinguished triangle
$$
\cf_1\ \to \ \cf_2 \ \to \ \cf_3 \ \to \cf_1[1]
\tag6.1
$$ and if any one of the main theorems holds for two of the $\cf_i$\,, then it
holds for the third. In the case of theorem 3.8, this follows from (3.5) and the
additivity of the integrals in 3.8. Theorems 5.24 and 5.27 involve alternating
sums of traces, and these, too, behave additively\footnote{In the case of
5.24, this depends on the functoriality of the morphism $\phi$, which insures
that the induced maps on cohomology determine a morphism of the long
exact cohomology sequences.} in distinguished triangles.  The additivity of
5.12 is less obvious; we establish it by appealing to the following
proposition.

\proclaim{6.2 Proposition} Either of the two formulas for the $c_{t,x}$ in
theorem 5.24, for any particular sheaf
$\cf\in\operatorname D_\GR(X)_{-\lambda}$\,, for all $x\in X$ and all regular
$t$ fixing $x$, is equivalent to the statement of theorem 5.12 for the sheaf
$\cf$.
\endproclaim

We postpone the proof to the end of this section. Let us summarize our
conclusions so far. 

\proclaim{6.3 Remark} It suffices to establish the main theorems for a
collection of sheaves which generate the K-group of $\operatorname
D_\GR(X)_{-\lambda}$\,.
\endproclaim

Standard sheaves are associated to $\GR$-orbits in the flag variety. Let
$S\subset X$ be an orbit and $\hat S\subset \hat X$ its inverse image in the
enhanced flag variety. We fix $\lambda\in\fh^*$ and consider an irreducible,
$\GR$-equivariant $(-\lambda-\rho)$-monodromic local system $\Cal L$ on
$S$ -- in other words, an irreducible, $\GR\times\fh$-equivariant local
system on
$\hat S$ whose monodromy along the fibers of $\hat S\to S$ agrees with the
monodromy of the multiple valued function $e^{-\lambda-\rho}$ on $H$. The
direct image $Rj_*\Cal L$ of $\Cal L$ under the inclusion $j$ of $S$ into $X$ is
an object of the twisted equivariant derived category $\operatorname
D_\GR(X)_{-\lambda}$\,. By definition, the sheaves of this type are the
standard sheaves in this category.

\proclaim{6.4 Lemma} The standard sheaves in $\operatorname
D_\GR(X)_{-\lambda}$ generate $K(\operatorname D_\GR(X)_{-\lambda})$.
More precisely, the standard sheaves generate $\operatorname
D_\GR(X)_{-\lambda}$ as a triangulated category.
\endproclaim

\demo{Proof} The support $\operatorname{Supp}(\cf)$ of any given $\cf\in
\operatorname D_\GR(X)_{-\lambda}$ is necessarily $\GR$-invariant. Let $S =
\cup\, S_j$ be the union of the orbits $S_j$  of maximal dimension in
$\operatorname{Supp}(\cf)$, \ $j:S \hookrightarrow X$ the inclusion. The
natural morphism $\cf \to Rj_*j^*\cf$ induces an isomorphism over $S$,
hence the third term in the distinguished triangle
$$
\cf \ \to \ Rj_*j^*\cf \ \to \ \cg \ \to \ \cf[1]
$$ has support of strictly lower dimension than $\operatorname{Supp}(\cf)$,
or vanishes identically. Arguing by induction, we may now assume that
$\cf=Rj_*j^*\cf$. The cohomology sheaves of $j^*\cf\in\operatorname
D_\GR(S)_{-\lambda}$ are $\GR$-equivariant,
$(-\lambda-\rho)$-monodromic local systems on $S$. As such, they must be
extensions of irreducible local systems on the various $S_j$\,. The lemma
follows.

\enddemo

In spite of the geometric terminology, irreducible  $\GR$-equivariant,
$(-\lambda-\rho)$-mono\-dromic local systems are essentially algebraic
objects. We recall how to construct them in terms of Lie theoretic data. This
involves, first of all, the enumeration of the
$\GR$-orbits in $X$ as in \cite{Ma,W}. Let
$T_\Bbb R$ be a Cartan subgroup, and $\tau_{x_0} : \ft \to \fh$ a concrete
isomorphism between the complexified Lie algebra $\ft$ of $T_\Bbb R$ and
the universal Cartan
$\fh$\,, corresponding to a fixed point ${x_0}\in X$ of $T_\Bbb R$. Note that
the
$\GR$-orbit
$$ S \ = \ S(T_\Bbb R,\tau_{x_0}) \ = \ \GR \cdot {x_0}
\tag6.5
$$   remains unchanged when $T_\Bbb R$ is replaced by a $\GR$-conjugate or
${x_0}$ by a
$N_{\GR}(T_\Bbb R)$-translate. The correspondence $(T_\Bbb
R,\tau_{x_0})\mapsto S$ sets up a bijection between the set of  $\GR$-orbits 
and the set of  pairs $(T_\Bbb R,\tau_{x_0})$, modulo the conjugacies just
mentioned. 

We fix the datum of a particular $\GR$-orbit $S= S(T_\Bbb R,\tau_{x_0})$ and
choose a character $\chi:T_\Bbb R \to \Bbb C^*$ such that $d\chi =
\tau_{x_0}^*(\lambda-\rho)$\,, as complex linear functions on $\ft$\,. A
parenthetical comment: if
$\GR$ is non-linear, contrary to our standing assumption, $T_\Bbb R$ may
fail to be abelian; in that case, we need to allow irreducible representations
$\chi$ of dimension greater than one, with differential equal to a direct sum
of copies  of $\tau_{x_0}^*(\lambda-\rho)$. The isotropy subgroup of $\GR$
at
${x_0}\in X$ splits into a direct product
$$ (\GR)_{x_0} \ = \ T_\Bbb R \cdot (N_{x_0}\cap \GR)\,,
\tag6.6a
$$ with $N_{x_0}=\exp\,[\fb_{x_0},\fb_{x_0}]$\,. The group $N_{x_0}\cap \GR$
is connected, so
$\chi$ lifts canonically to a character $\tilde\chi:(\GR)_{x_0} \to \Bbb C^*$,
with
$\tilde\chi = 1$ on $N_{x_0}\cap \GR$\,. When we identify the enhanced flag
variety $\hat X$ with $G/N_{x_0}$ as usual, the identity coset corresponds to
a point $\hat {x_0}\in \hat X$ lying over ${x_0}$. In particular, $\hat {x_0}$
lies in
$\hat S$, the inverse image of $S$ in $\hat X$. The group $\GR\times\fh$ acts
transitively on $\hat S$, with isotropy group
$$
\aligned (\GR\times\fh)_{\hat {x_0}} \ &= \ \{\,(g,\zeta)\in
(\GR)_{x_0}\times\fh
\ | \ g =
\exp(\tau_{x_0}^{-1}\zeta) \,\}
\\ &= \ \{\,(t,n,\zeta)\in T_\Bbb R \times (N_{x_0}\cap \GR) \times \fh \ | \ t
=
\exp(\tau_{x_0}^{-1}\zeta) \,\}
\endaligned
\tag6.6b
$$  at $\hat {x_0}$. A triple $(t,n,\zeta)$ lies in the identity component
precisely when $\tau_{x_0}^{-1}\zeta\in \ft_\Bbb R$\,, hence
$$
\aligned &\text{component group of\,\ }(\GR\times\fh)_{\hat {x_0}} \ \cong 
\\  &\{\,(t,\zeta)\in T_\Bbb R\times\fh \ | \ t =
\exp(\tau_{x_0}^{-1}\zeta) \,\}/\{\,(t,\zeta) \ | \ t =
\exp(\tau_{x_0}^{-1}\zeta) \,, \ \tau_{x_0}^{-1}\zeta\in\ft_\Bbb R\}\,.
\endaligned
\tag6.7
$$ Note that the character
$$
\psi\,:\,(\GR\times\fh)_{\hat {x_0}} \ \to \ \Bbb C^*, \qquad \psi(g,\zeta )=
\tilde
\chi(g) e^{-(\lambda-\rho)(\zeta)}\,,
\tag6.8
$$ is identically equal to 1 on the identity component of
$(\GR\times\fh)_{\hat {x_0}}$\,, and thus induces a character of the
component group. As such, it defines a $\GR\times\fh$-equivariant local
system $\Cal L_\chi$ on $\hat S$, of rank 1, which is therefore irreducible.

\proclaim{6.9 Lemma} The local system $\Cal L_\chi$ is
$(-\lambda-\rho)$-monodromic. Every irreducible
$\GR$-equivariant, $(-\lambda-\rho)$-monodromic local system on $S$ is of
this form, with uniquely determined $\chi$.
\endproclaim

\demo{Proof} Because of the $\GR$-equivariance, it suffices to check the
monodromicity condition on the fiber over ${x_0}$. On that fiber it has the
monodromy of the function $e^{-(\lambda-\rho)}$, hence the monodromy
of $e^{-(\lambda+\rho)}$ since $2\rho$ is integral. Conversely, the datum of
an irreducible
$\GR$-equivariant, $(-\lambda-\rho)$-monodromic local system on $S$ is
equivalent to that of an irreducible representation $\psi$ of the component
group (6.7), subject to the condition 
$$
\psi(e,\zeta) \ = \ e^{-(\lambda-\rho)(\zeta)} \ \ \text{for $\zeta\in\fh_\Bbb
Z$}
\tag6.10
$$ which reflects the monodromicity condition (recall: $\fh_\Bbb Z$
is the kernel of the exponential map $\exp:\fh \to H$\,). Note that
$\psi$ must be one dimensional since the component group is abelian. We
reconstruct the datum of the character $\chi:T_\Bbb R \to \Bbb C^*$ by
defining
$\chi(t) = \psi(t,\zeta) e^{(\lambda-\rho)(\zeta)}$ for any $\zeta\in\fh$ such
that $t= \exp(\tau_{x_0}^{-1}(\zeta))$ -- the condition (6.10) ensures that
the particular choice of $\zeta$ does not matter. When $t$ happens to lie the
identity component of $T_\Bbb R$, we can choose $\zeta$ to lie in
$\tau_{x_0}(\ft_\Bbb R)$, in which case $\psi(t,\zeta) = 1$. This shows that
$d\chi = \tau_{x_0}^*(\lambda-\rho)$, as required. 

\enddemo

In effect, the lemma gives us an explicit description of the standard sheaves
in $\operatorname D_\GR(S)_{-\lambda}$. To complete the verification 6.3,
we still need to prove proposition 6.2. This involves two major steps: the
reduction from the general, twisted situation to the untwisted case, and
secondly, the computation of the local fixed point contribution in the
untwisted case. This second ingredient already appears in Kashiwara's
announcement \cite{Ka4}.

\demo{Proof of proposition 6.2} Let us assume that the statement of
theorem 5.24 is satisfied for a particular   $\cf\in \operatorname
D_\GR(S)_{-\lambda}$ and all possible choices of a Cartan subgroup $T_\Bbb
R \subset G_\Bbb R$, of a regular element $t\in T_\Bbb R'$, and of fixed point
$x$ of $T$. We must  show that the character identity in theorem 5.12 holds,
and conversely, that the identity in theorem 5.12 implies the statement of
theorem 5.24 for all data $T_\Bbb R$, $t\in T_\Bbb R'$\,, and $x\in X^T$. The
virtual character
$\Theta = \Theta(\cf)$ is completely determined by its values on the regular
set, where it is real analytic. The value of the integral on the right of the
identity is also completely determined by the contribution lying over the
regular set: the inverse image of the singular set in $\tilde \GR$ has
codimension at least one, and can therefore be neglected in the integral. In
particular, instead of considering arbitrary test functions in the character
identity, it suffices to consider the case when the test function
$\phi$ is a delta function supported at an arbitrary regular point -- say $\phi
= \delta_t$ with $t\in T_\Bbb R'$ as before. In this situation, the left hand
side of the identity reduces to 
$$
\int_\GR\Theta\, \delta_t \,dg\ = \ \Theta(t)  \ = \ \sum_{x\in X^{T}}\ \frac {
c_{t,x}} {\tsize \prod_{\alpha \in {\Phi}^{ +}} (1-e^{-\alpha_x})(t )} 
\tag6.11a
$$ with the coefficients $c_{t,x}$ as in (4.11). The right hand side also splits
into a sum of terms indexed by the fixed points of $T$:
$$
\int_{c(\cf)}(q^*\delta_t)\tilde \omega \ = \ \sum_{x\in X^{T}}\ \frac {
\text{value of the cycle $c(\cf)$ at $(t,x)$}} {\tsize \prod_{\alpha \in {\Phi}^{
+}} (1-e^{-\alpha_x})(t )} \ .
\tag6.11b
$$ Here, as in section 4, we regard the cycle $c(\cf)$ as a section of the local
system $\Bbb C_\lambda$ via Poincar\'e duality. The value of the cycle is 
actually a number, since the local system $\Bbb C_\lambda$ is a subsheaf of
the sheaf of functions. At this point, the proposition comes
down to the assertion that the  two Lefschetz numbers 
$$
\aligned &\tsize\sum (-1)^i \,\operatorname{tr}(\,\phi_t\,:\,\Cal
H_{N'(t,x)}^i(\Bbb D
\cf(x))_x \,
\to
\, 
\Cal H_{N'(t,x)}^i(\Bbb D \cf(x))_x\otimes\Bbb C_\lambda\,)
\\ &\tsize\sum (-1)^i \,\operatorname{tr}(\,\phi_t:\Cal H_{\{x\}}^i((\Bbb D
\cf(x))|_{N''(t,x)})
\to 
\Cal H_{\{x\}}^i((\Bbb D \cf(x))|_{N''(t,x)})\otimes\Bbb C_\lambda\,)\,.
\endaligned
\tag6.12 
$$ 
coincide  and  equal the value of the cycle $c(\cf)$ at $(t,x)$.  The
verification of this assertion is our remaining task in this section. 

We use the notation established in (5.13-18). In particular, $N^+(t,x)\times
H$  is the inverse image of $N^+(t,x)$ in $\hat X$. The inclusion
$\{t\}\hookrightarrow \GR$\,, corresponding to a particular $t\in T_\Bbb R'$\,,
determines a Cartesian square
$$
\CD
\{t\}\times H\times ( N^+(t,x)\times H) @>>>\pro   
\\@AAA @AAA                                   
\\ \{t\}\times  \{t\} \times(\{x\}\times H)@>>>
\ \ \ \ \ \ \ \ \ {\hat G}_{\Bbb R}   \ \ \ \ \ \  \ \ \,;
\endCD
\tag6.13
$$  here we are using the formula (5.18) for the $T\times H$-action on
$N^+(t,x)\times H$. Combining this with (5.9) gives the commutative diagram

$$
\CD
\{t\}\times H\times ( N^+(t,x)\times H) @>>>\pro     @>>>           \hat X
\times
\hat X 
\\@AAA @AAA                                 @AA{\Delta}A   
\\ \{t\}\times  \{t\} \times(\{x\}\times H)@>>>{\hat G}_{\Bbb R}      @>>>     
\hat X
\endCD
\tag6.14
$$  involving three Cartesian squares: the right square, the left square, and
the square formed by the four terms on the perimeter -- the \lq\lq outer
square" for short. The outer square is formally analogous to the right square.
We can therefore apply  the fixed point formalism (5.10). The result is a map
\footnote{Recall our comment in section 5, to the effect that we no longer
make a notational distinction between $\cg\in \operatorname
D_\GR(X)_\lambda$ and its image in $\db {\hat X}$, as we had earlier in
section 5.}
$$
\operatorname {Hom}(\cg,\cg) \ \longrightarrow \
\Hom^0_{\hat t}(A_t,p^*\Bbb D_{\hat X}\otimes\Bbb C_\lambda) \ \cong \
\ohf {0}(\{(t,x)\}, \Bbb C_{-\lambda})\,,
\tag6.15a
$$ with the shorthand notation 
$$ A_t = \{t\}\times H\times ( N^+(t,x)\times H)\,, \qquad
\hat t = \{t\}\times  \{t\} \times(\{x\}\times H)\,.
\tag6.15b
$$ To keep the notation simple, we are using the symbols $a,p$ also for the
action and projection morphisms in the outer square. 

We claim: when we combine (6.15) with its precursor in (5.10), the resulting
diagram, 
$$
\CD
\operatorname {Hom}(\cg,\cg) @>>> \Hom^0_{\hat t}(A_t,p^*\Bbb
D_{\hat X}\otimes\Bbb C_\lambda) 
\\ @| @AAA
\\
\operatorname {Hom}(\cg,\cg) @>>> \ \ \Hom^0_{{\hat G}_{\Bbb
R}}(\pro , p^*\Bbb D_{\hat X}
\otimes
\Bbb C_\lambda) \ ,
\endCD
\tag6.16
$$ commutes. The reason is simply the functoriality of the fixed point
formalism, which follows from the functorial behavior of local cohomology.
We continue the bottom row in (6.16)  as in (5.10), then restrict to the open
subset
$\tilde\GR'$ and apply Poincar\'e duality on this open subset, which consists
of smooth points:
$$
\aligned
 \Hom^0_{{\hat G}_{\Bbb R}}(\pro , p^*\Bbb D_{\hat X}
\otimes
\Bbb C_\lambda)  @>{\ \sim\ }>> 
 \Hom^{inf}_{d+r}( {\hat G}_{\Bbb R} ,
\Bbb C_{-\lambda})  @>{\ \sim \ }>>
\\ 
\Hom^{inf}_d( {\tilde G}_{\Bbb R} , \Bbb C_{-\lambda}) @>>>
\Hom^{inf}_d( {\tilde G}'_{\Bbb R} , \Bbb C_{-\lambda})  @>{\ \sim \ }>> \oh^0 (
{\tilde G}'_{\Bbb R} , \Bbb C_{\lambda})\,.
\endaligned
\tag6.17a
$$ We analogously continue the top row, 
$$
\Hom^0_{\hat t}(A_t,p^*\Bbb D_{\hat X}\otimes\Bbb C_\lambda) @>{\ \sim \
}>> 
\ohf {0}(\{(t,x)\}, \Bbb C_{-\lambda}) @>{\ \sim \ }>>
\oh ^{0}(\{(t,x)\}, \Bbb C_{\lambda})\,.
\tag6.17b
$$ Putting (6.17a,b) together, we obtain the diagram
$$
\CD
 \Hom^0_{\hat t}(A_t,p^*\Bbb D_{\hat X}\otimes\Bbb C_\lambda) @>{\sim}>> 
\oh ^{0}(\{(t,x)\}, \Bbb C_{\lambda})
\\  @AAA @AAA
\\
 \Hom^0_{{\hat G}_{\Bbb R}}(\pro , p^*\Bbb D_{\hat X}
\otimes
\Bbb C_\lambda) @>{\sim}>> \oh^0 ( {\tilde G}'_{\Bbb R} , \Bbb C_{\lambda})
\,,
\endCD
\tag6.18
$$ in which the second vertical arrow  is evaluation of sections at
$(t,x)$. This diagram commutes because of the functoriality of the various
ingredients. Letting $\Bbb D\cf$ play the role of $\cg$, we conclude that the
value of the cycle $c(\cf)$ at $(t,x)$ is given by the image of the identity
morphism $1:\Bbb D \cf\to\Bbb D\cf$ under the chain of morphisms 
$$
\gathered
\operatorname {Hom}(\Bbb D\cf,\Bbb D\cf) @>>> \Hom^0_{\hat
t}(A_t,p^*\Bbb D_{\hat X}\otimes\Bbb C_\lambda)
\\
  @>{\sim}>>  \oh
^{0}(\{(t,x)\},
\Bbb C_{\lambda})   @>{\sim}>> \Bbb C\,.
\endgathered
\tag6.19
$$
It remains to be shown that this image of the identity coincides with the
Lefschetz number (6.12). The final isomorphism in (6.19), \lq\lq
evaluation at $t$",  makes the image of the identity morphism  a
specific number; here we are using the definition of $\Bbb C_\lambda$
as a subsheaf of the sheaf of functions. 

The morphisms in (6.19) were obtained by applying the fixed point formalism
(5.10) to the outer square in (6.14). Since the image of $\{t\}\times H\times (
N^+(t,x)\times H)$ lies in the open subset $( N^+(t,x)\times H)\times (
N^+(t,x)\times H)$ of $\hat X\times\hat X$, we may as well replace $\hat X$
by
$ N^+(t,x)\times H$ and $\cf$ by its restriction to $ N^+(t,x)\times H$.
Then, when we apply the fixed point formalism in the Cartesian square
$$
\CD
\{t\}\times H\times ( N^+(t,x)\times H)     @>>>           ( N^+(t,x)\times
H)\times ( N^+(t,x)\times H)
\\@AAA                                 @AA{\Delta}A   
\\ \{t\}\times  \{t\} \times(\{x\}\times H)      @>>>   \ \ \ \ \ \ \   
N^+(t,x)\times H \ \ \ \ \ \ \ ,
\endCD
\tag6.20
$$  we obtain the same maps as in (6.19). Next we restrict the twisted sheaf
$\cg = \Bbb D  \cf$ from $ N^+(t,x)\times H$ to $ N^+(t,x)\times \{e\}\cong
 N^+(t,x)$, as in (5.19b), resulting in the sheaf $\cg(x)\in \operatorname
D_{T_\Bbb R}( N^+(t,x))$.  Correspondingly, we take a slice of the diagram
(6.20) by replacing $ N^+(t,x)\times H$ with $  N^+(t,x)\times \{e\}\cong
 N^+(t,x)$ and $\{t\}\times H$ by $\{t\}\times \{t\}$. In this slice, the action
and projection maps induce 
$a,p : N^+(t,x)
\to N^+(t,x)$, with $a = \text{left translation by $t$}$, $p=\text{identity}$.
We use these maps to construct the Cartesian square
$$
\CD
 N^+(t,x)   @>>>            N^+(t,x)\times N^+(t,x)
\\@AAA                                 @AA{\Delta}A   
\\ \{x\}      @>>>   \ \ \ \ \ \ \    N^+(t,x) \ \ \ \ \ \ \ ,
\endCD
\tag6.21
$$ 
which maps into the Cartesian square (6.20) by inclusion. Recall the
construction of the isomorphism (5.21) by means of the commutative
diagram (5.20). We use the  notation  $\phi:t^*\cg(x) \to
\cg(x)$ for the induced morphism when we
restrict from $T_\Bbb R$ to $\{t\}$. The sheaf $\Bbb C_\lambda$
disappears at this point since $\Bbb C_\lambda|_{\{t\}} \cong \Bbb C$ by
evaluation at $t$. When we apply the fixed point formalism to the square
(6.21) instead of (5.9), $\phi$ induces a morphism
$$
\operatorname {Hom}(\Bbb D\cf(x),\Bbb D\cf(x)) \ @>>> \ \Hom^0_{
\{x\}}(N^+(t,x),\Bbb D_{N^+(t,x)})\   @>{\ \sim\ }>>  \ \Bbb C\,,
\tag6.22
$$ 
analogous to (6.19). 

We claim: The image of $1\in \operatorname
{Hom}(\Bbb D\cf(x),\Bbb D\cf(x))$ under the composite morphism in
(6.22) coincides with the image of the identity under the chain of
homomorphisms (6.19). The crux of the matter is the functoriality of the
fixed point formalism  with respect to non-characteristic maps
-- in our case, $N^+(t,x)
\hookrightarrow N^+(t,x)\times H$ -- a general fact which can be verified
by tracing through diagrams, though not entirely without effort. The
inclusion $N^+(t,x)
\hookrightarrow N^+(t,x)\times H$ relates the two chains of morphisms
(6.19) and (6.22), so our claim follows from the functoriality properties of
the fixed point formalism. 

To complete the proof, we  still must identify the Lefschetz numbers (6.12)
with the image of the identity under the  homomorphism (6.22). The
equality of the three quantities follows from an appropriate generalization
of the Lefschetz fixed point theorem. Generalizations that apply in our
situation have been given by Kashiwara \cite{Ka4,KSa} and
Goresky-MacPherson \cite{GM}. Kashiwara's fixed point formalism
expresses the global Lefschetz number as a sum over local contributions
corresponding to the components of the fixed point sets; in the case of
isolated fixed points -- which is the case of interest to us --  he gives an
explicit description of the local contributions. Specifically, the local
contribution corresponding to an isolated fixed point
$x$ is the image of the identity under the homomorphism (6.22);
further, Kashiwara identifies this local contribution with Lefschetz numbers
similar to (6.12), but specialized to the tangent space at $x$.
Goresky-MacPherson also express the global Lefschetz number as a sum of
local contributions. In the case of an isolated fixed point $x$, their local
contribution coincides precisely with either of the two local Lefschetz
numbers (6.12) \cite{GM, formulas for $A_4$ and $A_5$}. They also establish
the uniqueness of the local contributions, provided they are expressed in
terms of local data near the fixed point \cite{GM, \S 5.1}. In particular, their
local Lefschetz numbers coincide with Kashiwara's. That gives us the
conclusion we need. 
\enddemo
\vskip 1cm

\subheading{\bf 7. The case of the discrete series}
\vskip .5cm

In this section we establish our character formulas for standard sheaves
associated to certain $\GR$-orbits $S$ and certain twisting parameters
$\lambda$\,: orbits $S$ attached to a compact Cartan subgroup $T_\Bbb R$ --
these are necessarily open -- and any regular anti-dominant $\lambda$\,. Such
geometric data correspond to discrete series representations. This
particular case of theorems 5.12 and 5.24 was already established by
Kashiwara in his announcement \cite{Ka4}.

Let then $T_\Bbb R \subset \GR$ be a compact Cartan subgroup, and $S  = 
S(T_\Bbb R,\tau_{x_0})$ a $\GR$-orbit corresponding to the datum of $T_\Bbb
R$ and of a fixed point ${x_0}$ of $T_\Bbb R$, as in 6.5. Our hypotheses on $\GR$
imply that the compact Cartan subgroup $T_\Bbb R$ is connected, hence a
torus. As was argued in section 6, standard sheaves associated to $S$ and any
particular twisting parameter $\lambda$ correspond bijectively to
$\GR$-equivariant, irreducible,  $(-\lambda-\rho)$-monodromic local
systems on $S$, and those, in turn, correspond bijectively to characters
$\chi: T_\Bbb R \to \Bbb C^*$ whose differentials coincide with
$\lambda-\rho$. Because of the connectedness of $T_\Bbb R$, if an irreducible,
$\GR$-equivariant, $(-\lambda-\rho)$-monodromic local system on
$S$ exists, it is unique; moreover, such a local system exists precisely when
$\lambda-\rho$ is an $H$-integral weight \footnote{In the linear case, compact
Cartan subgroups of $\GR$ have the same weight lattice as the universal
Cartan
$H$; if, contrary to our assumptions, $\GR$ fails to be linear, the phrase
\lq\lq $H$-integral" should be replaced by \lq\lq $\GR$-integral".}. Thus,
without loss of generality, we assume 
$$
\lambda \, \in \, \Lambda+\rho \qquad (\,\Lambda = \text{weight lattice of
$H$\,)}\,.
\tag7.1
$$ Let $j:S\hookrightarrow X$ denote the inclusion. We can appeal to (2.12-13),
and conclude that $Rj_*\Bbb C_S\in \operatorname D_\GR(X) = \operatorname
D_\GR(X)_{-\lambda}$ is the only standard sheaf in
$\operatorname D_\GR(X)_{-\lambda}$ associated to the orbit $S$.

Recall (2.15) and the discussion below it: because of the integrality of
$\lambda +\rho$, the twisted sheaf $\Cal O(\lambda)$ becomes an actual sheaf
on $X$ with a $G$-action, and as such coincides with the sheaf of holomorphic
sections of the  $G$-equivariant line bundle $\bold L_{\lambda-\rho}$ on $X$.
Hence, with $\Cal F =  Rj_*\Bbb C_S$, the construction (2.16-17) produces the
virtual character
$$
\aligned
\Theta(Rj_*\Bbb C_S) \ &= \  \tsize \sum_p (-1)^p
\Theta(\operatorname{Ext}^p(\Bbb DRj_*\Bbb C_S,\Cal O_X(\lambda))) 
\\ &= \  \tsize \sum_p (-1)^p
\Theta(\operatorname{Ext}^p(Rj_!\Bbb C_S[2 \dim_\Bbb C X],\Cal O_X(\bold
L_{\lambda-\rho}))) 
\\  &= \  \tsize \sum_p (-1)^p
\Theta(\operatorname{Ext}^p(\Bbb C_S,j^!\Cal O_X(\bold L_{\lambda-\rho}))) 
\\   &= \  \tsize \sum_p (-1)^p
\Theta(\operatorname{Ext}^p(\Bbb C_S,j^*\Cal O_X(\bold L_{\lambda-\rho}))) 
\\    &= \  \tsize \sum_p (-1)^p
\Theta(\oh^p(S,\Cal O_X(\bold L_{\lambda-\rho}))) \,.
\endaligned
\tag7.2
$$ Here we have used the fact that an even shift in degrees does not affect the
Euler characteristic, the adjointness of $Rj_!$ and $j^!$, and the equality of
$j^!$ and $j^*$ for open embeddings.

Since $T_\Bbb R$ is a compact Cartan subgroup of the linear group $\GR$, 
$\tau_{x_0}^*\lambda$ is the full weight lattice of the torus $T_\Bbb R$, and
$\tau_{x_0}^*(\Lambda+\rho)$ the weight lattice of $T_\Bbb R$ shifted by the
half-sum of the positive roots. To each regular $\mu\in\tau_{x_0}^*(\Lambda
+\rho)$\,, Harish-Chandra
\cite{HC6} associates a discrete series representation  $\pi_\mu$\,, whose
character he denotes by the symbol $\Theta_\mu$\,; every discrete series
character is of this type for some regular  $\mu\in\tau_{x_0}^*(\Lambda
+\rho)$, and $\Theta_\mu = \Theta_\nu$ if and only if the parameters
$\mu,\nu$ are conjugate under the normalizer of $T_\Bbb R$ in $\GR$.  The
integrality condition (7.1) on the parameter $\lambda$ in (7.2) implies that
$\topsmash{2\frac {(\lambda,\alpha)} {(\alpha,\alpha)}}$ is an integer for each
$\alpha\in\Phi$. In addition to (7.1), we now impose the condition that
$\lambda$ is regular anti-dominant, in the sense that
$$ 2\frac{(\lambda,\alpha)} {(\alpha,\alpha)} \in \Bbb Z_{<0} \ \ \ \text{for all
$\alpha\in\Phi^+$}\,.
\tag7.3
$$ For every $\lambda\in\Lambda+\rho$ which satisfies this antidominance
condition, we define 
$$
\Theta(S,\lambda) \ = \ \Theta_{\tau_{x_0}^*\lambda}\,.
\tag7.4
$$ Two comments are in order. First, the open $\GR$-orbit $S$ determines the
pair
$(T_\Bbb R, {x_0})$ up to simultaneous $\GR$-conjugacy, so
$\Theta(S,\lambda)$ is canonically attached to $(S,\lambda)$. Secondly, we
insists on the anti-dominance condition (7.3) in the definition (7.4) because we
want
$\Theta(S,\lambda)$ to depend coherently on the parameter $\lambda$; see
section 8. This has the added advantage of making our parametrization of the
discrete series, in terms of pairs $(S,\lambda)$, one-to-one.

The  compact Cartan subgroup $T_\Bbb R$ lies in a maximal compact subgroup
of
$\GR$\,; we may as well assume that it is the one we had already chosen:
$T_\Bbb R \subset K_\Bbb R$\,.  Then
$K_\Bbb R/T_\Bbb R$ can be identified with the flag variety of $\fk$\,. In
particular, the (real) dimensions of $K_\Bbb R/T_\Bbb R$ and $\GR/K_\Bbb R$
are even. We set 
$$ s \ = \ \frac 1 2 \, \dim K_\Bbb R/T_\Bbb R \,, \qquad \qquad q \ = \ \frac 1 2
\, \dim G_\Bbb R/K_\Bbb R\,.
\tag7.5
$$ Note that the virtual character $\Theta(Rj_*\Bbb C_S)$ in (7.2) depends on
the choice of $\lambda$, even though $\lambda$ does not appear in the
notation: we have chosen to regard $Rj_*\Bbb C_S$ as an object in
$\operatorname D_\GR(X)_{-\lambda}\,$.

\proclaim{7.6 Theorem} Under the hypotheses just stated, $\Theta(\oh^s(S,\Cal
O_X(\bold L_{\lambda-\rho}))) = \Theta(S,\lambda)$, and $\oh^p(S,\Cal
O_X(\bold L_{\lambda-\rho})) = 0$ if $p\neq s$. In particular, $\Theta(Rj_*\Bbb
C_S) = (-1)^s
\Theta(S,\lambda)$.
\endproclaim

When $\lambda$ is not only anti-dominant but also \lq\lq very regular", i.e., if
$|(\lambda,\alpha)| \gg 0$ for all $\alpha \in\Phi$, this is the main result of
\cite{S1} in conjunction with the characterization of discrete series
representations by lowest $K$-type \cite{S2,V}. The general case
then follows from the Jantzen-Zuckerman translation principle, as in
\cite{S3}, for example. Alternatively, the theorem is a consequence of the main
result in \cite{KSd}.

In the remainder of this section, we shall verify the two character formulas in
the present situation, i.e., with $\operatorname{rk} \GR =
\operatorname{rk} K_\Bbb R$, $S\subset X$ open,
$\cf =  Rj_*\Bbb C_S$, and $\lambda$ satisfying the anti-dominance condition
(7.3). We first deal with the case of the integral formula (3.8).

The decomposition $\fg = \ft \oplus [\ft,\fg]$ permits us to view
$\tau_{x_0}^*\lambda$ as an element of $\fg^*$ -- an element of $i\fg_\Bbb
R^*$, in fact, since differentials of characters of $T_\Bbb R$ are purely
imaginary. Recall that $\Omega_\lambda$ denotes the $G$-orbit of
$\tau_{x_0}^*\lambda$ in $\fg^*$. We define
$$
\Omega(S,\lambda) \ = \ \text{$\GR$-orbit of $\tau_{x_0}^*\lambda$ in
$i\fg_\Bbb R^*$}\,.
\tag7.7
$$ As in the case of $\Theta(S,\lambda)$, the notation is justified because $S$
determines $T_\Bbb R$ and $\tau_{x_0}$ up to simultaneous $\GR$-conjugacy.
Let us note that the union of the $\Omega(S,\lambda)$, corresponding to all the
open
$\GR$-orbits in $X$, coincides with $\Omega_\lambda \cap i\fg_\Bbb R^*$.
Recall that the complex coadjoint orbit $\Omega_\lambda$ comes equipped
with a canonical algebraic symplectic form $\sigma_\lambda$. Since
$\lambda$ is regular, the orbit $\Omega_\lambda$ has complex dimension $2n$,
where
$n = \dim_\Bbb C X$, as before. Also recall our convention (3.6) -- without the
choice of
$i=\sqrt{-1}$ -- for the Fourier transform $\hat\phi$.  The following result is
due to Rossmann
\cite{R1}. 

\proclaim{7.8 Theorem} The two form $\,-i\,\sigma_\lambda$ restricts to a
real, non-degenerate form on the real submanifold
$\Omega(S,\lambda)\subset\Omega_\lambda$. Orient $\Omega(S,\lambda)$ so
as to make the top exterior power of $\,-i\,\sigma_\lambda$ positive. Let
$\theta(S,\lambda)$ be the character on the Lie algebra corresponding to 
the discrete series character
$\Theta(S,\lambda)$. Then, for $\phi\in C_c^\infty(\gr)$, 
$$
\int_\gr \, \theta(S,\lambda)\,\phi\, dx \ = \ \frac 1 {(2\pi i)^n
n!}\,\int_{\Omega(S,\lambda)}\, \hat\phi\,
\sigma_\lambda^n\ ;
$$ the integral on the right converges absolutely.
\endproclaim
 
The definition of $\hat \phi$ involves the choice of Euclidean measure $dx$ on
$\gr$; we are using the same measure on the left of the above identity, of
course. 

We shall use the open embedding theorem of \cite{SV4, theorem 4.2} to relate
the inverse image of $\Omega(S,\lambda)$ under the twisted moment map
$\mu_\lambda$ to the characteristic cycle of the sheaf $Rj_*\bc_S$\,. The crux
of the matter is the following technical result. The statement involves a real
valued function
$f_\lambda$ on $S$ which we now define. The parameter $\lambda$ lies in
$\Lambda +\rho$, and $\rho$ is integral, if not for $G$, then at least for some
2-fold covering  of the complex group $G$. We conclude that  the  line bundle
$\bold L_\lambda$ exists as an algebraic line bundle on $X$, equivariant with
respect to $G$ or a 2-fold covering of $G$. By assumption, $\GR$ acts on $S$
with compact isotropy groups. Replacing $\GR$ by its inverse image in the
2-fold covering of $G$ if necessary, we find that
$\bold L_\lambda|_{S}$ admits a $G_\Bbb R$-invariant Hermitian metric
$\gamma_{nc}$. Analogously, there exists a Hermitian metric $\gamma_c$ on
$\bold L_\lambda$ invariant  under $U_\Bbb R$ or its inverse image in the
2-fold covering. As the quotient of two real algebraic  metrics on a line
bundle, 
$$ f_\lambda \ = \  \frac {\gamma_{nc}}{\gamma_c}
\tag7.9
$$ is a positive, real algebraic function on $S$. In the case of $\lambda =
-2\rho$, this function was first used in \cite{S1} to study the Dolbeault
cohomology of $\GR$-equivariant line bundles. 

\proclaim{7.10 Lemma} If $\lambda\in\Lambda+\rho$ satisfies the
anti-dominance condition (7.3), some positive integral power $f_\lambda^m$
of the function
$f_\lambda$ extends real algebraically to all of
$X$, and this extension vanishes on the boundary of $S$. The image
$\mu_\lambda(d\log f_\lambda)$ of $d\log f_\lambda$, viewed as a
submanifold of $T^*S$, coincides with $\Omega(S,\lambda)$. More precisely, 
$$
\mu_\lambda\ : \ d\log f_\lambda \ \longrightarrow  \ \Omega(S,\lambda)
$$ is a real algebraic isomorphism which preserves or reverses orientation,
depending on whether the integer $s$ is even or odd; here we orient 
$\Omega(S,\lambda)$ as in (7.8) and $d\log f_\lambda \cong S$ via the
complex structure of $X$.
\endproclaim

Before embarking on the proof -- which is lengthy, though not difficult -- we
show how to deduce theorem 3.8 in the current setting. We apply lemma 3.19
with $C_1 =
 d\log f_\lambda$, with
$C_2 = \occ(Rj_*\Bbb C_S)$,  and with
$$
\tilde C \ = \ \text{image of} \  \ (0\,,1)\times S \to \ct\ , \ \ \ (t,x)\mapsto
t\,d\log f_\lambda(x)\,.
\tag7.11
$$ 
We orient $C_1$  via $|C_1|=d\log f_\lambda \cong S$ and the complex
structure on $S$, and orient $\tilde C$ by means of $|\tilde C|\cong (0,1)\times
S$. Since $f_\lambda^m$ is real algebraic, so is $d\log f_\lambda$. Thus $C_1$
is a real algebraic cycle, and $\tilde C$ is a semi-algebraic chain.  The support
$|C_1|$ coincides with $\mu_\lambda^{-1}(\Omega(S,\lambda))$, so the values
of $\mu_\lambda$ on $|C_1|$ lie in $\Omega(S,\lambda)\subset i\frak g_\Bbb
R^*$. It follows that the real part of $\mu$ remains bounded on $|C_1|$ and
$|\tilde C|$. The open embedding theorem 4.2 of \cite{SV4}, coupled with
\cite{SV4, proposition 3.25}, implies 
$$
\partial \tilde C \ = \ d\log f_\lambda \ - \ \occ(Rj_*\Bbb C_S) \ = \ C_1 -
C_2\,.
\tag7.12
$$ 
We have verified the hypotheses of lemma 3.19 in the present setting, hence
$$
\int_{C_1}\mu_\lambda^*\hat\phi\,(-\sigma +\pi^*\tau_\lambda)^n =
\int_{C_2}\mu_\lambda^*\hat\phi\,(-\sigma +\pi^*\tau_\lambda)^n\,.
\tag7.13
$$
Since $C_1 = (-1)^s \,\mu_\lambda^{-1}\Omega(S,\lambda)$ by lemma
7.10,  taking into account proposition 3.3, we find
$$
\aligned
\int_{\Omega(S,\lambda)} \hat\phi\,\sigma_\lambda^n\ &= \
\int_{\mu_\lambda^{-1}\Omega(S,\lambda)}\mu_\lambda^*\hat\phi\, (-\sigma
+
\pi^*\tau_\lambda)^n \ 
\\
\vspace{2.5\jot} &= \ (-1)^s\int_{\occ(Rj_*\Bbb C_S)}\mu_\lambda^*\hat\phi\,
(-\sigma +
\pi^*\tau_\lambda)^n\,,
\endaligned
\tag7.14
$$ for every test function $\phi\in C_c^\infty(\gr)$. We combine this with
theorems 7.6 and 7.8, and conclude theorem 3.8 in the case of the discrete
series:
\proclaim{7.15 Proposition} Under the hypotheses stated at the beginning of
this section, the character formula (3.8) holds.
\endproclaim

\demo{Proof of lemma 7.10} The function $f_\lambda$ depends
multiplicatively on the parameter $\lambda$\,; in particular, $f_\lambda^m =
f_{m\lambda}$\,. Thus, replacing $\lambda$ by an appropriate positive integral
multiple $m\lambda$, we may as well assume that the parameter
$\lambda$ is integral even with respect to any particular finite quotient of
$G$. This observation allows us to  replace $G$ by $G/Z_{[G,G]}$ \ (\,
$Z_{[G,G]}=$ \,center of the commutator subgroup of $G$\,). Because of (7.3) we
can define
$$ V_\lambda \ = \ \text{irreducible $G$-module with lowest weight
$\lambda$}\,.
\tag7.16
$$ For $x\in X$, let $\frak n_x = [\fb_x,\fb_x]$ denote the unipotent radical of
the isotropy algebra at $x$. Because of our convention for ordering the roots,
the space of $\fn_x$-invariants $V_\lambda^{\fn_x}$ is the lowest weight
space for the action of any concrete Cartan subgroup $\fc\subset\fb_x$\,,
$\fc\cong
\fb_x/\fn_x \cong \fh$. As $x$ varies over $X$, the one dimensional subspaces
$V_\lambda^{\fn_x}\subset V_\lambda$ constitute a line bundle, and 
$$
\bold L_\lambda \ = \ \text{line bundle with fiber $V_\lambda^{\fn_x}$ at
$x$}\,.
\tag7.17
$$ In fact, $\Cal O(\bold L_\lambda)$ is the Beilinson-Bernstein localization of
$V_\lambda$ at $\lambda$. Since $U_\Bbb R$ is compact, there exists an
essentially unique $U_\Bbb R$-invariant metric $h_c$ on $V_\lambda$\,.
Renormalizing this metric, if necessary, we can make the identification
$$
\gamma_c \ = \ \text{restriction of $h_c$ to the lines $V_\lambda^{\fn_x}$}
\tag7.18
$$ via (7.17). 

We need an analogous description of $\gamma_{nc}$\,. Let $\theta: \fg \to \fg$
denote the Cartan involution corresponding to the maximal compact subgroup
$K_\Bbb R\subset G_\Bbb R$, and $\frak p_\Bbb R$ the $(-1)$-eigenspace of
$\theta$ in $\gr$\,. Then 
$$
\GR \ = \ \exp(\frak p_\Bbb R)\cdot K_\Bbb R \qquad \text{(Cartan
decomposition)}.
\tag7.19
$$ Since $\KR$ contains the compact Cartan subgroup $T_\Bbb R$, the Cartan
involution is inner; see, for example, \cite{Hel}. Specifically, there exists
$$ r_\theta\in T_\Bbb R\cap [G,G] \,, \ \ \text{such that} \ \ \theta =
\operatorname{Ad}r_\theta\,.
\tag7.20a
$$ Since we had made $[G,G]$  center free, this identity uniquely determines
$r_\theta$\,. But $\theta^{-1} = \theta$, hence
$$ r_\theta \ = \ r_\theta^{-1}
\tag7.20b
$$
 We claim: the identity
$$ h_{nc}(u,v) \ = \ h_c(u,r_\theta v) \qquad (\,u,\,v\in V_\lambda\,)
\tag7.21
$$ defines a $\GR$-invariant indefinite Hermitian form. Indeed, 
$$ h_{nc}(v,u) \ = \ h_c(v,r_\theta u) \ = \ h_c(r_\theta v, u) \ = \ 
\overline {h_c(u,r_\theta v)} \ = \ \overline {h_{nc}(u, v)}
$$ because $r_\theta  =  r_\theta^{-1}\in T_\Bbb R \subset U_\Bbb R$\,. Given
any $g\in\GR$, we write $g = pk$ with $p\in\exp(\frak p_\Bbb R)$ and $k\in
K_\Bbb R\,$. The Lie algebra $\frak u_\Bbb R$ of $U_\Bbb R$  contains $i\frak
p_\Bbb R$ -- in fact, $\frak u_\Bbb R = \kr\oplus i\frak p_\Bbb R$ \cite{Hel}
-- hence both $\frak p_\Bbb R$ and $\exp \frak p_\Bbb R$  operate on
$V_\lambda$ by symmetric operators, relative to $h_c$\,. Thus
$$
\aligned h_{nc}(gu,gv) \ &= \ h_c(pku,r_\theta pkv) \ = \ h_c(pku,r_\theta
pr_\theta^{-1}r_\theta kv)  \ = \ h_c(pku, p^{-1}r_\theta kv)
\\ &= \ h_c(ku,pp^{-1}r_\theta kv) \ = \ h_c(u,k^{-1}r_\theta
kr_\theta^{-1}r_\theta v)  \ = \ h_c(u,k^{-1} kr_\theta v)
\\ &=  \ h_c(u,r_\theta v) \ = \ h_{nc}(u,v)\,.
\endaligned
$$ This establishes our claim. The Cartan involution preserves the root space
decomposition 
$$
\fg \ = \ \ft\, \oplus \,(\,\oplus_{\alpha \in\Phi} \fg^\alpha\,)
$$ and each root space is contained either in $\fk$ or in $\frak p = \fk^\perp$,
hence $e^\alpha(r_\theta) = \pm 1$. Since $G$ is center free, the weights of
$V_\lambda$ lie in the root lattice. We conclude that $r_\theta$ acts as
multiplication by $\pm 1$ on each $\ft$-weight space in $V_\lambda$\,. In
particular, the $\GR$-invariant indefinite Hermitian form $h_{nc}$ is  either
strictly positive or strictly negative on all the lines $V_\lambda^{\fn_x}$,
$x\in S$. Renormalizing $h_{nc}$ by an appropriate positive or negative factor,
we get the description 
$$
\gamma_{nc} \ = \ \text{restriction of $h_{nc}$ to the lines
$V_\lambda^{\fn_x}$, $x\in S$}
\tag7.22
$$ of  $\gamma_{nc}$ analogous to (7.18). 

The lines $V_\lambda^{\fn_x}$  vary algebraically with $x$, and $h_c$ is
positive definite for every $x\in X$. It follows that the ratio $h_{nc}/h_c$ of
the two metrics on the lines $V_\lambda^{\fn_x}$ is a globally defined, real
algebraic function on $X$, which agrees with $f_\lambda$ on $S$. Now let
$x$ be a point in the boundary of $S$; we claim:
$$ h_{nc} \ \ \text{vanishes on} \ \ V_\lambda^{\fn_x}\,.
\tag7.23
$$ Since $h_{nc}$ is $\GR$-invariant, we may replace $x$ by any
$\GR$-translate. In this way we can arrange that the point $x$ is fixed by a
$\theta$-stable Cartan subalgebra $\fc_\Bbb R$ of $\gr$ \cite{Ma}.
The corresponding Cartan subgroup $C_\Bbb R\subset \GR$ is non-compact
since otherwise $x$ would have to lie in an open $\GR$-orbit, and consequently
$\theta$ does not act as the identity on the complexified Cartan subalgebra
$\fc\subset\fg$. Both $V_\lambda^{\fn_x}$ and $r_\theta V_\lambda^{\fn_x}$
are weight spaces of
$\fc$, and the corresponding weights are necessarily distinct because the
lowest weight of $V_\lambda$ is regular, and $\theta =
\operatorname{Ad}r_\theta$ normalizes $\fc$ and acts as a non-trivial
element of the Weyl group, as we have just seen. On the other hand, as a
$\theta$-stable,
$\gr$-real Cartan subalgebra, $\fc$ is also $\frak u_\Bbb R$-real. This makes
the $\fc$-weight space decomposition of $V_\lambda$ perpendicular with
respect to $h_c$. The $h_c$-perpendicularity of the two lines
$V_\lambda^{\fn_x}$ and
$r_\theta V_\lambda^{\fn_x}$ implies (7.23), which in turn implies the
vanishing at $x$ of the natural extension of $f_\lambda$ to $X$. This
establishes the first part of the lemma.

Both $S$ and $\Omega(S,\lambda)$ are homogenous spaces for $\GR$, the
former with isotropy subgroup $T_\Bbb R$ at ${x_0}$, the latter with isotropy
subgroup
$T_\Bbb R$ at $\tau^*_{x_0}\lambda$\,. To see this in the case of $S$, we recall
(6.6a) and note that $\fn_{x_0}\cap \bar\fn_{x_0}=0$ ($\bar\fn_{x_0}=$
complex conjugate of
$\fn_{x_0}$ with respect to $\gr$) since all roots are imaginary on the Lie
algebra  of the compact Cartan $T_\Bbb R$. In the case of
$\Omega(S,\lambda)$, we appeal to the regularity of $\tau^*_{x_0}\lambda$
and the equality $T_\Bbb R =
\GR\cap T$. We conclude that
$$ I \ : \ S\  @>{\ \sim\ }>> \ \Omega(S,\lambda)\,, \qquad I(g{x_0}) =
\operatorname{Ad} g\,(\tau^*_{x_0}\lambda) \ \ \text{if $g\in\GR$}\,,
\tag7.24
$$ is a diffeomorphism. To calculate  $\mu_\lambda^{-1}\circ I$, we write
$$
\operatorname{Ad} g\,(\tau^*_{x_0}\lambda) \ = \ I(g{x_0}) \ = \
\mu_\lambda(\tilde x ,
\xi)
\qquad \text{with $\tilde x\in X$ and
$\xi\in T^*_{\tilde x} X$}\,.
\tag7.25a
$$ We identify $\fh^* \cong (\fb_{\tilde x}/\fn_{\tilde x})^*\cong \fn_{\tilde
x}^\perp/\fb_{\tilde x}^\perp\,$. Then, by definition of $\mu_\lambda$,
$$
\mu_\lambda(\tilde x ,\xi) \,\in\,\fn_{\tilde x}^\perp\,, \qquad
\mu_\lambda(\tilde x ,
\xi) \equiv \lambda \ \text{mod} \,\fb_{\tilde x}^\perp\,.
\tag7.25b
$$ On the other hand, 
$$
\operatorname{Ad} g\,(\tau^*_{x_0}\lambda) \,\in\,\fn_{g{x_0}}^\perp\,,
\qquad
\operatorname{Ad} g\,(\tau^*_{x_0}\lambda) \equiv \lambda \ \text{mod}
\,\fb_{g {x_0}}^\perp\,.
\tag7.25c
$$ Since $\lambda$ is regular, (7.25a-c) imply $g{x_0} = \tilde x$. We choose
$u\in U_\Bbb R$ so that 
$$ g{x_0} = \tilde x = u{x_0}\,.
\tag7.26a
$$ The definition of $\mu_\lambda$ gives the identity
$$
\mu_\lambda(\tilde x,\xi) \ = \ \operatorname{Ad} u\,(\tau_{x_0}^*\lambda) +
\mu(\xi)\,.
\tag7.26b
$$ In effect, we have calculated $\mu_\lambda^{-1} \circ I$\,:
$$ (\mu_\lambda^{-1} \circ I)(g{x_0}) \ = \ (g{x_0}, \xi) \,,  \ \  \text{with}\ \
\mu(\xi)
\ =
\ \operatorname{Ad} g\,(\tau_{x_0}^*\lambda) - \operatorname{Ad}
u\,(\tau_{x_0}^*\lambda)\,;
\tag7.26c
$$ here $u$ is determined by (7.26a), uniquely up to right multiplication by an
element of $T_\Bbb R$\,. In particular, when we compose
$\mu_\lambda^{-1}\circ I$ with the projection $\pi:\ct \to X$ we obtain the
identity on $S$. At this point, the remaining assertions of the lemma come
down to 
$$
\mu(d\log f_\lambda|_{g{x_0}}) \ = \ \ \operatorname{Ad}
g\,(\tau_{x_0}^*\lambda) -
\operatorname{Ad} u\,(\tau_{x_0}^*\lambda)\,,
\tag7.27a
$$ whenever $g\in\GR$, $u\in U_\Bbb R$ and $g{x_0}=u{x_0}$ as before,  and
$$ I \ : \ S \ @>{\ \sim \ }>> \Omega(S,\lambda) \ \ \text{is orientation
preserving if and only if $s$ is even}
\tag7.27b
$$ when $S$ and $\Omega(S,\lambda)$ are oriented as in the statement of the
lemma. 

We deal with (7.27a) first. In the statement of the lemma, as in (7.27a), we
regard the differential $d\log f_\lambda$ as a real algebraic section of the
 cotangent bundle; this involves an explicit isomorphism $\ct
\cong T^*X^\Bbb R$, which we normalize as in \cite{KSa, Chapter 11}. Hence, by
the definition of the moment map,
$$
\frac d {dt}
\log f_\lambda(\exp(tZ)g{x_0})|_{t=0} \  = \  
2\,\operatorname{Re}\langle\,\mu(d\log f_\lambda|_{g{x_0}})\,,\,Z\,\rangle\,, 
\tag7.28
$$ for every  $Z\in\fg$\,. Let $\ell$ be a generator of
$V_\lambda^{\fn_{g{x_0}}}$. Then, by (7.18) and (7.22), 
$$
\aligned &\frac d {dt}
\log f_\lambda(\exp(tZ)g{x_0})|_{t=0}\ = 
\\ &\ \ \  \ \frac{h_{nc}(Z\ell,\ell) + h_{nc}(\ell,Z\ell)}{h_{nc}(\ell,\ell)} \ -
\ \frac{h_{c}(Z\ell,\ell) + h_{c}(\ell,Z\ell)}{h_{c}(\ell,\ell)}\,.
\endaligned 
\tag7.29
$$ Both summands on the right are homogenous of degree zero in $\ell$, so we
may use two different generators of $V_\lambda^{\fn_{g{x_0}}}$ in the two
summands. Let us use $\ell = g\ell_0$ in the first instance, with $\ell_0\in
V_\lambda^{\fn_{{x_0}}}$, and $\ell = u\ell_0$ in the second instance:$$
\aligned &\frac d {dt}
\log f_\lambda(\exp(tZ)g{x_0})|_{t=0}\ = 
\\ &\ \ \ \frac{h_{nc}(Zg\ell_0,g\ell_0) +
h_{nc}(g\ell_0,Zg\ell_0)}{h_{nc}(\ell_0,\ell_0)} \ -
\ \frac{h_{c}(Zu\ell_0,u\ell_0) + h_{c}(u\ell_0,Zu\ell_0)}{h_{c}(\ell_0,\ell_0)}\,.
\endaligned 
\tag7.30
$$ Note that $h_{nc}(Zg\ell_0,g\ell_0) =
h_{nc}((\operatorname{Ad}g^{-1}Z)\ell_0,\ell_0)$ depends only on the
$\ft$-component of $\operatorname{Ad}g^{-1}Z$ in the decomposition 
$\fg = \ft \oplus [\ft,\fg]$, since the $\ft$-weight spaces in $V_\lambda$ are
$h_{nc}$-perpendicular. On the other hand, $Y\ell_0 =
\langle\,\tau_{x_0}^*\lambda\,,\, Y\,\rangle\,\ell_0$ if $Y\in\ft$\,. Thus
$$
\frac{h_{nc}(Zg\ell_0,g\ell_0)}{h_{nc}(\ell_0,\ell_0)} \ = \
\frac{h_{nc}((\operatorname{Ad}g^{-1}Z)\ell_0,\ell_0)}{h_{nc}(\ell_0,\ell_0)} \
= \ \langle\,\tau_{x_0}^*\lambda\,,\, (\operatorname{Ad}g^{-1}Z)\,\rangle\,.
$$ 
Arguing similarly in the case of the other terms on the right of (7.30), we
find
$$
\aligned
\frac d {dt} \log f_\lambda(\exp(tZ)g{x_0})|_{t=0}\ &= \
2\,\operatorname{Re}\,\langle\,\tau_{x_0}^*\lambda\,,\,
(\operatorname{Ad}g^{-1}Z) - (\operatorname{Ad}u^{-1}Z)\,\rangle
\\ &=\  2\,\operatorname{Re}\,\langle\,\operatorname{Ad}
g\,(\tau_{x_0}^*\lambda) -
\operatorname{Ad} u\,(\tau_{x_0}^*\lambda)\,,\,Z\,\rangle\,.
\endaligned
\tag7.31
$$ Taken together, (7.28) and (7.31) imply (7.27a).

We verify the orientation statement (7.27b) by reducing it to the  special
cases of $\GR = SU(1,1)$ and $\GR=SU(2)$. To begin with, we need to establish
the statement only at a single point -- specifically, at the fixed point
${x_0}\in S$ for
$T_\Bbb R$. We identify the holomorphic tangent space $T_\lambda
\Omega_\lambda$ of the complex orbit at $\lambda$ with 
$$
\fg/\ft \ \cong \ [\ft,\fg] \ = \ \ \underset {\alpha\in\Phi(\fg,\ft)} \to \oplus
\ \fg^\alpha\,.
\tag7.32a
$$ There are analogous descriptions for the real tangent space $T_\lambda
\Omega(S,\lambda)$, 
$$
\gr/\ft_\Bbb R \ \cong\ T_\lambda \Omega(S,\lambda) \ \hookrightarrow \
\Bbb C
\otimes_\Bbb R T_\lambda \Omega(S,\lambda)\ \cong \ \fg/\ft\,,
\tag7.32b
$$ and of the real tangent space $T_{x_0} S^\Bbb R$ of the open orbit
$S\subset X$ considered as a real manifold,
$$
\gr/\ft_\Bbb R \ \cong\ T_{x_0} S^\Bbb R \ \hookrightarrow \ \Bbb C
\otimes_\Bbb R T_{x_0} S^\Bbb R\ \cong \ \fg/\ft\,.
\tag7.32c
$$ Via these identifications, the differential of the map $I$ at ${x_0}$ becomes
the identity, and the symplectic form $\sigma_\lambda$ is given by the
formula
$$
\sigma_\lambda(Z_1,Z_2) \ = \ \langle\,\lambda\,,\,[Z_1,Z_2]\,\rangle\,,
\qquad Z_1,Z_2 \in [\ft,\fg]\,.
\tag7.33
$$ Note that $\langle\,\lambda\,,\,[\ft,\fg]\,\rangle = 0 $, hence
$$
\sigma_\lambda(\fg^\alpha,\fg^\beta) \ = \ 0 \ \ \ \text{unless} \ \alpha +\beta
= 0\,.
\tag7.34
$$ The isomorphism $\tau_{x_0}^*:\Phi \cong \Phi(\fg,\ft)$ induces a positive
root system $\Phi^+(\fg,\ft)$. Each root $\alpha\in\Phi(\fg,\ft)$ is either
compact or noncompact, in the sense that the subalgebra
$$
\fg_\alpha \ = \ \fg^\alpha \oplus \fg ^{-\alpha} \oplus [\fg^\alpha,
\fg ^{-\alpha}] \ \cong \ \frak s\frak l(2,\Bbb C)
\tag7.35
$$ of $\fg$ intersects $\gr$ in a copy of either $\frak s \frak u (2)$ or $\frak s
\frak u (1,1)$. Let $G_\alpha$ denote the connected subgroup of $G$ with Lie
algebra $\fg_\alpha$, and $G_{\Bbb R,\alpha}$ the real form
$\GR\cap G_\alpha$\,. Then $G_{\Bbb R,\alpha}\cdot {x_0}$\,, the $G_{\Bbb
R,\alpha}$-orbit of ${x_0}$ in $S$, is isomorphic as homogenous space to $\Bbb
P^1$ if $\alpha$ is compact, or to \,$\Delta$\,=\,unit disc if
$\alpha$ is noncompact. The three spaces $\Omega_\lambda$\,,
$\Omega(S,\lambda)$, and
$S$ split locally into products of the orbits $G_\alpha\cdot \lambda$\,,
$G_{\Bbb R,\alpha}\cdot\lambda$\,, and $G_{\Bbb R,\alpha}\cdot {x_0}$\,,
respectively, corresponding to the various positive roots $\alpha\in
\Phi^+(\fg,\ft)$. These splittings are compatible with the definition (7.24) of
$I$ and with symplectic form (7.34) -- note: $\lambda$ induces an
antidominant, regular weight for each $\fg_\alpha^*$ by restriction. The
integer $s$ equals the number of positive compact roots. These
considerations reduce the statement (7.27b) to the special cases of $\GR=
SU(2)$ and $\GR=SU(1,1)$. In these two cases is can be verified by direct
calculation. This completes the proof of lemma 7.10.

\enddemo

We conclude this section with the verification of the fixed point formulas
5.12, 5.24, 5.27 in the case of the discrete series -- i.e., under the hypotheses
enunciated at the beginning of this section. Kashiwara already verified his
conjecture for the discrete series (with anti-dominant
$\lambda$, as in our current setting). Our argument differs from his only by
replacing a reference to a lemma of \cite{OM} by  a short calculation.

Discrete  series characters -- both on the group and the Lie algebra -- are
uniquely determined, among all invariant eigendistributions, by a) their
restriction to the elliptic set, and b) being tempered  \cite{HC5}. The fixed
point formalism discussed in section 4 does produce invariant
eigendistributions. Thus we only need to verify those properties of the
coefficients $c_{g,x}$ and $d_{E,x}$\,, corresponding to $\Cal F = Rj_*\Bbb
C_S$\,, which embody the properties a) and b).

Because of 7.1, a globally well defined section $e^{\lambda-\rho}$ of $\Bbb
C_\lambda$ exists on the universal Cartan $H$, and thus on all of the spaces
mapping naturally to $H$. Harish-Chandra's formula for the discrete series
characters on the elliptic set can be stated as follows:
$$
\Theta(S,\lambda)(g) \, = \, (-1)^s \sum_{x\in X^g\cap
S}\frac{e^{\lambda-\rho}(g,x)}{\prod_{\alpha\in\Phi^+}(1-e^{-\alpha})(g,x)}\ \
\ 
\text{($g$ regular elliptic)};
\tag7.36
$$ the appearance of $(-1)^s$, as opposed to $(-1)^q$ in Harish-Chandra's
original formula, is caused by anti-dominance of $\lambda$\,: in the definition
of the Weyl denominator, Harish-Chandra uses the positive root system which
makes the parameter $\lambda$ dominant. The alternating nature of the
summation in Harish-Chandra's formula translates into the dependence of the
denominator on the fixed point $x$. We now suppose that $g$ lies in the
compact Cartan subgroup $T_\Bbb R$ chosen at the beginning of the section,
and write $t$ instead of $g$ for a generic element of $T'_\Bbb R$, as in
(5.22-25). There are no real roots on $\ft_\Bbb R$\,, so we can choose $\Psi' =
\Phi^+$ in the definition of $\fn'(t,x)$. Hence $N'(t,x)=N^+(t,x)$ is an open
Schubert cell containing the fixed point $x$. In particular, 
$$
\aligned &\chi(\Cal H^*_{N'(t,x)}(\Bbb D Rj_*\Bbb C_S)_x) \ = \ \chi(\Cal
H^*_{N'(t,x)}(j_!\Bbb D_S)_x) \ =
\\ &\qquad\chi(\Cal H^*(j_!\Bbb C_S[2n])_x) \ = \ \cases 1 \ \ \ &\text{if $x\in
S$} \\ 0 \ \ \ &\text{if $x\notin S$\,.} \endcases
\endaligned
\tag7.37
$$ This, in conjunction with 5.25 and 7.6, establishes a).

The temperedness condition b) comes down to the vanishing of certain
coefficients $c_{t,x}$ corresponding to regular semisimple
$t$. We now deviate from the earlier notation of this section: to be consistent
with section 5, $T_\Bbb R$ will be an arbitrary Cartan subgroup,
$t$ a regular element of $T_\Bbb R$, and $x$ a fixed point of $t$. Non-zero
terms in the local expression\footnote{with $\zeta=0$, and re-expressed in
terms of the constants $d_{t,x}$ of (5.25)} (4.3) for
$\Theta = (-1)^s\Theta(Rj_*\Bbb C_S)$ violate temperedness precisely when
they are indexed by a pair $(t,x)\in \tilde G_\Bbb R$ with
$|e^{\lambda_x}(t)| > 1$\,; note: $e^\rho$ is defined up to sign on $\tilde G$, and
hence so is
$e^\lambda=e^{\lambda-\rho}\,e^\rho$. Hence, in view of (5.25), the condition
b) is equivalent to 
$$
\chi(\Cal H^*_{\{x\}}((\Bbb D Rj_*\Bbb C_S)|_{N''(t,x)})) \ = \ 0 \ \ \ \
\text{if
$\ \ |e^{\lambda_x}(t)| > 1$}\,.
\tag7.38
$$ The definition (5.17) of $N''(t,x)$ depends on the choice of the subset $\Psi''$
in $\Phi^+$. The conditions (5.15) allow us to pick 
$$
\Psi'' \ = \ \{\,\alpha\in\Phi^+ \ | \ \text{$\alpha$ is real valued on
$\ft_\Bbb R$ and $e^\alpha(t) >1$}\,\}\,.
\tag7.39
$$ To verify the temperedness condition b), we only need to show that
$$ N''(t,x)\cap S =\emptyset \qquad \text{if\ \ \  $|e^{\lambda_x}(t)|>1$}\,,
\tag7.40
$$ since then $j_!\Bbb C_S|_{N''(t,x)}=0$, hence
$$
\Cal H^*_{\{x\}}((\Bbb D Rj_*\Bbb C_S)|_{N''(t,x)}) \ = \ \Cal H^*_{\{x\}}((j_!\Bbb
C_S[2n])|_{N''(t,x)}) \ = \ 0\,.
$$ The condition (7.40) holds vacuously when the Cartan subgroup $T_\Bbb R$
containing $t$ is compact. Let us assume, then, that $T_\Bbb R$ is a
non-compact Cartan subgroup of a group $\GR$ which does contain some
compact Cartan. In this situation the fixed points of $t\in T'_\Bbb R$ cannot lie
in open $\GR$-orbits. Thus the non-vacuous case of (7.40) is a consequence of:

\proclaim{7.41 Lemma} Let $x\in X$ be a point not lying in any open
$\GR$-orbit and $t\in\GR$ a regular semisimple element fixing $x$, such that
 $|e^{\lambda_x}(t)|\geq 1$. Then $N''(t,x)\cap S =\emptyset$, provided
$\Psi''$ is chosen as in (7.39).
\endproclaim

\demo{Proof} We replace $\lambda$ by a suitable positive integral multiple
 so that the function $f_{\lambda}$ of lemma 7.10 extends to all of
$X$. Let $\zeta\in\fn''(t,x)$ be such that $(\exp\zeta)x\in S$. Because of the
choice of $\Psi$, this implies $\operatorname{Ad}(t^{-n})\zeta \to 0$ as
$n \to +\infty$. We now use the notation of the proof of (7.10): 
$$
\aligned (\exp\zeta)x\in S \ \ &\implies \ \ f_\lambda((\exp\zeta)x)\neq 0
\\ &\implies \ \ h_{nc}\neq 0  \ \ \text{on} \ \ 
\,V^{\fn_{\exp\zeta x}}_\lambda\,=\,\exp\zeta \,V^{\fn_x}_\lambda
\\ &\implies \ \ h_{nc}(\exp \zeta \,v_0,\exp \zeta \,v_0)\neq 0 \ \ \text{for}
\ \ v_0\neq 0 \ \ \text{in} \ \ V_\lambda^{\fn_x}
\\ &\implies  \ \ h_{nc}(t^{-n}\exp \zeta\, v_0,t^{-n}\exp \zeta\, v_0)\neq 0 \,;
\endaligned
$$ 
here we are using the $\GR$-invariance of $h_{nc}$. But 
$$
\aligned &h_{nc}(t^{-n}\exp \zeta\, v_0,t^{-n}\exp \zeta\, v_0) \, =
\\ &h_{nc}(\exp (\operatorname{Ad}(t^{-n})\zeta)\, t^{-n}\,v_0,\exp
(\operatorname{Ad}(t^{-n})\zeta)\, t^{-n}\,v_0)\,=
\\ &|e^{\lambda_x}(t^{-n})|^2\,h_{nc}(\exp (\operatorname{Ad}(t^{-n})\zeta)\,
v_0,\exp (\operatorname{Ad}(t^{-n})\zeta)\, v_0)\,.
\endaligned
$$ Since $\operatorname{Ad}(t^{-n})\zeta \to 0$ and $h_{nc}(v_0,v_0)=0$ --
recall: $x$ does not lie in an open orbit, hence $h_{nc} = 0$ on
$V_\lambda^{\fn_x}$ -- this is possible only if $|e^{\lambda_x}(t^{-n})|\to
\infty$ as $n\to+\infty$. In other words, the hypothesis $(\exp\zeta)x\in S$
forces $|e^{\lambda_x}(t)|<1$.
\enddemo  

The two statements (7.37) and (7.40) now imply the fixed point formula (5.24),
and hence also (5.12), for the discrete series characters $\Theta(S,\lambda)$.
In the case of the characters $\theta(S,\lambda)$ on the Lie algebra, the
argument is virtually identical. Thus:

\proclaim{7.42 Proposition} Under the hypotheses stated at the beginning of
this section, the character formulas (5.12), (5.24), and (5.27) hold.
\endproclaim
\vskip 1cm

\subheading{\bf 8. Coherent continuation}
\vskip .5cm

Both of our character formulas are compatible with the
representation-theoretic process of coherent continuation \cite{S3}. In
the case of the integral formula 3.8 this is not so easy to see a priori,
though it  follows easily from the formula once it has been proved. In this
section, we shall use the process of coherent continuation to establish our
character formulas for the coherent continuations of all discrete series
characters. In effect, we shall remove the positivity condition (7.3) as
hypothesis from propositions 7.18 and 7.45.

Let us recall the notion of coherent continuation. A family of invariant
eigendistributions $\Theta(\lambda)$ parametrized by a coset
$\lambda_0+\Lambda$ of the weight lattice $\Lambda$ depends coherently on
$\lambda$ provided 
$$
\aligned
\text{a)}\ \ &\Cal Z(\fg) \ \, \text{operates on} \ \, \Theta(\lambda)\ \,
\text{according to the infinitesimal}\\
&\text{character} \ \, \chi_\lambda \ \, \text{for every
$\lambda\in\lambda_0+\Lambda$}\,;
\\
\text{b)} \ \ &\text{the coefficients $c_{g,x} = c_{g,x}(\lambda)$ in the local
expressions (4.3)}
\\
&\text{for the} \ \,\Theta(\lambda)\ \,\text{satisfy}\ \,
c_{g,x}(\lambda)=e^{\lambda-\lambda_0}(g)\,c_{g,x}(\lambda_0)\,.
\endaligned
\tag8.1
$$
Note that $\lambda-\lambda_0$ is an integral weight, so
$e^{\lambda-\lambda_0}(g)$ is well defined. The  condition b) can be stated
differently, as follows. Let
$\phi$ be a finite dimensional character. Then 
$$
\phi \ = \ \sum_{\mu\in\Lambda} \ n_\mu(\phi) \, e^\mu \,, 
\tag8.2
$$
with only finitely many $n_\mu(\phi)\neq 0$\,. On any concrete Cartan, this
formula has on obvious meaning; since it is symmetric under the Weyl group,
it makes sense also to state it in universal terms, as above. If the family
$\Theta(\lambda)$ satisfies (8.1a), the second condition (8.1b) is equivalent
to 
$$
\phi\,\theta(\lambda) \ = \ \sum_{\mu\in\Lambda}\ n_\mu(\phi)\,
\Theta(\lambda+\mu)\,,
\tag8.3
$$
for every $\lambda\in\lambda_0+\Lambda$ and every finite dimensional
character $\phi$ \cite{S3}.

Let $\Theta(\lambda)$ be a coherent family parametrized by
$\lambda_0+\Lambda$. If $\lambda_0$ is regular, the single member
$\Theta(\lambda_0)$ of the family determines all the others -- this follows
immediately from (8.1b), coupled with the uniqueness of the coefficients
$c_{g,x}(\lambda)$ in the case of regular parameter $\lambda$. Thus:

\proclaim{8.4 Observation} Two coherent families parametrized by the same
coset $\lambda_0+\Lambda$ coincide as soon as they agree at a single regular
parameter $\lambda\in \lambda_0+\Lambda$.
\endproclaim 

In the definition (2.17), we have attached the invariant
eigendistribution $\Theta(\cf)$ to
$\cf\in\operatorname
D_\GR(X)_{-\lambda}$ and the datum of a specific $\lambda$\,, though the
dependence on $\lambda$ does not come out in the notation -- recall that 
$\operatorname D_\GR(X)_{-\lambda}$ depends only on the image of $\lambda$ 
in $\fh^*/\Lambda$. In the present section, we need to make the dependence on
$\lambda$ explicit. Thus, for $\lambda\in \lambda_0+\Lambda$ and
$\cf\in\operatorname D_\GR(X)_{-\lambda_0}$\, we shall write
$\Theta_\lambda(\cf)$ for the invariant eigendistribution (2.17)
corresponding to the representative $\lambda$ of the coset
$\lambda_0+\Lambda$. 

\proclaim{8.5 Theorem} The invariant eigendistributions
$\Theta_\lambda(\cf)$, $\lambda\in\lambda_0+\Lambda$\,, constitute a 
coherent family.
\endproclaim

This is a standard fact. As specific reference let us mention \cite{SW},
where the coherence is proved for standard sheaves; that is enough, of
course, in view of lemma 6.4. Alternatively, one can use the corresponding
fact about the Beilinson-Bernstein construction -- see \cite{Mi}, for example
-- and carry it over to our setting via the main result in \cite{KSd}. 

Kashiwara's fixed point formalism, as discussed in section 5, associates
another eigendistribution to the datum of a $\lambda\in\fh^*$ and
$\cf\in\operatorname D_\GR(X)_{-\lambda}$, thus a family
$\tilde\Theta_\lambda(\cf)$ parametrized by the $\lambda$-coset of the weight
lattice. In effect, $\tilde\Theta_\lambda(\cf)$ is the invariant
eigendistribution described by the right hand side of the equation in theorem
5.12. 

\proclaim{8.6 Lemma} The family $\tilde\Theta_\lambda(\cf)$ is coherent
\endproclaim

\demo{Proof} The character cycle $c(\cf)$ is a cycle of top degree with
values in the local system $\Bbb C_{-\lambda}$, which is completely
determined by the section $c:(g,x)\mapsto c_{g,x}$ of the dual local system
$\Bbb C_\lambda$ over $\tilde G'_\Bbb R$ -- cf. (4.12). The coherence
condition (8.1b) for the family $\tilde\Theta_\lambda(\cf)$ is equivalent to
the following statement: when $\lambda$ is replaced by the translate
$\lambda+\mu$ by some $\mu\in\Lambda$, the section $c:(g,x)\mapsto
c_{g,x}(\lambda)$ gets multiplied by $e^\mu$. This multiplicative behavior of
$c$ is clear from the construction; the dependence of the character cycle
on the specific parameter $\lambda$, rather than on the coset
$\lambda+\Lambda$, appears at exactly one point, the passage from $\phi'$ to
$\phi''$ in (5.7-8) which is multiplicative in the sense mentioned earlier.
\enddemo

The general principle (8.4), theorem 8.5, and lemma 8.6 allow us to remove
the positivity hypothesis (7.3) in proposition 7.45. For later reference, we
state 

\proclaim{8.7 Proposition} The character formulas (5.12), (5.24), and (5.27)
are satisfied by every standard sheaf associated to an open $\GR$-orbit,
provided the group $\GR$ contains a compact Cartan subgroup.
\endproclaim

According to (8.4) the invariant eigendistributions (7.2) constitute a coherent
family $\Theta_\lambda(Rj_*\Bbb C_S)$, parametrized by
$\lambda\in\Lambda-\rho$\,. Let us write $\theta_\lambda(Rj_*\Bbb C_S)$ for
the corresponding family on the Lie algebra. The conditions (8.1) and (8.2)
have obvious analogues on the Lie algebra. Thus it makes sense to talk about
coherent families of invariant eigendistributions on the Lie algebra. The
coherence of $\Theta_\lambda(Rj_*\Bbb C_S)$ implies coherence also for the
family $\theta_\lambda(Rj_*\Bbb C_S)$. Unlike the $\Theta_\lambda(Rj_*\Bbb
C_S)$, which are defined only for $\lambda\in\Lambda-\rho$, the
$\theta_\lambda(Rj_*\Bbb C_S)$ can be given meaning for any
$\lambda\in\Bbb R \otimes_\Bbb Z\Lambda$, as follows. The integral on the left 
in (7.14) converges for every regular $\lambda\in\Bbb R\otimes_\Bbb
Z\Lambda$, and thus defines a family of invariant eigendistributions
$$
\phi \ \ \longmapsto \ \ \int_{\Omega(S,\lambda)}\hat\phi\,\sigma_\lambda^n
\qquad\qquad (\,\phi\in C_c^\infty(\gr)\,),
\tag8.8
$$
parametrized by the regular set in $\Bbb R\otimes_\Bbb Z\Lambda$\,.
According to Rossmann's theorem 7.8, this distribution coincides with
the  discrete series character $\theta_\lambda(Rj_*\Bbb C_S)$ when
$\lambda\in\Lambda-\rho$ satisfies the antidominance condition (7.3):
$$
\gathered
\int_ \gr \phi \,\theta_\lambda(Rj_*\Bbb
C_S)\,  dx \ = \ \frac 1
{(2\pi i)^n n!}\int_{\Omega(S,\lambda)}\hat\phi\,\sigma_\lambda^n
\\
\vspace{1.5\jot}
\qquad\qquad\text{for} \ \ \lambda\in\Lambda-\rho \ \ \text{regular
antidominant}\,.
\endgathered
\tag8.9
$$
Since the right hand side is well defined for any regular antidominant
$\lambda\in\Bbb R\otimes_\Bbb Z\Lambda$, so is the family
$\theta_\lambda(Rj_*\Bbb C_S)$. At this point, then, the family is well defined
for all $\lambda\in\Bbb R\otimes_\Bbb Z\Lambda$ which either lie in
$\Lambda-\rho$, or are antidominant regular.
The following result is due to Harish-Chandra \cite{HC5, lemma 32}.

\proclaim{8.10 Lemma} The coefficients $d_{E,x}$ in the local
expressions (4.4) for the invariant eigendistributions (8.8) depend only on
the Weyl chamber in which $\lambda$ lies.
\endproclaim

Because of the equality (8.9), we conclude that the coefficients in the local
expressions for the $\theta_\lambda(Rj_*\Bbb
C_S)$ are also independent of $\lambda$. Thus we can coherently continue the  
$\theta_\lambda(Rj_*\Bbb C_S)$ from any antidominant regular
$\lambda_0\in\Bbb R\otimes_\Bbb Z\Lambda$ to $\lambda_0+\Lambda$, and
therefore to all of $\Bbb R\otimes_\Bbb Z\Lambda$. In effect, this family is
coherent in the strongest possible sense: defined for all $\lambda\in\Bbb
R\otimes_\Bbb Z\Lambda$, with coefficients $d_{E,x}$ independent of $\lambda$.
The constancy of the $d_{E,x}$ implies that the values of the function
$\theta_\lambda(Rj_*\Bbb C_S)$ on the regular set in $\gr$ depend real
analytically on $\lambda$.  This makes the family of
distributions $\theta_\lambda(Rj_*\Bbb C_S)$ weakly analytic: the integral of
the family against any test function is real analytic in the parameter
$\lambda$.

There is a second family attached
to the datum of the sheaf $\cf=Rj_*\Bbb C_S$\,, namely the right hand side of
the equation in theorem 3.8. Let us denote this family by
$\tilde\theta_\lambda(Rj_*\Bbb C_S)$. Proposition 7.15 asserts:
$$
\theta_\lambda(Rj_*\Bbb C_S) \ = \ \tilde\theta_\lambda(Rj_*\Bbb C_S) \ \
\text{when  $\lambda\in\Lambda-\rho$ is regular antidominant}\,.
\tag8.11  
$$
The definition of 
$\tilde\theta_\lambda(Rj_*\Bbb C_S)$ involves integration over the
characteristic cycle of the sheaf $Rj_*\Bbb C_S$ and makes sense for every
$\lambda\in\fh^*$, not just for $\lambda\in\Lambda-\rho$. The coadjoint orbit
$\Omega(S,\lambda)$ as defined in (7.7) also has meaning for every
$\lambda\in\Bbb R\otimes_\Bbb Z\Lambda$. Our proof of the identity (7.14)
uses the fact that a positive integral multiple of $\lambda$ satisfies the
positivity condition (7.3) -- the integrality of $\lambda-\rho$ plays no role.
Hence (7.14) remains valid for every antidominant regular $\lambda\in
\Bbb Q\otimes_\Bbb Z\Lambda$. Combining this with (8.9), we conclude
$$
\theta_\lambda(Rj_*\Bbb C_S) \ = \ \tilde\theta_\lambda(Rj_*\Bbb C_S) \ \
\text{for  $\lambda\in
\Bbb Q\otimes_\Bbb Z\Lambda$  regular antidominant}\,.
\tag8.12
$$
Proposition 3.7 asserts that the family
$\tilde\theta_\lambda(Rj_*\Bbb C_S)$, $\lambda\in\fh^*$, is weakly
holomorphic in $\lambda$, in the sense that its value  on any particular test
function depends holomorphically on $\lambda$. In particular, it is weakly
(real) analytic when restricted to  $\Bbb R\otimes_\Bbb Z\Lambda$. Two
weakly analytic families which coincide on a large enough set must coincide,
hence, by (8.12):

\proclaim{8.13 Proposition} The character formula (3.8)
is satisfied by every standard sheaf associated to an open $\GR$-orbit,
provided the group $\GR$ contains a compact Cartan subgroup.
\endproclaim
\vskip 1cm

\subheading{\bf 9. Induction}
\vskip .5cm

In the previous two sections we proved our character formulas for standard
sheaves associated to open orbits of groups $\GR$ which contain a compact
Cartan subgroup. We shall now extend the validity to standard sheaves
attached to a larger class of orbits -- roughly speaking orbits which fiber
over a closed orbit in a generalized flag variety, with fibers of the type we
have considered so far. In representation theoretic terms, we extend the
validity of the character formulas by the process of parabolic induction.
This has a counterpart on the level of standard sheaves, as we shall explain
next. 

Let $S = S(T_\Bbb R,\tau_{x_0}) = \GR \cdot {x_0}$ be a $\GR$-orbit attached
to a concrete Cartan subgroup  $T_\Bbb R\subset \GR$ and fixed point
${x_0}$ of
$T_\Bbb R$, as in (6.5). The isomorphism  $\tau_{x_0}:\ft\to\fh$ determined
by the fixed point ${x_0}$ pulls back the universal positive root system
$\Phi^+$ to a positive root system  $\Phi^+_S(\fg,\ft)$ in $\Phi(\fg,\ft)$, the
root system for
$(\fg,\ft)$. We call $\alpha\in\Phi(\fg,\ft)$ real, imaginary, or complex
depending on whether $\alpha = \bar\alpha$ (= complex conjugate of
$\alpha$), $\alpha = -\bar\alpha$, or $\alpha \neq \pm\bar\alpha$. The
following two properties of
$\Phi^+_S(\fg,\ft)$ are equivalent:
$$
\aligned &\text{a)\ \, for every complex root $\alpha\in\Phi^+_S(\fg,\ft)$, the
root
$\bar\alpha$ is also positive\,;}
\\ &\text{b)\ \, for every complex simple root $\alpha\in\Phi^+_S(\fg,\ft)$,
the root
$\bar\alpha$ is positive\,.}
\endaligned
\tag9.1
$$ Though phrased in terms of $T_\Bbb R$ and ${x_0}$, both a) and b) are
really properties of the orbit $S$. The preceding statements are easy to
verify -- see \cite{S4, lemma 6.14} for example\footnote{\cite{S4} treats
the case of $K$-orbits, but also supplies the translation via Matsuki duality
between
$K$-orbits and $\GR$-orbits.}. We shall call the orbit $S$ \lq\lq maximally
real" if it satisfies the equivalent conditions (9.1a,b). It is these orbits that
will be considered in the current section. 

Let us suppose, then, that the orbit $S = S(T_\Bbb R,\tau_{x_0}) = \GR \cdot
{x_0}$   is maximally real. As explained in \S 6, standard sheaves attached to
the orbit
$S$ are determined by the following sets of data: a linear function
$\lambda\in \fh^*$ and a character $\chi:\TR \to \Bbb C^*$, such that
$d\chi=\tau_{x_0}^*(\lambda-\rho)$; the pair $(\lambda,\chi)$ induces an
irreducible, $\GR$-equivariant, $(-\lambda-\rho)$-monodromic local system
$\Cal L_\chi$ on $S$. The direct image $\cf= Rj_*\Cal L_\chi$ via the
embedding $j:S\to X$ is the standard sheaf we shall work with. 

We need to assemble various pieces of known structural information. To
begin with, the Cartan subgroup $\TR$ has a direct product decomposition 
$$
\TR \ = \  C_\Bbb R \cdot A_\Bbb R\,, \qquad \text{with $\CR$ compact and
$\AR$ connected and split}\,.
\tag9.2a
$$ The fact that $\TR$ is fixed by an anti-holomorphic involution of the
complex torus $T$ implies 
$$
\gathered
\CR \ = \ C_\Bbb R^0 \cdot F \qquad \text{(direct product), }\\
\text{with $F\subset\GR\cap\exp(i\frak a_\Bbb R)$\,, $F\cong \Bbb Z/2\Bbb
Z
\times\dots\times\Bbb Z/2\Bbb Z$\,.}
\endgathered
\tag9.2b
$$ In particular, 
$$
\aligned
\chi \ &= \ \chi_C\cdot\chi_A\,, \qquad\text{with $\chi_C:\CR\to\Bbb C^*$
and
$\chi_A:\AR\to\Bbb C^*$}\,,
\\
\chi_C\ &= \ \chi_C^0\cdot \chi_F\,, \qquad\text{with $\chi_C^0:C_\Bbb
R^0\to\Bbb C^*$ and
$\chi_F:F\to\Bbb C^*$}\,.
\endaligned
\tag9.3
$$ Since $C_\Bbb R^0$ and $\AR$ are connected, the compatibility conditions 
$$ d\chi_C^0 \ = \ \tau_{x_0}^*(\lambda-\rho)|_{\fc_\Bbb R}\,,\qquad
d\chi_A\ =
\
\tau_{x_0}^*(\lambda-\rho)|_{\fa_\Bbb R}
\tag9.4
$$ completely determine $\chi_C^0$ and $\chi_A$. Following standard
notation, we write the centralizer of $\fa$ as a direct product 
$$
\gathered Z_G(\fa) \ = \ M\cdot A\,, \qquad\text{(direct product),}
\\
\text{with $M$ connected and defined over $\Bbb R$}
\endgathered
\tag9.5a
$$ and correspondingly,
$$
\gathered Z_{\GR}(\fa) \ = \ M_\Bbb R\cdot F\cdot A_\Bbb R\,,
\qquad\text{(direct product),}
\\
\text{where $M_\Bbb R= M\cap\GR$}\,.
\endgathered
\tag9.5b
$$ We should warn the reader that $\AR$ was defined to be the identity
component of $A\cap\GR$\,; in fact, $A\cap\GR = \AR\cdot F$.  Because of our
hypothesis (9.1)
$$
\fv \ = \ \oplus\ \{\, \fg^{-\alpha} \ | \ \alpha\in\Phi^+_S(\fg,\ft), \
\alpha|_{\fa}
\neq 0\ \}
\tag9.6
$$ is a subalgebra of $\fn_{x_0}$, defined over $\Bbb R$, normalized by $MA$
but linearly disjoint from $\fm\oplus\fa$\,; it is the nilradical of the parabolic
subalgebra
$$
\fp \ = \ \fm\oplus\fa\oplus\fv \qquad \text{(semidirect product)}\,.
\tag9.7a
$$  Thus $P$, the normalizer of $\frak p$ in $G$, is a parabolic subgroup,
defined over
$\Bbb R$, with Levi decomposition 
$$ P\ = \  M\cdot A\cdot V \qquad \text{(semidirect product)}\,.
\tag9.7b
$$  Here $V\subset G$ denotes the connected subgroup with Lie algebra
$\fv$\,; its group of real points $V_\Bbb R = V\cap \GR$ is also connected.
Unlike the complex parabolic $P$, 
$$ P_\Bbb R \ =_{\text{def}} \ P\cap \GR \ = \ M_\Bbb R \cdot A_\Bbb R\cdot
F\cdot  V_\Bbb R
\tag9.7c
$$ is not connected in general. In fact, the decomposition (9.7c) is
topologically direct, so $F$ can be identified with the component group of
$P_\Bbb R$. 

The $G$-conjugates of $\fp$ constitute the underlying set of a generalized
flag variety $Y$, the base of a $G$-equivariant fibration
$$ X \ @>{\ \ \ }>> \ Y \qquad \text{with fiber $X_M\cong$ flag variety of $M$}\,.
\tag9.8
$$ We identify $X_M$ concretely with the fiber through the point ${x_0}$.
Then (9.8) induces a $\GR$-equivariant fibration of orbits
$$ S \ @>{\ \ \ }>> \ S_Y \qquad \text{with fiber $S_M=$ $M_\Bbb R$-orbit
through
${x_0}$ }\,.
\tag9.9
$$ The $\GR$-orbit $S_Y\subset Y$ can be identified with $\GR/P_\Bbb R$,
hence is a compact real form of the complex manifold $Y$. On the other hand,
the
$M_\Bbb R$-orbit $S_M\subset X_M$ is associated to the compact Cartan
subgroup $C^0_\Bbb R\subset M_\Bbb R$\,, in the sense that $C^0_\Bbb R$
has
${x_0}$ as fixed point, hence is open. Since $M_\Bbb R$ stabilizes $X_M$,
$$
\cf_M \ = \ \cf|_{X_M}
\tag9.10
$$ is an object in $\operatorname D_{M_\Bbb R}(X_M)_{-\lambda_M}$, the
twisted $M_\Bbb R$-equivariant derived category with twisting parameter
$\lambda_M=_{\text{def}}$ restriction of $\lambda$ to the universal Cartan
$\fh_M$ for $M$. Two comments are in order. We use the point ${x_0}$ to
identify
$\fh\cong \ft$ and $\fh_M \cong \fc$, and this allows us to regard $\fh_m$ as
lying in $\fh$; the same choice of ${x_0}$ in the description (2.9) of the
enhanced flag variety provides a distinguished embedding $\hat
X_M\hookrightarrow
\hat X$ compatible with the embedding $X_M\hookrightarrow X$ and the
embedding of groups $M_\Bbb R\times H_M\hookrightarrow \GR\times H$.
The category  $\operatorname D_{G_\Bbb R}(X)_{-\lambda}$ is built from
$(-\lambda-\rho)$-monodromic sheaves, whereas the definition of
$\operatorname D_{M_\Bbb R}(X_M)_{-\lambda_M}$ involves
$(-\lambda_M-\rho_M)$-monodromic sheaves. Thus, to see that restriction
from $X$ to $X_M$ maps $\operatorname D_{G_\Bbb R}(X)_{-\lambda}$ to
$\operatorname D_{M_\Bbb R}(X_M)_{-\lambda_M}$, we need to know 
$$
\rho_M \ = \text{restriction of $\rho$ to $\fh_M$ via $\fh_M\subset\fh$}\,;
\tag9.11
$$ that follows from our hypothesis (9.1). Tracing through the definitions, one
finds:

\proclaim{9.12 Lemma} The restricted sheaf $\cf_M$ is the standard sheaf
associated to the open orbit $S_M\subset X_M$ and the data
$(\lambda_M,\chi^0_C)$.
\endproclaim

The character $\Theta(\cf)$ and its $M$-analogue $\Theta_M(\cf_M)$ are
related by parabolic induction. Recall the notion  of normalized\footnote{i.e.,
with built-in $\rho$-shift} parabolic induction,
$$ (\tau,\chi_F,\nu) \ \mapsto \ I_{P_\Bbb R}^\GR(\tau\otimes\chi_F\otimes
e^\nu) \,,
\tag9.13
$$ which associates an admissible representation of $\GR$, of finite length,
to any triple $ (\tau,\chi_F,\nu)$ consisting of an admissible $M_\Bbb
R$-representation $\tau$ of finite length, a character $\chi_F: F\to \Bbb C^*$,
and a linear function $\nu\in\fa^*$.   The induction functor
$I_{P_\Bbb R}^\GR$ can also be applied to virtual representations, and hence
to virtual characters. 

\proclaim{9.14 Proposition} $\Theta(\cf)= I_{P_\Bbb
R}^\GR(\Theta_M\cdot \chi_F\cdot e^\nu)$\,, with $\Theta_M=
\Theta_M(\cf_M)$\,, $\chi_F$\,,  $\nu = \lambda|_{\fa}$ as in (9.3), and
$m=\dim _\Bbb C Y$. 
\endproclaim

In view of (9.12), when $\lambda_M$ is regular anti-dominant  with respect to
$M$, this is a statement about standard representations, which can be found
in
\cite{SW}; alternatively -- still in the anti-dominant situation -- this follows
from  the duality theorem of \cite{HMSW1}.  Parabolic induction is
compatible with coherent continuation -- see, for example, \cite{HS} -- so
the anti-dominance assumption can be dropped. Note that the standard
sheaf $\cf$ is the direct image of a locally constant sheaf on the orbit $S$, in
degree zero; dually, $\Bbb D \cf$ is concentrated in degree $-(m+2\dim
_\Bbb C X_L)$. On the other hand,   the sheaf of hyperfunctions on $S_Y$ is the
local cohomology sheaf along $Y$ of $\Cal O_Y$, in degree $m$, so the sign
changes one might expect cancel. 

We should remark that the restriction operation  $\operatorname D_{G_\Bbb
R}(X)_{-\lambda}\to \operatorname D_{M_\Bbb R}(X_M)_{-\lambda_M}$ as in
(9.10)  factors,
$$
\operatorname D_{G_\Bbb R}(X)_{-\lambda}\to \operatorname
D_{(MA)_\Bbb R}(X_{MA})_{-\lambda}\to \operatorname D_{M_\Bbb
R}(X_M)_{-\lambda_M}\,.
\tag9.15
$$ The first of these arrows has a right adjoint, \lq\lq sweeping out" along
the orbit $S_Y\subset Y$, which is closed. On the level of representations,
this corresponds to parabolic induction, whether or not the sheaf in question
is a standard sheaf. We have chosen to go down all the way to
$\operatorname D_{M_\Bbb R}(X_M)_{-\lambda_M}$ for technical reasons
only: the group $M_\Bbb R$ has a compact Cartan subgroup, so that the
discussion of \S\S 7,8 applies directly. 

There exists an explicit formula for the induced character $I_{P_\Bbb
R}^\GR(\Theta_M\cdot \chi_F\cdot e^\nu)$, in terms of the inducing data
\cite{Hi}, which we shall recall later. As one consequence of this formula, the
restriction of the induced character  to a small neighborhood of the identity
depends only on $\Theta_M$ and $\nu$, not on $\chi_F$.   Every connected
component of  $\fg'_\Bbb R$, the set of regular semisimple elements in
$\gr$\,, meets any neighborhood of $0$.  It follows that the pullback of the
induced character to the Lie algebra does not depend on $\chi_F$.
Accordingly, we suppress the symbol $\chi_F$ in the formula 
$$
\theta(\cf)\ = \ I_{P_\Bbb R}^\GR(\theta_M\cdot e^\nu) \qquad
(\,\theta_M=\theta_M(\cf_M)\,, \ \nu = \lambda|_\fa\,)\,,
\tag9.16
$$ which is the Lie algebra analogue of (9.14). 

\proclaim{9.17 Proposition} The character formula (3.8) is satisfied by every
standard sheaf associated to a maximally real orbit.
\endproclaim

\demo{Proof} In effect, we must show that the character formula (3.8) is
satisfied by the sheaf $\cf = Rj_*\Cal L_\chi$ constructed at the beginning of
this section, since every standard sheaf associated to a maximally real orbit
is of this form. Just as in the group case, parabolic induction of characters on
the level of the Lie algebra is given by an explicit formula\footnote{indeed,
the formula on the Lie algebra can be deduced from that on the group}. This
formula can be applied to any invariant eigendistribution on $\fm_\Bbb R$,
whether or not it comes from a character on the group $M_\Bbb R$.
Moreover, when the induction process is applied to a holomorphic family of
invariant eigendistributions on $\fm_\Bbb R$, the resulting  invariant
eigendistribution on $\gr$ depends holomorphically on the parameter of the
inducing family, and also on the parameter $\nu$. In section 8 we showed that
$\theta_M =
\theta_M(\cf_M)$ -- which is specified by the parameter $\lambda|_\fc$ and
an open $M_\Bbb R$-orbit  in $X_M$ -- can be continued to a holomorphic
family parametrized by $\fc^*$. It follows that the induced character (9.16)
depends holomorphically on $\lambda$. The characteristic cycle $\occ(\cf)$
depends only on the orbit $S$, so the family of invariant eigendistributions
$$
\phi \ \ \mapsto \ \  \frac 1 {(2\pi i)^n n!}
\int_{\occ(\cf)} \mu^*_\lambda \hat \phi \ (-\sigma + \pi^* \tau_\lambda)^n
\qquad (\,\lambda\in \fh^*\,)
\tag9.18
$$ is holomorphic in the parameter $\lambda$ -- cf. proposition 3.7. Hence it
suffices to prove the statement of the proposition under the following
hypotheses:
$$
\text{ a) \ $\lambda$ is regular; \ \ \ \ \ \ \ b) \ $\nu\in i\fa_\Bbb R^*$}\,,
\tag9.19
$$ though we cannot assume that $\lambda|_\fc$ is integral. 

Let us describe the characteristic cycle $\occ(\cf)$ in terms of
$\occ(\cf_M)$. Locally near any point in $Y$, the fibration (9.8) is trivial; thus
locally over $Y$ the sheaf $\cf = Rj_*\Cal L_\chi$ is the exterior product of
two sheaves: the direct images of the constant sheaves $\Bbb C_{S_M}$ and
$\Bbb C_{S_Y}$ under, respectively, the open embedding
$S_M\hookrightarrow X_M$ and the closed embedding $S_Y\hookrightarrow
Y$. Correspondingly,  $\occ(\cf)$ is -- again locally over $Y$ -- the product
of the characteristic cycle of the two sheaves. The orientation of the 
product cycle is independent of the order of the product, since both factors
are even dimensional cycles. The orientation of the characteristic cycle of
the direct image of $\Bbb C_{S_M}$ was pinned down in \S 7. We orient the
conormal bundle $T^*_{S_Y}Y$ by means of the {\it negative\/} of the
imaginary part of the holomorphic symplectic structure on $T^*Y$. We
claim:  $T^*_{S_Y}Y$,  oriented as above,  coincides with the characteristic  
cycle of the direct
image of $\Bbb C_{S_Y}$. In effect, this is a statement about
$\Bbb R^m\subset \Bbb C^m$; our conventions, as set up in \cite{SV4},
reduce this further to the case of $\Bbb R\subset \Bbb C$, where it can be
checked directly.

We may as well suppose that the Cartan subgroup $\TR\subset \GR$ is
invariant under the Cartan involution. The Levi component $MA$ of $P$ is
then also invariant under the Cartan involution. Thus $M\cap U_\Bbb R$ -- the
intersection of $M$ with the compact real form $U_\Bbb R\subset G$ which
was used to define the twisted moment map $\mu_\lambda$ -- is a compact
real form of $M$. We use it to define the twisted moment map
$$
\mu_{M,\lambda-\nu}\ : T^*X_M \ \longrightarrow \Omega_{M,\lambda-\nu}
\tag9.20
$$ for $M$ and the twisting parameter $\lambda|_\fc$\,, which we tacitly
identify with $\lambda-\nu$; here $\Omega_{M,\lambda-\nu}\subset \fm^*$
denotes the $M$-orbit of $\lambda|_\fc=\lambda-\nu$. 

Next, we associate a cycle $C$ in $\Omega_\lambda$ to the cycle
$\mu_{M,\lambda-\nu}\occ(\cf_M)$ in $\Omega_{M,\lambda-\nu}$; the cycle
$C$ will turn out to be the $\mu_\lambda$-image of $\occ(\cf)$. As a set, 
$$ C\ = \  \KR\text{-orbit of } \
\left(\mu_{M,\lambda-\nu}\occ(\cf_M)+\nu+(\gr+\fp)^\perp\right)\,.
\tag9.21
$$ Here $(\gr+\fp)^\perp$ refers to the annihilator in $\fg^*$ of
$(\gr+\fp)^\perp$, relative to the real part of the pairing $\fg\times\fg^* \to
\Bbb C$; under the isomorphism $\fg^*\cong \fg$ induced by the Killing form, 
$(\gr+\fp)^\perp$ corresponds to $i\fv_\Bbb R$. Locally near $\lambda$, $C$
splits into a product $C\cong \mu_{M,\lambda-\nu}\occ(\cf_M) \times C'$,
where $C'$ is a submanifold of $\fm^\perp\cap i\fg^*_\Bbb R$, passing
through
$\nu$ (recall (9.19b)!), with tangent space $(\fm+\fa)^\perp\cap i\fg^*_\Bbb
R$ at $\nu$. Note that $(\fm+\fa)^\perp\cap i\fg^*_\Bbb R$ corresponds to
$i\fv_\Bbb R \oplus i(\fv_{\text{opp}})_\Bbb R$\,. Once this tangent space is
oriented, its  $\KR$-translates become consistently oriented, providing an
orientation of $C'$. Finally,  we orient 
$$ T_\nu C' \ \cong \  (\fm+\fa)^\perp\cap i\fg^*_\Bbb R
\tag9.22
$$ by means of the imaginary part of the holomorphic symplectic form of
$\Omega_\lambda$. This form is non-degenerate on the image of
$(\fm+\fa)^\perp$ in $T_\lambda\Omega_\lambda$, and defined over $\Bbb
R$ relative to the natural real structure of $(\fm+\fa)^\perp$;  therefore its
imaginary part does specify an orientation of the space (9.22). Since
$\mu_{M,\lambda-\nu}\occ(\cf_M)$ is a cycle, the orientation of $C'$ now
determines $C$ as a cycle.

\proclaim{9.23 Lemma} $\mu_\lambda\occ(\cf)\ = \  C$.
\endproclaim

We shall prove the lemma after completing the proof of proposition 9.17.
Because of the hypothesis (9.19a), we can rewrite the integral (9.18) and its
$M$-analogue as integrals over the cycles
$\mu_\lambda(\occ(\cf))\subset\Omega_\lambda$, respectively
$\mu_{M,\lambda-\nu}(\occ(\cf_M))\subset\Omega_{M,\lambda-\nu}$, as in
(3.9). The integral over $\mu_{M,\lambda-\nu}(\occ(\cf_M))$ represents the
virtual character $\theta_M(\cf_M)$ -- that is the main result of the previous
section. We need to show that the integral over the cycle
$\mu_\lambda(\occ(\cf))$ represents the virtual character (9.16). Lemma
9.23 reduces the problem to the following assertion: if a virtual character
$\theta_M$ on the Lie algebra $\fm_\Bbb R$ is represented as an integral of
the type (3.9) over the cycle $C_M$ in $\Omega_{M,\lambda-\nu}$, then the
induced character $I_{P_\Bbb R}^\GR(\theta_M\cdot e^\nu)$ is represented
by the integral over the cycle $C$ defined below (9.21), with the unspecified
cycle $C_M$ playing the role of
$\mu_{M,\lambda-\nu}(\occ(\cf_M))$. This last assertion is established by
Rossmann \cite{R2, \S 2}, who generalizes an argument  of Duflo \cite{D1}.
Rossmann specifies the orientation of the \lq\lq induced cycle" in
$\Omega_\lambda$ only implicitly; a careful examination of his proof shows
that he uses the symplectic form of $\Omega_\lambda$, exactly as we have
done in the description of $C$. At this point the proof of proposition 9.17 is
complete, except for the verification of  (9.23). 

\enddemo

\demo{Proof of Lemma 9.23} The two cycles are $\KR$-invariant, as is the
moment map $\mu_\lambda$, and $\KR$ acts transitively on $S_Y$. It
therefore suffices to identify the two cycles over a single point in $S_Y$.
The embedding
$(\fm\oplus\fa)^*\hookrightarrow \fg^*$ induced by the Killing form which
was used to identify $X_M$ with a specific fiber of $X\to Y$ can also be used
to embed $T^*X_M\hookrightarrow T^*X$ compatibly with
$X_M\hookrightarrow X$ and $MA$-equivariantly. The twisted moment maps
$\mu_\lambda$ and
$\mu_{M,\lambda-\nu}$ fit into a commutative diagram
$$
\CD T^*X_M @>{\ \ \ \mu_{M,\lambda-\nu}\ \ \ }>> \Omega_{M,\lambda-\nu} \\
@VVV  @VVV
\\ T^*X @>{\pretend \mu_\lambda\haswidth{\ \ \ \mu_{M,\lambda-\nu}\ \ \
}}>>
\Omega_\lambda
\endCD
\tag9.24
$$ in which the right vertical arrow is induced by $\fm^* \hookrightarrow
\fg^*$, followed by translation by $\nu$. Over the base point $y_0$ in $S_Y$,
our passage from $\occ(\cf_M)$ to $\occ(\cf)$ as explained above amounts to
adding the inverse image in $T^*X$ of the conormal space $(T^*_{S_Y} Y)_y$ 
to
$\occ(\cf_M)$ and taking the $\KR$-saturation. Via the twisted moment map
$\mu_\lambda$, the inverse image of $(T^*_{S_Y} Y)_{y_0}$ maps
isomorphically to the $\nu$-translate of $(\gr+\fp)^\perp$, so the passage
from
$\mu_{M,\lambda-\nu}(\occ(\cf_M))$ to $C$ is the isomorphic
$\mu_\lambda$-image of the passage from $\occ(\cf_M)$ to $\occ(\cf)$,
except possibly for the orientation. We used the symplectic form  of
$\Omega_\lambda$ to orient  $(\gr+\fm+\fa)^\perp$  and the symplectic
form of $T^*Y$ to orient $T^*_{S_Y}Y$. We also regarded $X$, locally over
${y_0}$, as isomorphic to the  product $X_M\times Y$; this allows us to
regard $T^*Y$ as a factor of $\ct$ and to restrict the symplectic form of
$\ct$ to $T^*Y$, where it must agree with its symplectic form. In other
words, the symplectic form of
$\ct$ can be used to orient $T^*_{S_Y}Y$. Proposition 3.3 relates the
symplectic form of $\Omega_\lambda$ to that of $\ct$\,: they correspond to
each other via $\mu_\lambda$, except for a sign change and the addition of
the two form $\pi^*\tau_\lambda$. Since that form is pulled back from the
base, it can be checked that  $(-\sigma + \pi^*\tau_\lambda)^m$ agrees with
$(-\sigma)^m$ on $T^*_{S_Y}Y$ -- what matters here is that the (real)
$2m$-manifold $T^*_{S_Y}Y$ is the conormal bundle of a (real)
$m$-manifold. Since we had oriented $T^*_{S_Y}Y$ by the imaginary part  of
$-\sigma$, we conclude that the two cycles
$\mu_\lambda\occ(\cf)$ and $C$ agree, as asserted by the lemma. 

\enddemo

In effect, we have shown that the integral formula 3.8 is compatible with
parabolic induction. We need to do the same for the fixed point formula:

\proclaim{9.25 Proposition} The character formulas (5.12), (5.24), and (5.27)
are satisfied by every standard sheaf associated to a maximally real orbit.
\endproclaim

\demo{Proof} Because of proposition 6.2, we only need to show that one of
the two formulas in theorem 5.24 for the $c_{t,x}$ is compatible with
parabolic induction -- that not only implies the compatibility with induction
of the formula (5.12), but also of (5.27), which is really a special case of
(5.24). 

We begin by recalling the explicit formula for the character of an induced
representation in terms of the inducing character \cite{Hi}. Let $\Theta_M$,
$\chi_F$, and
$\nu\in\fa^*$  be inducing data as in (9.14). Then $I_{P_\Bbb
R}^\GR(\Theta_M\cdot \chi_F\cdot e^\nu)$ is supported on the union of the
$\GR$-conjugates of $Z_{\GR}(\fa) \ = \ M_\Bbb R\cdot F\cdot A_\Bbb R$, so
it suffices to describe this induced character on any Cartan subgroup
 of $\GR$ which is contained in $M_\Bbb R\cdot F\cdot A_\Bbb R$. Since we
have used the symbol $\TR$ for the specific Cartan subgroup used in
describing the orbit $S$, let us write $\ttr$ for  a typical Cartan subgroup 
contained in $M_\Bbb R\cdot F\cdot A_\Bbb R$.  Then, for any regular $t\in
\ttr$,
$$
\aligned &\left(I_{P_\Bbb R}^\GR(\Theta_M\cdot \chi_F\cdot e^\nu)\ 
\tsize\prod_{\alpha\in\Phi^+} |e^{\alpha/2}-e^{-\alpha/2}| \
\right)(t) \ = \ 
\\ 
&\qquad\sum_{g\in R/\ttr} \ \left((\Theta_M\cdot \chi_F\cdot e^\nu)\ 
\tsize\prod_{\alpha\in\Phi_M^+} |e^{\alpha/2}-e^{-\alpha/2}|
\ \right)(gtg^{-1})\,,
\\ &\text{where} \ \ R=\{g\in\GR \ | \ g\ttr g^{-1}\subset MA\}\,.
\endaligned
\tag9.26
$$ 
In this formula, the particular choices of positive root system $\Phi^+$
for the pair $(\fg,\tilde\ft)$ and $\Phi_M^+$ for $(\fm\oplus\fa,\tilde\ft)$ do
not matter since we are taking absolute values. On the right in (9.26), and in
various formulas below, we evaluate the functions $\Theta_M$, $\chi_F$, and $
e^\nu$ at points $gtg^{-1}$, with $g\in R$ and $t\in\ttr$, in which case 
$gtg^{-1}$ lies in the direct product $M_\Bbb R\cdot F\cdot A_\Bbb R$; thus
$\Theta_M(gtg^{-1})$ denotes the value of $\Theta_M$ at the $M_\Bbb
R$-component of $gtg^{-1}$, and with the analogous convention in the other
two cases. 

Let us re-write both sides of (9.26) in terms of the constants $c_{t,x}$ as in
(4.3). With $t\in \ttr$ regular and $\zeta\in\tilde\ft_\Bbb R$ sufficiently
small, 
$$ I_{P_\Bbb R}^\GR\left(\Theta_M\cdot \chi_F\cdot e^\nu\right)(t
\exp(\zeta)) \ = \ \sum_{x\in X^{\ttr}}\ \frac { c_{t,x}\, e^{\lambda_x(\zeta) -
\rho_x(\zeta)}} {\tsize \prod_{\alpha \in {\Phi}^{ +}} (1-e^{-\alpha_x})(t
\exp(\zeta))}\,.
\tag9.27
$$ On the other hand, 
$$
\Theta_M (t \exp(\zeta)) \ = \ \sum_{x\in X_M^{\ttr}}\
\frac { c^M_{t,x}\, e^{(\lambda-\nu)_x(\zeta) - \rho_{M,x}(\zeta)}} {\tsize
\prod_{\alpha
\in {\Phi}_M^{ +}} (1-e^{-\alpha_x})(t \exp(\zeta))}\,,
\tag9.28
$$ hence
$$
\aligned &\left(\Theta_M\cdot \chi_F\cdot
e^\nu\right)(gt\exp(\zeta)g^{-1})\ = \ 
\\ &\sum_{x\in X_M^{g\ttr g^{-1}}}\
\frac { e^\nu (gtg^{-1})\,\chi_F(gtg^{-1})c^M_{gtg^{-1},x}\,
e^{\lambda_{x}(g\zeta g^{-1}) -
\rho_{M,x}(g\zeta g^{-1})}} {\tsize
\prod_{\alpha
\in {\Phi}_M^{ +}} (1-e^{-\alpha_{x}})(gtg^{-1} \exp(g\zeta g^{-1}))}\,.
\endaligned
\tag9.29
$$ Unlike in (9.26),  $\Phi^+$ and $\Phi^+_M$ in (9.27-29) refer to the universal
positive root systems. The notation $c^M_{gtg^{-1},gx}$ is slightly
misleading, since $gtg^{-1}$ may lie in $M_\Bbb R\cdot A_\Bbb R \cdot F$
rather than $M_\Bbb R$ so we need to take the $M_\Bbb R$-component of
$gtg^{-1}$. Note that $\chi_F$ is constant on connected components, which
explains why we can omit the factor $\chi_F(g\zeta g^{-1})$. We should
remark that $\nu$ is used in two different but compatible senses. On the one
hand, as a linear function on $\fa$, on the other as a linear function on the
universal Cartan $\fh$, identified with the universal Cartan of $\fm\oplus\fa$,
so that $\fa\subset\fh$. When we use a fixed point $x$ of $\ttr$ in $X_M$ to
identify $\fh\cong\tilde \ft$, the restriction of $\lambda_x$ to
$\fa\subset\tilde \ft$ agrees with $\nu=\nu_x$. 

As $g$ ranges over $R/\tilde T_\Bbb R$ and $x$ over the fixed points of
$g\ttr g^{-1}$ in $X_M$, the translates $g^{-1}x$ range over the fixed points
of
$\ttr$ in $X$ which lie in $p^{-1}S_Y$, the inverse image of the closed orbit
$S_Y\subset Y$. Given such a fixed point $z=g^{-1}x$, the original fixed point
$x$ and $g\in\GR$ are determined by $z$ up to left translation, respectively
left multiplication, by an element of $\GR\cap MA = M_\Bbb R\cdot F\cdot\AR$.
Note also that $\lambda_{gx}(g\zeta g^{-1}) = \lambda_x(\zeta)$, with
similar identities for $\rho_{M,x}$, $\alpha_x$, etc.  Thus, substituting (9.29)
into the quantity on the right in (9.26), we find 
$$
\aligned &\left(I_{P_\Bbb R}^\GR\left(\Theta_M \cdot \chi_F
\cdot e^\nu\right)\  \tsize\prod_{\alpha\in\Phi^+}
|e^{\alpha/2}-e^{-\alpha/2}| \
\right)(t
\exp(\zeta)) \ = \ 
\\
\vspace{2.5\jot} &\sum_{x\in X^{\ttr}\cap p^{-1}S_Y}\ { e^\nu
(gtg^{-1})\,\chi_F(gtg^{-1})\,c^M_{gtg^{-1},gx}\, e^{\lambda_{x}(\zeta ) -
\rho_{M,x}(\zeta )}} \ \times
\\ &\qquad\times \ |e^{\rho_{M,x}}(t\exp(\zeta))|\prod_{\alpha
\in {\Phi}_M^{ +}}\frac { |1-e^{-\alpha_{x}}|}{ (1-e^{-\alpha_{x}})}(t \exp(\zeta
)) \,;
\endaligned
\tag9.30
$$ in this formula, the $g = g(x)\in\GR$ are chosen so that $g^{-1}x \in X_M$,
in which case it is unique up to left multiplication by some $m\in
M_\Bbb R\cdot F\cdot\AR$ -- that makes the various terms involving $g$ well
defined. We continue by substituting (9.27) into the left hand side of (9.30),
$$
\aligned &\sum_{x\in X^{\ttr}}\  { c_{t,x}\, e^{\lambda_x(\zeta) -
\rho_x(\zeta)}}  \ |e^{\rho_{x}}(t\exp(\zeta))|
\prod_{\alpha
\in {\Phi}^{ +}}\frac { |1-e^{-\alpha_{x}}|}{ (1-e^{-\alpha_{x}})}(t \exp(\zeta )) \
= \ 
\\
\vspace{2.5\jot} &\sum_{x\in X^{\ttr}\cap p^{-1}S_Y}\ { e^\nu
(gtg^{-1})\,\chi_F(gtg^{-1})\,c^M_{gtg^{-1},gx}\, e^{\lambda_{x}(\zeta ) -
\rho_{M,x}(\zeta )}} \ \times
\\ &\qquad\times \ |e^{\rho_{M,x}}(t\exp(\zeta))|\prod_{\alpha
\in {\Phi}_M^{ +}}\frac { |1-e^{-\alpha_{x}}|}{ (1-e^{-\alpha_{x}})}(t \exp(\zeta
)) \,.
\endaligned
\tag9.31
$$ Comparing the two sides in (9.31), we obtain the following re-statement of
the induced character formula (9.26) in terms of the $c_{t,x}$. First,
$$ c_{t,x} \ = \ 0 \ \ \ \text{unless} \ \ x\in p^{-1}S_Y\,;
\tag9.32a
$$ in words, the induced character is supported on the union of the
conjugates of the inducing subgroup $P_\Bbb R$. Secondly,
$$
\gathered c_{t,x} = (-1)^{k_x} e^\nu(gtg^{-1})
\frac{\chi_F(gtg^{-1})}{|e^{(\rho_x-\rho_{M,x})}(t)|} c^M_{gtg^{-1},gx} \ \ \
\text{if $x\in X^{\ttr}\cap p^{-1}S_Y$}\,,
\\
\vspace{2\jot}
\text{with $k_x=\#\{\alpha\in\Phi^+-\Phi^+_M \ | \ 0<e^{\alpha_x}(t)<1\}$\ \ 
and \ \ $g^{-1}x\in X_M$\,.}
\endgathered
\tag9.32b
$$ We should explain that non-real roots occur in pairs, so the quotient of the
products on the two sides of (9.31) is $\pm 1$, with the sign governed by the
number of positive real roots such that $e^{-\alpha_x}(t)>1$ , i.e., by the
number $k_x$.  Also,
$|e^{(\rho_x-\rho_{M,x})}(\exp\zeta)|=e^{(\rho_x-\rho_{M,x})}(\exp\zeta)$
since $\rho_x-\rho_{M,x}$ is real valued. 

To complete the proof, we must show that the relation (9.32) between the
$c_{t,x}$ and $c^M_{t,x}$ is consistent with one of the two formulas in (9.24)
when we express $c_{t,x}$ in terms of $\cf$ and $c^M_{t,x}$ in terms of
$\cf_M$. The sheaf $\cf$ is supported on $p^{-1}(S_Y)$ so any quantity
attached to $\cf$ by a local construction -- such as either of the two
Lefschetz numbers in (5.24) -- vanish at points not in   $p^{-1}(S_Y)$. Thus we
only need to check the consistency of (9.32b). Let us reduce the problem to
the case when the fixed point $x$ lies in $X_M$, so that we can choose $g=e$;
the legitimacy of this reduction reflects the fact that the choice of the base
point $y_0$ over which the reference fiber $X_M$ lies is arbitrary. 

After the reduction, we are dealing with a  point $x\in X_M$ fixed by the
Cartan subgroup $\tilde T_\Bbb R\subset M_\Bbb R\AR F$ and a regular element
$t\in \tilde T_\Bbb R$. The open Schubert cell $N^+(t,x)$ defined in (5.17)
splits naturally into a product
$$
\gathered V^+\times N^+_M(t,x) \ @>{\ \sim\ }>>\  N^+(t,x) \,, \ \ \ \ (v,\tilde
x)\mapsto v\cdot\tilde x\,,
\\
\text{with \ \ } V^+ \ = \ \exp(\oplus_{\alpha\in \Phi^+-\Phi_M^+}\
\fg^{\alpha_x})\,, \ \ \  N^+_M(t,x) \ = \  N^+(t,x)\cap X_M\,.
\endgathered
\tag9.33
$$  The group $V^+$ is the opposite to the unipotent radical $V$ of $P$; its
orbit
$V^+\cdot y_0$ is the open Schubert cell in $Y$ around $y_0$. Note that
$V^+_\Bbb R\cdot y_0$, the open Schubert cell in $S_Y\cong \GR/P_\Bbb R$,
coincides with the intersection of the complex Schubert cell $V^+\cdot y_0$
with the $\GR$-orbit $S_Y$. Thus the $S\cap N^+(t,x)$ splits into a product
$$ S\cap N^+(t,x) \ \cong \ V^+_\Bbb R \times (S_M\cap N^+_M(t,x))\,,
\tag9.34
$$ compatibly with (9.33). The choice of the base point $y_0$ allows us to
regard $A$ as a subgroup of the universal Cartan $H$, in fact, $H\cong A\times
H_M$. This splitting, too, is compatible with (9.33), as is the splitting
$$ N^+(t,x)\times H \ \cong \ (V^+\times A)\times (N^+_M(t,x)\times H_M)
\tag9.35
$$ of $N^+(t,x)\times H$, the inverse image of $N^+(t,x)$ in the enhanced flag
variety $\hat X$. The monodromy data defining the sheaf $\cf$ split into a
product corresponding to the splitting (9.35), of  the two characters
$$
\aligned
\chi_ A \cdot \chi_F \,&: \,A\cap\GR \ \longrightarrow \ \Bbb C^*\,,
\\
\chi_C^0\,&: \, C^0_\Bbb R  \ \longrightarrow \ \Bbb C^*\,;
\endaligned
\tag9.36
$$ recall that $A\cap \GR = \AR \cdot F$. The monodromic sheaf $\cf$,
regarded as a sheaf on  $N^+(t,x)\times H$, decomposes into an exterior
product of monodromic sheaves on  the two factors, with monodromy data
(9.36). According to (9.12), one of these is the
$(-\lambda_M-\rho_M)$-monodromic sheaf $\cf_M$, which corresponds to
the character $\chi_C^0$, though now we think of it as a  sheaf on 
$N^+_M(t,x)\times H_M$. The other is an $A\cap\GR$-equivariant,
$A$-monodromic sheaf $\cf_V$ on $V^+$, regarded as a sheaf on $V^+\times
A$.  This sheaf $\cf_V$ is the direct image of the
$A\cap\GR$-equivariant, $(-\lambda - \rho)|_\fa$-monodromic local system
$\Cal L_{\chi,V}$ on $V_\Bbb R^+$ associated to the character
$\chi|_{A\cap\GR} =
\chi_ A
\cdot
\chi_F$. Note that 
$$
d(\chi_A\cdot\chi_F) \ = \ d(\chi|_{A\cap\GR}) \ = \ (d\chi)|_\fa \ = \
\nu-(\rho_x-\rho_{M,x})\,,
\tag9.37
$$
provided $x\in X_M$ -- recall (9.11) and the definition of $\nu$ in (9.14).

Kashiwara's fixed point formalism is functorial with respect to exterior
products of sheaves. Let us argue that the subspaces $N'(t,x)$\,, $N''(t,x)$ of
$N^+(t,x)$ defined in (5.17) split compatibly with the product
decompositions (9.33-34); what makes this true is the fact that both are
orbits through $x$ of subgroups of $\exp(\fn^+(t,x))$ whose Lie algebras are
sums of root spaces -- any such group is the semidirect product of its
intersections with $V^+$ and $\exp(\fn_M^+(t,x))$. Recall the passage $\cf
\rightsquigarrow \cf(x)$ described in (5.19), with now $\tilde T_\Bbb R$
taking the place of $\TR$.  We use the analogous notation for $\cf_M$ and
$\cf_V$. The K\"unneth isomorphism
$$
\Cal H^*_{N'(t,x)}(\Bbb D\cf(x))_x \ \cong \ \Cal H^*_{V^+\cap
\exp \fn'(t,x)}(\Bbb D\cf_V(e))_e
\otimes \Cal H^*_{N'_M(t,x)}(\Bbb D\cf_M(x))_x
\tag9.38
$$ 
gives us the first Lefschetz number in theorem 5.24 as the product of the
two Lefschetz numbers
$$
\gathered
\tsize\sum (-1)^i \,\operatorname{tr}(\,\phi_t:\Cal H^i_{V'}(\Bbb D\cf_V(e))_e
\,
\to
\, 
\Cal H^i_{V'}(\Bbb D\cf_V(e))_e\otimes\Bbb C_{\nu}\,)\,,
\\
\tsize\sum (-1)^i \,\operatorname{tr}(\,\phi_t:\Cal H^i_{N'_M(t,x)}(\Bbb
D\cf_M(x))_x
\,
\to
\, 
\Cal H^i_{N'_M(t,x)}(\Bbb D\cf_M(x))_x\otimes\Bbb C_{\lambda_M}\,)\,,
\endgathered
\tag9.39
$$ with $V' = V^+\cap
\exp \fn'(t,x)$. Here $\Bbb C_\nu$ is the local system on $A$ with generating
section
$e^{\nu-\rho+\rho_M}$  -- recall (9.11) and the notation $\nu=\lambda|_\fa$
-- and $\Bbb C_{\lambda_M}$ the local system on $H_M$ with generating
section $e^{\lambda_M - \rho_M}$. Their exterior product, we remark, is the
local system $\Bbb C_\lambda$ on $H$. 

In the previous section we proved the special instance of theorem 5.24
involving the group $M_\Bbb R$ and the sheaf $\cf_M$. In effect, this
identifies the second of the two Lefschetz numbers as the constant
$c^M_{t,x}$ in (9.32b). The constant $c_{t,x}$ on the left in (9.23b)
corresponds to the induced character (9.27); because of proposition 9.14,
$c_{t,x}$ is the constant that corresponds to $\Theta(\cf)$. Thus, to
complete the proof, we must show:
$$
\gathered (-1)^{k_x} \,e^\nu(t)\,
\frac{\chi_F(t)}{|e^{(\rho_x-\rho_{M,x})}(t)|} \ = 
\\
\vspace{1.5\jot}
\tsize\sum (-1)^i \,\operatorname{tr}(\,\phi_t:\Cal H^i_{V'}(\Bbb D\cf_V(e))_e
\,
\to
\, 
\Cal H^i_{V'}(\Bbb D\cf_V(e))_e\otimes\Bbb C_{\nu}\,)\,;
\endgathered
\tag9.40
$$
as was pointed out, we are free to assume $g=e$ in (9.32b), so $g$ no longer
appears in (9.40). The first of the two Lefschetz numbers in (9.39) only
depends on the component of $t$ in $A\cap\GR$, so we assume $t\in
A\cap\GR$ from now on. 

The standard sheaf $\cf$ was constructed as the direct image of the twisted
local system $\Cal L_\chi$ on the $\GR$-orbit $S\subset X$. As was remarked
earlier,   $\cf_V$ is  the direct image of the local system $\Cal L_{\chi,V}$
under the closed embedding $V_\Bbb R^+ \hookrightarrow V^+$. The
operations of Verdier duality and local cohomology along a closed subspace
commute with direct image under a closed embedding. This allows us to
reduce (9.40) to a statement about  
the action of $t$ on the space $V^+_\Bbb R$ and the twisted local system $\Cal
L_{\chi,V}$\,:
$$
\gathered (-1)^{k_x} \,e^\nu(t)\,
\frac{\chi_F(t)}{|e^{(\rho_x-\rho_{M,x})}(t)|} \ = 
\\
\vspace{1.5\jot}
\tsize\sum (-1)^i \,\operatorname{tr}(\,\phi_t:\Cal H^i_{V_\Bbb R'}(\Bbb D\Cal
L_{\chi,V})_e
\,
\to
\, 
\Cal H^i_{V_\Bbb R'}(\Bbb D\Cal L_{\chi,V})_e\otimes\Bbb C_{\nu}\,)\,;
\endgathered
\tag9.41
$$
here $V_\Bbb R'$ denotes the intersection $V'\cap V_\Bbb R^+=
V^+_\Bbb R\cap\fn'(t,x)$ and $\phi_t$ the morphism induced by the $\phi$ 
corresponding to $\Bbb D \Cal L_{\chi,V}$. Let $\Cal L^*_{\chi,V} = \Cal
H\text{\it om}(\Cal L_{\chi,V}, \Bbb C_{V^+_\Bbb R})$ denote the local system
dual to $\Cal L_{\chi,V}$. Then 
$$
\Bbb D\Cal L_{\chi,V} \ = \ \Cal L^*_{\chi,V}\otimes \Bbb D_{V^+_\Bbb R}\,,
\tag9.42
$$
as equivariant twisted sheaves.  Correspondingly, the quantity on the right of
(9.41) becomes a product of two numbers. First, the action of $t$ via $\phi$ on
the stalk of $\Cal L^*_{\chi,V}$ at the fixed point $e$\,:  the structure of
$(A\cap\GR)$-equivariant,
$(-\lambda-\rho)|_\fa$-monodromic sheaf on $\Cal L^*_{\chi,V}$ gives us a
morphism $\phi: t^* \Cal L^*_{\chi,V} \to \Cal L^*_{\chi,V}\otimes \Bbb
C_\nu$, in complete analogy to (5.22); here, as in section 5, we think of
$t$ as variable in $A\cap\GR$. At the fixed point $e$, 
$\phi$ induces a morphism of stalks $\phi:  ({\Cal L^*_{\chi,V}})_e \to ({\Cal
L^*_{\chi,V}})_e\otimes \Bbb C_\nu$.  The sheaf $\,\Cal L^*_{\chi,V}$ has a one
dimensional  stalk at
$e$, and any section of $\Bbb C_\nu$ has a numerical value at any particular
$t\in A\cap\GR$; thus $t$ acts on $({\Cal L^*_{\chi,V}})_e$ by a scalar. We claim:
$$
\text{via}\ \phi \,,\ t\in A\cap \GR \ \ \text{acts on} \ \ ({\Cal L^*_{\chi,V}})_e \
\
\text{as the scalar}\ 
\  \,e^\nu(t)\,
\frac{\chi_F(t)}{|e^{(\rho_x-\rho_{M,x})}(t)|}\,.
\tag9.43
$$
Secondly, we apply  the fixed point formalism to the dualizing sheaf  $\,\Bbb
D_{V_\Bbb R^+}$. When we regard  $\Bbb
D_{V_\Bbb R^+}$ as an $A\cap
\GR$-equivariant sheaf, we get a morphism $\phi:t^*\Bbb D_{V_\Bbb R^+} \to
\Bbb D_{V_\Bbb R^+}$ -- this time to $\Bbb D_{V_\Bbb R^+}$ itself since the
twisting is trivial. The morphism $\phi$ induces a morphism $\phi_t$ on the
stalk of the local cohomology sheaf at the fixed point $e$, in analogy to (5.23).
We claim:
$$
\tsize\sum (-1)^i \,\operatorname{tr}(\,\phi_t:\Cal H^i_{V'}(\Bbb D_{V_\Bbb
R^+})_e
\,
\to
\, 
\Cal H^i_{V'}(\Bbb D_{V_\Bbb R^+})_e\,)\ 
= \ 
(-1)^{k_x}\,.
\tag9.44
$$
Together, (9.43-44) will imply (9.41),  so we only need to verify these two
assertions.

To clarify the reason for (9.43-44), let us look more generally at the case of a
Lie group $L$ acting transitively on a manifold $M$. An
$L$-equivariant local system $\Cal E$ on $M$ is determined by the datum of a
representation $\tau$ of the component group $\pi_0(L_m)$, of the isotropy
group $L_m$ at some $m\in M$. On the other hand, the formalism of
equivariant sheaves involves a distinguished isomorphism $\phi: a^*\Cal E \to
p^*\Cal E$ on $L\times M$. Since $L_m$ fixes $m$, both $\phi$ and its
inverse $\phi^{-1}: p^*\Cal E \to a^*\Cal E$ induce maps
$\phi, \phi^{-1} :L_m \to \operatorname{Aut}(\Cal E_m)$, which are  related to
$\tau$ by
$$
\phi^{-1}(\ell) \ = \ \tau(\ell)\,,\qquad \text{for} \ \ \ell \in L_m\,.
\tag9.45
$$
Indeed, $\phi_{\ell,m}$, for $(\ell,m) \in L\times M$,  maps $\Cal E_{\ell m}$ to
$\Cal E_m$, hence the cocycle condition on $\phi$ becomes
$$
\phi^{-1}_{\ell_2\ell_1,m} \ = \ \phi^{-1}_{\ell_2,\ell_1m} \circ
\phi^{-1}_{\ell_1,m}
\tag9.46
$$
whereas $\phi$ has the wrong variance to be multiplicative. The
correspondence (9.45) is obtained by restricting $\phi^{-1}$ to
$L_m\times\{m\}$. 

We apply this discussion to the case
$L=\GR\times\fh$ and $M = \hat S$, the inverse image in $\hat X$ of the orbit
$S$ to which we had associated the local system $\Cal L_{\chi}$. Since we are 
treating $\Cal L_{\chi}$ as an $\GR\times \fh$-equivariant sheaf, we shall
follow the conventions of (5.5-8) and use the symbol $\phi'$ for the map
which relates $\e^*a^*\Cal L_{\chi}^*$ to  $\e^*p^*\Cal L_{\chi}^*$\,:
$$
\phi': \e^*a^*\Cal L_{\chi}^* \to \e^*p^*\Cal L_{\chi}^*\,,
\tag9.47a
$$
as in (5.6a). Then, as was argued in section 5, the product
$e^{\lambda-\rho}\phi'$ descends to a map 
$$
\phi=e^{\lambda-\rho}\phi':  a^*\Cal L_{\chi}^* \to  p^*\Cal L_{\chi}^*\,.
\tag9.47b
$$
In the construction of the local system $\Cal L_\chi$ before the statement of
lemma 6.9 the action of the component group $\pi_0((\GR\times\fh)_{\hat
X_0})$ on the stalk of $\Cal L_\chi$ at $\hat x^0$ was written as $\psi$, and
was related to the character $\chi$ to which $\Cal L_\chi$ is associated by
the formula
$$
\psi = \tilde\chi e^{-(\lambda-\rho)}\,,
\tag9.48
$$
where $\tilde\chi$ is the lifting of $\chi$ from $\TR$ to $(\GR)_{x_0}$ -- cf.
(6.6a) and (6.8). According to (9.45), specialized to the present situation,
$\psi$ coincides with the inverse of the  character $\phi'$, but the one
 corresponding to $\Cal L_\chi$ rather than its dual as in (9.47a). Thus, with
$\phi'$ as in (9.47a), we get the equality $\psi = (\phi')^*$ ($=\phi'$, when
$\phi'$ is viewed as a character on the fixed point set), hence
$$
\chi = (\psi e^{\lambda-\rho})|_\TR  = (\phi'  e^{\lambda-\rho})|_\TR =
\phi|_\TR\,.
\tag9.49
$$
We restrict this identity further to $A\cap\GR= A\cap \TR$ and observe that
$$
\chi(t) = e^\nu(t)\,
\frac{\chi_F(t)}{|e^{(\rho_x-\rho_{M,x})}(t)|} \ \ \ \ \text{for $t\in
A\cap\GR$}\,;
\tag9.50
$$
note: $d\chi|_{\fa_\Bbb R} = (\lambda-\rho)|_{\fa_\Bbb R}=
\nu-(\rho-\rho_M)|_{\fa_\Bbb R}$\,, from which (9.50) follows when  $t\in
(A\cap\GR)^0$. But $\chi_F= \chi|_F$\,, so (9.50) follows for every $t\in
A\cap\GR$.  At this point, (9.49-50) imply (9.43).

The verification of (9.44) is simpler since the equivariant sheaf \,$\Bbb
D_{V_\Bbb R^+}$ is untwisted. In fact, (9.44) can be deduced
directly from \cite{KSa, 9.6.14}; however, the present instance of
Kashiwara-Schapira's result is elementary so we shall argue
directly. Concretely,
$$
\Cal H^*_{V'}(\Bbb D_{V_\Bbb R^+}) = or_{V'\cap  V_\Bbb R^+} [\dim_\Bbb R
V'\cap  V_\Bbb R^+]
\tag9.51
$$
is the orientation sheaf of $V'\cap  V_\Bbb R^+$ in degree $-(\dim_\Bbb R
V'\cap  V_\Bbb R^+)$. The exponential map maps $\fv^+_\Bbb R \cap \fn'(t,x)$
homeomorphically and $\ttr$-equivariantly onto $V'\cap V^+_\Bbb R$. Recall
the condition (5.15) on the subset $\Psi'\subset \Phi^+$. For
$\alpha\in\Psi'$, $\alpha\notin\Phi_M$, we have the following possibilities:
$$
\aligned
&\text{a)} \ \ \alpha\in \Psi'\,, \ \alpha = - \bar\alpha \ \implies \ 
\fg^\alpha\cap \fv^+_\Bbb R\cap\fn'(t,x) = 0 
\\
&\text{b)}  \ \ \alpha\in \Psi'\,, \ \alpha = \bar\alpha \ \implies \ 
\dim_\Bbb R(\fg^\alpha\cap \fv^+_\Bbb R\cap\fn'(t,x)) = 1
\\
&\text{c)} \ \ \alpha,\bar\alpha\in \Psi'\,, \ \alpha \neq \pm \bar\alpha \
\implies
\ 
\dim_\Bbb R((\fg^\alpha\oplus{\fg^{\bar\alpha}})\cap \fv^+_\Bbb
R\cap\fn'(t,x)) = 2
\\
&\text{d)} \ \ \alpha\in \Psi',\bar\alpha\notin\Psi'\,, \ \alpha \neq \pm
\bar\alpha \
\implies
\ 
(\fg^\alpha\oplus{\fg^{\bar\alpha}})\cap \fv^+_\Bbb
R\cap\fn'(t,x) = 0\,;
\endaligned
\tag9.52
$$
moreover, $\fv_\Bbb R^+\cap\fn'(t,x)$ splits into the direct sum of the
intersections appearing on the right in (9.52). Thus 
$$
\dim_\Bbb R(\fv_\Bbb R^+\cap\fn'(t,x)) \equiv \#\{\alpha\in\Psi'\mid
\alpha = \bar\alpha\}\pmod 2\,.
\tag9.53
$$
In the situation (9.52c), $t$ acts on $\fg^\alpha\oplus{\fg^{\bar\alpha}}$ with
eigenvalues $e^\alpha (t)$, $e^{\bar
\alpha}(t)$, and hence preserves the orientation of $(
\fg^\alpha\oplus{\fg^{\bar\alpha}})\cap \fv^+_\Bbb R\cap\fn'(t,x)$.  On the
other hand, in the situation (9.52b), $t$ preserves or reverses the orientation
of $\fg^\alpha\cap \fv^+_\Bbb R\cap\fn'(t,x)$ depending on the sign of
$e^\alpha(t)$. Hence
$$
\gathered
\text{$t$ preserves or reverses the orientation
of $\fv^+_\Bbb R\cap\fn'(t,x)$ }
\\
\text{depending on the parity of $\#\{\alpha\in\Psi'\mid
\alpha = \bar\alpha, e^\alpha(t)<0\}$}\,.
\endgathered
\tag9.54
$$
Finally, since complex roots in $\Phi-\Phi_M$ occur in pairs, and since $t$
must be regular, 
$$
\aligned
k_x &=  \#\{\alpha\in\Phi^+-\Phi^+_M\mid
0< e^\alpha(t)<1\} 
\\
&\equiv \#\{\alpha\in\Psi'\mid
\alpha = \bar\alpha, e^\alpha(t)>0\} \pmod 2\,.
\endaligned
\tag9.55
$$
In view of (9.51), $\phi_t:\Cal H^*_{V'}(\Bbb D_{V_\Bbb R^+})_e
\to
\Cal H^*_{V'}(\Bbb D_{V_\Bbb R^+})_e$ is the geometric action of $t$ on the
orientation of $\fv^+_\Bbb R\cap\fn'(t,x)$ (in degree $-\dim_\Bbb R
\fv^+_\Bbb R\cap\fn'(t,x)$\,).  At this point  (9.53-55) imply (9.44), and that
completes the proof.
\enddemo
\vskip1cm

\subheading{\bf 10. Intertwining Functors}
\vskip .5cm

So far, we have verified the main theorems for all standard sheaves
associated to maximally real orbits. In the present section we shall use the
mechanism of intertwining functors to extend the results to standard
sheaves attached to all orbits; as was argued in \S 6, that will complete our
proofs of theorems 3.8, 5.12, 5.24, and 5.27. 

We begin by recalling certain facts about the orbit structure. Let
$\TR\subset\GR$ be a Cartan subgroup, $x_0\in X$ a fixed point of $\TR$, and
$\tau_{x_0}:\ft\to\fh$ the corresponding isomorphism. As was mentioned in
\S 6, these data determine a $\GR$-orbit 
$$
S_0\  = \ S(\TR,\tau_{x_0})\,,
\tag10.1
$$
and the correspondence between pairs $(\TR,\tau_{x_0})$ becomes bijective
when $\TR$ is taken modulo $\GR$-conjugacy and, once $\TR$ is specified,
$\tau_{x_0}$ modulo $N_\GR(\TR)$-conjugacy. Let us consider a particular
orbit $S_0$ as in (10.1). The isomorphism $\tau_{x_0}$ identifies the
universal root system $\Phi$ with the concrete root system $\Phi(\fg,\ft)$.
We use this identification to transfer the universal positive root system
$\Phi^+$ to a  positive root system $\Phi^+_{S_0}(\fg,\ft)$ in  $\Phi(\fg,\ft)$;
this positive root system does depend only on the orbit, since $S_0$
determines 
$(\TR,\tau_{x_0})$ up to conjugacy.  Recall the notions of real, imaginary,
and complex root, which were defined in \S 9. The integer 
$$
c(S_0) \ = \ \#\{\,\alpha\in\Phi^+\mid\text{$\tau_{x_0}^*\alpha$ is
complex and
$\overline{\tau_{x_0}^*\alpha}\in -\Phi^+_{S_0}(\fg,\ft)$}\,\}
\tag10.2
$$
measures the extent to which   $S_0$ fails to be maximally real: if
$c(S_0)=0$, it is maximally real, so we do not need to deal with this case any
further. 

Let us suppose then that $c(S_0)>0$. According to (9.1), there must exist at
least  one  simple root $\alpha\in\Phi^+$ such that 
$\tau_{x_0}^*\alpha\in\Phi^+_{S_0}(\fg,\ft)$ is complex and has a negative 
complex conjugate
$\overline{\tau_{x_0}^*\alpha}\in-\Phi^+_{S_0}(\fg,\ft)$. Let us fix such a 
simple root
$\alpha$. It
determines a
$G$-equivariant fibration 
$$
X \ @>{\ \pi \ }>> \  X_\alpha \qquad\text{with fiber $\Bbb P^1$}
\tag10.3
$$
over the generalized flag variety $X_\alpha$. The Cartan subgroup $\TR$ has
exactly two fixed points on the fiber $\pi^{-1}(\pi(x_0))$; these
correspond to the fixed points of a Cartan subgroup of the group $SL(2,
\Bbb C)$ acting on $\Bbb P^1$. One of the two fixed points is 
$x_0$. We call the  second fixed point $x_1$. Let $S_1$ be the $\GR$-orbit
through $x_1$. Then 
$$
\aligned
\text{a)} \ \ &\pi^{-1}(\pi(S_0)) \ = \ S_0\cup S_1\,;
\\
\text{b)} \ \ &\tau_{x_1}\!\!\circ\tau_{x_0}^{-1}\ :\ \fh\to\fh \ \ \ \text{is
reflection about the root $\alpha$}\,;
\\
\text{c)} \ \ &c(S_1) \ = \ c(S_0) -1\,,
\endaligned
\tag10.4
$$
as is shown in \cite{HMSW2, S4}, for example. 

Let $\cf_0\in \operatorname D_\GR(X)_{-\lambda}$ be a standard sheaf
attached to the orbit $S_0$. We had argued in \S 6 that $\cf_0$ is the direct
image $R{j_0}_* \Cal L_{\chi_0}$ of an irreducible, $\GR$-equivariant local
system
$\Cal L_{\chi_0}$ on the orbit $S_0$. This local system corresponds to the
datum of a character 
$$
\chi_0\ : \  \TR \longrightarrow \Bbb C^*\,, \ \ \ \text{with} \ \ d\chi_0 =
\tau_{x_0}^*(\lambda-\rho)\,.
\tag10.5
$$
We define a new character $\chi_1$ by the formula
$$
\chi_1\ = \ e^{-\alpha_{x_1}}\,\chi_0\ : \ \TR \longrightarrow \Bbb C^*\,,
\tag10.6a
$$
where $\alpha_{x_1}= \tau_{x_1}^*\alpha$, as before.
Then, because of (10.4b), 
$$
d\chi_1\ =  \ - \tau_{x_1}^*\alpha +
\tau_{x_0}^*(\lambda-\rho) \ = \ 
\tau_{x_1}^*(-\alpha+s_\alpha(\lambda-\rho))\,,
$$
hence
$$
d\chi_1\ =  \ \tau_{x_1}^*(s_\alpha\lambda-\rho)\,.
\tag10.6b
$$
Since $\chi_1$ satisfies this condition, it determines an irreducible
$\GR$-equivariant local system $\Cal L_{\chi_1}$ on $S_1$, whose direct
image is a standard sheaf $\cf_1\in\operatorname
D_\GR(X)_{-s_\alpha\lambda}$ attached to the orbit $S_1$. We shall argue by
induction on the integer (10.2), and thus assume 
$$
\text{theorems 3.8, 5.12, 5.24, and 5.27 are satisfied by the sheaf}\,\ \cf_1\,.
\tag10.7
$$
We shall deduce the validity of these theorems for $\cf_0$, which will
then imply the theorems in full generality. 

The two standard sheaves $\cf_0\,,\,\cf_1$ are related geometrically by an
intertwining functor. Intertwining functors were introduced by
Beilinson-Bernstein \cite{BB2} in the context of $\Cal D$-modules. Their
formalism carries over readily to the setting of constructible
(untwisted) sheaves on $X$ -- see \cite{SV4}. In the case of twisted sheaves on
$X$, the construction of intertwining functors takes place on the generalized
flag variety; this was worked out by Beilinson-Bernstein in an earlier,
preprint version of \cite{BB3}. Here we shall discuss their construction, in
slightly modified form, in the case we need: the intertwining functor
corresponding to the simple root $\alpha$. 

Following the conventions in \cite{SV4}, we let $Y_\alpha$ denote the variety of
pairs $(x',x'')$ in $X\times X$ in relative position $s_\alpha$, and $p,q$ for the
natural projections 
$$
X \ @<{\ p \ }<< Y_\alpha @>{\ q \ }>> X
\tag10.8
$$
to the two factors. Both are fibrations with fiber $\Bbb C$. We shall put (10.8)
into a commutative diagram 
$$
\CD
\hat X \ @<{\ \hat p \ }<< \hat Y_\alpha @>{\ \hat q \ }>> \hat X
\\
@VVV @VVV @VVV
\\
X \ @<{\ p \ }<< Y_\alpha @>{\ q \ }>> X
\endCD
\tag10.9
$$
with certain properties that we shall explain next. Recall that the two outer
vertical arrows exhibit $\hat X$ as principal $H$-bundle over $X$. To begin
with, 
$$
\hat Y_\alpha \longrightarrow Y_\alpha \ \ \ \text{is a principal $H$-bundle}\,.
\tag10.10a
$$
Thus $H$ operates on the three spaces in the top row in (10.9). With respect to
to these actions, 
$$
\hat p(h\cdot \hat y) \ = \ h\cdot \hat p(\hat y) \ \ \ \text{and}\ \ \ \hat
q(h\cdot
\hat y)
\ =
\ s_\alpha(h)\cdot \hat q(\hat y)
\tag10.10b
$$
for all $\hat y\in \hat Y_\alpha$ and $h\in H$. Lastly, 
$$
\text{the two outer squares in (10.9) are Cartesian.}
\tag10.10c
$$
These are the formal properties that will matter to us. 

We shall not phrase our construction of the commutative diagram (10.9) in
invariant terms, though this could be done along the lines of the discussion of
the enhanced flag variety in \S 2. We make the identifications
$$
X\ = \ G/TN\,, \ \ \ \hat X \ = \ G/N\,,
\tag10.11a
$$
where $T$ is a concrete Cartan subgroup, $N$ a maximal unipotent subgroup
normalized by $T$. By letting $N$ correspond to the negative roots as usual,
we obtain an explicit identification 
$$
H \ \cong \ T
\tag10.11b
$$
between the concrete Cartan $T$ and the universal Cartan $H$. In the
construction of the diagram (10.9), we shall think of (10.11b) as an equality,
and accordingly shall identify the root system $\Phi(\fg,\ft)$ with the
universal root system $\Phi$.  

We choose a particular representative $s$ for the Weyl
reflection about the simple root $\alpha$, 
$$
s\in N_G(T)\,,\qquad \text{so that \, $\operatorname{Ad}s\,$ induces
$s_\alpha\in W$}\,,
\tag10.12
$$
and set 
$$
\fn^\alpha  \ = \ \fn\,\cap\, (\operatorname{Ad}s) \fn\ = \
\bigoplus_{\beta\in\Phi^+,\, \beta\neq\alpha}
\fg^{-\beta}\,, \
\ \ \ \ \ N^\alpha \ = \ \exp(\fn^\alpha)\,.
\tag10.13
$$
The Cartan subgroup $T$ normalizes $N^\alpha =  N\cap sNs^{-1}$. The group
$TN^\alpha$ is precisely the simultaneous stabilizer of the identity coset
and its $s$-translate in $G/TN = X$, so $G/TN^\alpha$ is the $G$-orbit through
$(eB,sB)\in G/TN\times G/TN= X\times X$ -- in other words, $G/TN^\alpha$
is the variety $Y_\alpha$ of pairs $(x',x'')$ in $X\times X$ in relative position
$s_\alpha$\,:
$$
Y_\alpha \ = \ G/TN^\alpha\,;
\tag10.14a
$$
in terms of this identification, the projections $p,q:Y_\alpha=G/TN^\alpha
\to G/TN = X$   are given by
$$
p(gTN^\alpha) \ = \ gTN \,, \ \ \ \ q(gTN^\alpha) \ = \ gsTN \,.
\tag10.14b
$$
At this point it is a simple matter to finish the construction of the
commutative diagram (10.9). We set
$$
\gathered
\hat Y_\alpha \ = \ G/N^\alpha
\\
\hat p(gN^\alpha) \ = \ gN \,, \ \ \ \ \hat q(gN^\alpha) \ = \ gsN \,.
\endgathered
\tag10.15
$$
The commutativity of the diagram and the properties (10.10a-c) are readily
verified.  

The Weyl reflection $s_\alpha$ induces an intertwining functor $I_\alpha$
from the (untwisted) bounded derived category of sheaves on $X$ to itself;
for the discussion of the untwisted case, we shall rely on \cite{SV4, \S 7}.
Analogously we define a twisted intertwining functor
$$
\hat I_\alpha \ : \ \operatorname D_\GR(X)_{-\lambda} \ \longrightarrow \
\operatorname D_\GR(X)_{-s_\alpha\lambda}\,, \ \ \ \ \hat I_\alpha(\cf) \ = \
R\hat q_*\hat p ^*(\cf)[1]\,.
\tag10.16
$$
By definition, $\hat I_\alpha$ operates on the category of $\GR\times
\fh$-equivariant sheaves. Because of (10.10b), it sends 
$(-\lambda-\rho)$-monodromic sheaves to
$s_\alpha(-\lambda-\rho)$-monodromic sheaves. But
$s_\alpha(-\lambda-\rho)$ differs from $-s_\alpha\lambda-\rho$ by the
root $\alpha$, and the monodromicity condition depends only on the
twisting parameter modulo the weight lattice. Thus $\hat I_\alpha$ affects
the twisting as claimed in (10.16). 

\proclaim{10.17 Lemma} The intertwining functor $\hat I_\alpha$
maps the standard sheaf $\cf_1$ to $\cf_0[1]$, the standard sheaf $\cf_0$
with a shift in degree.
\endproclaim

\demo{Proof} By equivariance, $\hat I_\alpha(\cf_1)$ is supported on a union
of  $\GR\times
H$-orbits in $\hat X$. Every such orbit is the inverse image $\hat S$ of a
unique $\GR$-orbit $S\subset X$. Let then $ j: \hat S \hookrightarrow \h
X$ be the inclusion of an orbit. We shall show 
$$
 j^!\hat I_\alpha(\cf_1) \ = \ 
\cases 0 \ \ \ &\text{if} \, \ \h S\neq \h S_0
\\
\Cal L_{\chi_0}[1] \ \ \ &\text{if}\, \  \h S=\h S_0\,.
\endcases
\tag10.18
$$
Let us argue first that (10.18) implies the conclusion of the lemma. We apply
Verdier duality, which reverses the roles of stars and shrieks: to establish
that
$\Cal G=\Bbb D
\hat I_\alpha(\cf_1)$ is the lower shriek extension (i.e., extension by zero) of
a sheaf $\Cal L$ on the orbit $\h S_0$, it manifestly suffices to check that it
restricts to
$\Cal L$ on $S_0$, and that its restriction to all the other orbits is zero.  

To verify (10.18), we put the inclusions  $ j: \hat S \hookrightarrow \h
X$,  $ j_1: \hat S_1 \hookrightarrow \h
X$ into the commutative diagram
$$
\CD
\h S_1 @<{\h p_1}<< \h p^{-1}(\h S_1) @<{\tilde j}<< \h p^{-1}(\h S_1)\cap
\h q ^{-1}(\h S) @>{\h q_S}>> \h S
\\
@V{j_1}V{\cap}V @V{\tilde j_1}V{\cap}V  @V{\tilde j_1}V{\cap}V 
@VjV{\cap}V
\\
\h X @<{\h p}<< \h Y_\alpha @=  \h Y_\alpha @>{\h q}>> \h X\ .
\endCD
\tag10.19
$$
We need to identify
$$
 j^!\hat I_\alpha(\cf_1) \ = \ j^!R\h q_*\h p^*{Rj_1}_*(\Cal L_{\chi_1})[1]\,.
 \tag10.20a
$$
First we apply base change in the square on the left:
$$
\h p^*{Rj_1}_*(\Cal L_{\chi_1}) \ = \ {R\tilde j_1}_*\h p^*_1(\Cal
L_{\chi_1})\,.
\tag10.20b
$$
Next, we consider the subdiagram with arrows $\tilde j$, $\h q \circ\tilde
j_1$, $\h q_S$, and $j$, which is a Cartesian square. Thus, by base change,
$$
j^!R(\h q \circ\tilde j_1)_*(\h p_1^*\Cal L_{\chi_1})\ = \ {R{\h q_S}}_*\tilde
j^!(\h p_1^*\Cal L_{\chi_1})\,.
\tag10.20c
$$
Combining (10.20a-c), we find
$$
j^!\hat I_\alpha(\cf_1) \ = \ {R{\h q_S}}_*\tilde
j^!(\h p_1^*\Cal L_{\chi_1})[1]\,.
\tag10.21
$$
We claim:
$$
\h q \ : \ \h p^{-1}\h S_1 @>{\ \ \sim \ \ }>> \h S_0 \ \ \ \text{is a $(\GR\times
H)$-equivariant isomorphism}\,,
\tag10.22
$$
relative to the natural action of $\GR$ on the two spaces, the natural action
of $H$ on $\h p^{-1}\h S_1$, and the action of $H$ on $\h S_0$ obtained by
composing the natural action with the Weyl reflection $s_\alpha$ -- cf.
(10.10b).  Let us assume this for the moment. 

First suppose that the $\GR$-orbit $S$ is unequal to $S_0$. Then, by (10.22),
the intersection $\h p^{-1}(\h S_1)\cap
\h q ^{-1}(\h S)$ is empty. This intersection appears in the diagram (10.19).
Hence, in (10.21), the right hand side is the direct image of a sheaf on the
empty set -- the zero sheaf. Conclusion: $j^!\hat I_\alpha(\cf_1) = 0$ when
$S\neq S_0$. 

Now we suppose that $S=S_0$. Then, in (10.19), $\h p^{-1}(\h S_1)\cap
\h q ^{-1}(\h S_0)$ coincides with $\h p^{-1}(\h S_1)$, $\tilde j$ reduces to
the identity, and $\h q_S$ is the isomorphism (10.22). Thus, in the present
situation,
$$
j^!\hat I_\alpha(\cf_1) \ = \ \hat q_{S*} \h p_1^*\Cal L_{\chi_1}[1]
\tag10.23
$$
is -- up to a shift in degree -- the pullback to $\h S_0\simeq \h p^{-1}(\h
S_1)$ of a
$(\GR\times
\fh)$-equivariant local system under the $(\GR\times
\fh)$-equivariant map $\h p_1$. As explained in section 6, these equivariant
local systems are specified by characters of the Cartan subgroup $\TR$,
which fixes the two base points $x_0,x_1$. The correspondence
between equivariant local systems and characters involves the
identification between the complexification $T$ of $\TR$ and the universal
Cartan $H$. The two identifications $T\cong H$ determined by the base points
$x_0,x_1$ are related by the Weyl reflection $s_\alpha$. On the other hand,
the isomorphism $\h S_0\simeq \h p^{-1}(\h
S_1)$ is $H$-equivariant only when the $H$-action on one of the two spaces is
twisted by $s_\alpha$. These two occurrences of $s_\alpha$ cancel. It follows
that $j^!\hat I_\alpha(\cf_1)[-1]$ is the equivariant local system attached to
the character $\chi_1= e^{-\alpha}\chi_0$; note: our labeling of the twists
involves the shift by $\rho$, and $s_\alpha\rho = \rho -\alpha$ -- this
accounts for the factor $e^{-\alpha}$ in the relation between $\chi_0$ and
$\chi_1$. 

This completes the proof of the lemma, except the verification of (10.22).
The orbits $\h S_0$, $\h S_1$ fiber as principal $H$-bundles over the
$\GR$-orbits $S_0$ and $S_1$, respectively. Thus (10.22) follows formally
from (10.10b) and 
$$
 q \ : \  p^{-1} S_1 @>{\ \ \sim \ \ }>> S_0 \ \ \ \text{is a $\GR$-equivariant
isomorphism}\,.
\tag10.24
$$
This assertion, in turn, reduces to the following statement: let $F\cong
\Bbb P^1$ be the fiber of the fibration (10.3) which contains $x_0$ and $x_1$;
then 
$$
\aligned
\text{a)} \ \ &F\cap S_1 \ = \ \{x_1\}\,,
\\
\text{b)} \ \ &F\cap S_0 \ = \ F-\{x_1\}\,.
\endaligned
\tag10.25
$$
That, in effect, is a paraphrase of (10.4).
\enddemo
 
\proclaim{10.26 Lemma} $\Theta(\cf_1) \ = \ -\Theta(\cf_0)$.
\endproclaim
 This fact is established, explicitly or implicitly, in various places. See for
example \cite{HMSW2,S4}. 

The discussion in this section up to this point is relevant to the proof both
of the integral formula (3.8) and the fixed point theorems 5.12, 5.24, and
5.27. We now turn specifically to the integral formula. In \cite{SV4} we
remarked that the (untwisted) intertwining functors on the bounded derived 
category of (semi-algebraically) constructible sheaves,
$$
I_w \,:\, \operatorname D^b_c(X) \ \longrightarrow \ \operatorname
D^b_c(X)\,,\qquad w\in W\,, 
\tag10.27
$$
induce an action of $W$ on $\Cal L^+(X)$, the group of $\Bbb R^+$-conical,
semi-algebraic Lagrangian cycles in $\ct$ so that the map
$\occ:\operatorname D^b_c(X)\to \Cal L^+(X)$ becomes
$W$-equivariant. 

\proclaim{10.28 Lemma} $\occ(\h I_\alpha\cf)= I_\alpha\occ(\cf)$, for 
any $\lambda\in\fh^*$  and $\cf\in
\operatorname D_\GR(X)_{\lambda}$ .
\endproclaim
\demo{Proof}We think of $\cf$ as an element of the $K$-group
of $\operatorname D^b_c(Sh_{X,\lambda})$. This $K$-group is generated by
locally constant twisted sheaves on contractible closed subsets of $X$. For
sheaves of this special type, the assertion of the lemma follows from the
properties of the diagram (10.9).
\enddemo

Recall the definition of the group $\ohf {2n}(T^*_\GR X,\Bbb Z)$, which contains
the characteris\-tic cycle  $\occ(\cf_1)$, and of the differential forms
$\sigma$, $\tau_\lambda$ in section 3. 

\proclaim{10.29 Lemma} For all $\phi \in C^\infty_c(\gr)$ and $\lambda \in
\fh^*$, 
$$
\int_{I_\alpha \occ(\cf_1)} \mu^*_ \lambda\hat \phi \ (-\sigma + \pi^*
\tau_ \lambda)^n
\ =
\
\int_{\occ(\cf_1)}
\mu^*_{s_\alpha\lambda} \hat \phi \ (-\sigma + \pi^*
\tau_{s_\alpha\lambda})^n
$$
\endproclaim
Before giving the proof, let us remark that lemmas 10.17, 10.28, and
10.29 imply 
$$
\gathered
 \frac 1 {(2\pi i)^n n!}\int_{\occ(\cf_0)}
\mu^*_\lambda \hat \phi \ (-\sigma + \pi^* \tau_\lambda)^n\ = \\
\qquad\qquad - \frac 1 {(2\pi i)^n n!}\int_{\occ(\cf_1)}
\mu^*_{s_\alpha\lambda} \hat \phi \ (-\sigma + \pi^*
\tau_{s_\alpha\lambda} )^n\,,
\endgathered
\tag10.30
$$
for any $\phi\in C^\infty_c(\gr)$. The induction hypothesis identifies the right
hand side of (10.30) as the integral of $\phi$ against the distribution
$-\theta(\cf_1)$. Thus, by 10.26, 
$$
 \frac 1 {(2\pi i)^n n!}\int_{\occ(\cf_0)}
\mu^*_\lambda \hat \phi \ (-\sigma + \pi^* \tau_\lambda)^n\ = \
\int_{\gr}\theta(\cf_0)\,\phi\,dx \,.
$$
This is the assertion of theorem 3.8 for the sheaf $\cf_0$, thus completing our
inductive proof of the theorem.

\demo{Proof of lemma 10.29} By proposition 3.7, we may as well assume that
$\lambda$ is regular. According to \cite{SV4, theorem 9.1}, the $W$-action
(10.27) agrees with the geometric action described in \cite{SV4, \S 8},
which was originally defined by Rossmann \cite{R3}. The path 
$$
\gamma(t) \ = \ ts_\alpha\lambda\,,\qquad 0\leq t\leq 1
\tag10.31
$$
satisfies the hypothesis of \cite{SV4, lemma 9.4}. Thus
$$
I_\alpha\occ(\cf_1) \, = \,  \lim_{t\to 0^+}
(\mu_{t\lambda}^{-1}\circ\mu_{ts_\alpha\lambda})\occ(\cf_1)\, = \, -\partial C
+(\mu_{\lambda}^{-1}\circ\mu_{s_\alpha\lambda})\occ(\cf_1)\,,
\tag10.32a
$$
where $C$ is the semi-algebraic chain 
$$
C \ = \ \{\,(\mu_{t\lambda}^{-1}\circ\mu_{ts_\alpha\lambda})\occ(\cf_1)\mid
0\leq t \leq 1\,\}\,;
\tag10.32b
$$
for the relation between the boundary of the chain (10.32b) and the limit
(10.32a) and the convention for orienting the chain, see \cite{SV4, \S 3}. The
chain $C$ satisfies the hypothesis of lemma 3.19 because the curve (10.23) is
compact. We conclude
$$
\gathered
\int_{I_\alpha \occ(\cf_1)} \mu^*_ \lambda\hat \phi \ (-\sigma + \pi^*
\tau_ \lambda)^n
\ = \ 
\int_{(\mu_{\lambda}^{-1}\circ\mu_{s_\alpha\lambda})\occ(\cf_1)} \mu^*_
\lambda\hat \phi \ (-\sigma + \pi^*
\tau_ \lambda)^n
\\
=\ \int_{\occ(\cf_1)} \mu^*_{
s_\alpha\lambda}\hat \phi \
(\mu_{\lambda}^{-1}\circ\mu_{s_\alpha\lambda})^*(-\sigma + \pi^*
\tau_ \lambda)^n\,.
\endgathered
\tag10.33
$$
The complex coadjoint orbit $\Omega_\lambda$ and the canonical symplectic
form $\sigma_\lambda$ on it depend only on the $W$-orbit of $\lambda$, not
on $\lambda$ itself. Hence, by proposition 3.3,
$$
\mu_\lambda^*\sigma_\lambda \ = \ -\sigma + \pi^*\tau_\lambda\ \ \ 
\text{and}\ \ \   \mu_{s_\alpha\lambda}^*\sigma_\lambda \ = \ -\sigma +
\pi^*\tau_{s_\alpha\lambda}\,,
\tag10.34
$$
which implies $(\mu_{\lambda}^{-1}\circ\mu_{s_\alpha\lambda})^*
(-\sigma + \pi^*\tau_ \lambda) = -\sigma + \pi^*\tau_ {s_\alpha\lambda}$. In
view of (10.33), this completes the proof of lemma 10.29.
\enddemo

Let us carry out the induction step for the fixed point formulas (5.24); we had
argued earlier that (5.24) implies also theorems 5.12 and 5.27. According to
the induction hypothesis (10.7), we know the assertion of (5.24) for the
particular sheaf $\cf_1$. We must establish it also for $\cf_0 =
\hat I_\alpha\cf_1[-1]$. As in the statement of (5.24) we fix a regular
semisimple element $t\in \GR'$ and a fixed point $x\in X$ of $t$. Then $x$ is a
fixed point also of the Cartan subgroup $\ttr= N_\GR(t)$. Let us denote the
sections
$c_{t,x}$ corresponding to $\cf_0,\cf_1$ by $c_{t,x}(\cf_0), c_{t,x}(\cf_1)$.
Note that $\cf_0\in\operatorname D_\GR(X)_{-\lambda}$, whereas
$\cf_1\in\operatorname D_\GR(X)_{-s_\alpha\lambda}$. Hence
$c_{t,x}(\cf_0)$ is a (local) section of $\Bbb C_\lambda$ -- i.e., a multiple of a
branch of the multiple valued function $e^{\lambda_x-\rho_x}$ -- and $
c_{t,x}(\cf_1)$ is a section of $\Bbb C_{s_\alpha\lambda}$  -- i.e., a multiple of a
branch of the multiple valued function
$e^{s_\alpha\lambda_x-\rho_x}$. Comparing the local expressions (4.11) for
the two invariant eigendistributions $\Theta(\cf_0)$ and $\Theta(\cf_1)$, and
taking into account lemma 10.26, we find 
$$
e^{\alpha_x}c_{t,x}(\cf_0) \ = \  c_{t,s_\alpha x}(\cf_1) \qquad
\text{for all}\ \  x\in X^\ttr\,.
\tag10.35
$$
Here $s_\alpha$ is shorthand notation for the point obtained as follows: we
identify the abstract Weyl group $W$ with the concrete Weyl group
$W(G,\tilde T)$ via
$\tau_x$; via this identification, $W$ permutes the fixed points of $\tilde T$.

Theorem 5.24 expresses the two constants
$c_{t,x}(\cf_0)$ and
$ c_{t,s_\alpha x}(\cf_1) $ in (10.35) as fixed point expressions applied not to the
sheaves
$\cf_0$ and
$\cf_1$ themselves, but rather their duals $\cg_0 = \Bbb D \cf_0$, $\cg_1 =
\Bbb D \cf_1$.  Let us denote the fixed point expressions in the second line of
5.24 by 
$$
\gathered
d_{t,x}(\cg_0)=\tsize\sum (-1)^i \operatorname{tr}(\phi_t:\Cal
H_{\{x\}}^i(\cg_0(x)|_{N''(t,x)})
\to
\Cal H_{\{x\}}^i(\cg_0(x)|_{N''(t,x)})\otimes\Bbb C_\lambda)
\\
d_{t,x}(\cg_1)=\tsize\sum (-1)^i \,\operatorname{tr}(\phi_t:\Cal
H_{\{x\}}^i(\cg_1(x)|_{N''(t,x)})
\to
\Cal H_{\{x\}}^i(\cg_1(x)|_{N''(t,x)})\otimes\Bbb C_{s_\alpha\lambda})\,;
\endgathered
$$
note that these constants have a different meaning from the $d_{t,x}$ in section
5. Our induction hypothesis is contained in the equation 
$$
d_{t,x}(\cg_1)=c_{t,x}(\cf_1)
\tag10.36a
$$
for all  regular $t$ and all fixed points $x$ of $t$. The conclusion we want is 
$$
d_{t,x}(\cg_0)=c_{t,x}(\cf_0)\,,
\tag10.36b
$$
again for all $t$ and $x$. Dualizing the relation $\hat I_\alpha\cf_1 = \cf_0[1]$
-- cf. (10.16)  and lemma 10.17 -- we find 
$$
\hat J_\alpha\cg_1 = \cg_0[-1]\,,\qquad \text{where} \ \ \hat J_\alpha=R\h
q_!\h p ^*[1]\,.
\tag10.37
$$
In view of (10.35) and (10.37), the inductive conclusion (10.36b) will follow
from the identity
$$
e^{\alpha_x}d_{t,x}(\h J_\alpha\cg) \ = \  -d_{t,s_\alpha x}(\cg) 
\tag10.38
$$
for $\cg=\cg_1$; the shift by $[-1]$ accounts for the minus sign. 
Our remaining task is to verify (10.38) for an arbitrary twisted sheaf
$\cg\in\operatorname D_\GR(X)_{\lambda}$\,.

\proclaim{10.39 Lemma} In the $K$-group $K(\operatorname
D_\GR(X)_{\lambda})$\,, the  square  $\h J_\alpha\circ \h J_\alpha$ of the
operator $\h J_\alpha$ coincides with the identity.
\endproclaim

In the untwisted case this follows from \cite{SV4, (7.14)}. The argument 
there does not obviously carry over to the twisted case, so we argue ab initio.

\demo{Proof} With $\h p$ and $\h q$ as in (10.9), $\h J_\alpha\circ \h J_\alpha
= R\h q_!\h p^* R\h q_!\h p^*[2]$. Applying base change in the fiber square
$$
\CD
\h Z @>{\Hat{\Hat q}}>> \h Y_\alpha
\\
@V{\Hat{\Hat p}}VV @V{\h p}VV
\\
\h Y_\alpha @>{\h q}>> \h X\,,
\endCD
\tag10.40
$$
we can rewrite the previous identity as $\h J_\alpha\circ \h J_\alpha =  R\h
q_!R{\Hat{\Hat q}}_!{\Hat{\Hat p}}^*\h p^*[2]$. Let $\h\pi_1$,$\h\pi_2$ denote
the projections from $\h X\times\h X$ to the two factors, and $\pi :\h Z \to \h
X\times\h X$ given by $\h \pi (\h z) = (\h p({\Hat{\Hat p}}(\h z)), \h q({\Hat{\Hat
q}}(\h z)))$.  Then
$$
\h J_\alpha\circ \h J_\alpha =  R\h
q_!R{\Hat{\Hat q}}_!{\Hat{\Hat p}}^*\h p^*[2] =
R\h\pi_{2!}R\h\pi_!\h\pi^*\h\pi_1^*[2]\,.
\tag10.41
$$
The projection formula for $\pi :\h Z \to \h
X\times\h X$ asserts that $R\h\pi_!\h\pi^*\Cal S =  (R\h\pi_!\Bbb C_{\h
Z})\otimes \Cal S$, for any sheaf $\Cal S$ on $\h X\times\h X$. We use this
formula to transform (10.41) into
$$
(\h J_\alpha\circ \h J_\alpha)(\cg) =
R\h\pi_{2!}\left((R\h\pi_!\Bbb C_{\h Z})\otimes\h\pi_1^*(\cg)\right)[2]\,.
\tag10.42
$$
The geometry of the various fibrations gives
$$
\h\pi^{-1}(\h x _1,\h x_2) = \cases \, \Bbb C &\ \ \text{if} \ \ \h x_1=\h x_2
\\
\,\Bbb C^* &\ \ \text{if $\ \ \h x_1$ and $\h x_2$ are in position $\alpha$}
\\
\,\emptyset  &\ \ \text{otherwise}\,.
\endcases
\tag10.43
$$
In the $K$-group, this gives the equality $R\h\pi_!\Bbb C_{\h Z}[2] = \Bbb
C_{\Delta \h X}$\,, and this equality persists even in the $K$-group of the
$\GR$-equivariant derived category. Thus $(\h J_\alpha\circ \h J_\alpha)(\cg) =
R\h\pi_{2!}\left((R\h\pi_!\Bbb C_{\h
Z})\otimes\h\pi_1^*(\cg)\right)[2]=R\h\pi_{2!}\left(\Bbb
C_{\Delta \h X}\otimes\h\pi_1^*(\cg)\right) = \cg$ in $K(\operatorname
D_\GR(X)_{\lambda})$\,.
\enddemo

Lemma 10.39 allows us to treat the points $x$ and $s_\alpha x$ in a
symmetric fashion. Thus, without loss of generality, we assume 
$$
|\alpha_x(t)| \leq 1\,.
\tag10.44
$$
Recall:  the submanifolds $N''(t,x)$, $N''(t,s_\alpha x)$ depended on the choice of
$\Psi''\subset \Phi^+$. At $x$, we make the minimal choice 
$$
\Psi''(x) = \{\beta\in\Phi^+ \mid |\beta_x(t)| >1\}\,,
\tag10.45a
$$
and at $s_\alpha x$, we also include $\alpha$ itself even if (10.44) is an equality,
i.e., 
$$
\Psi''(s_\alpha x) = \Psi''(x)\cup\{\alpha\}\,.
\tag10.45b
$$
Because $\alpha$ is simple, and because of (10.44), this choice is consistent with
(5.15b). 

In defining the morphism (5.22), we think of $t\in \tilde T_\Bbb R$ as variable,
and this is the reason for the appearance of the factor $\Bbb C_\lambda$. Now
that we regard $t$ as fixed, we can evaluate sections of $\Bbb C_\lambda$ --
which was constructed as a subsheaf of $\Cal O_{\tilde G}$ -- at  $t$, giving us
morphisms
$$
\phi(\cg,s_\alpha x) : t^*\cg(s_\alpha x) \to \cg(s_\alpha x)\,,\qquad \phi(\h
J_\alpha\cg, x): t^*(\h J_\alpha\cg)( x) \to(\h J_\alpha\cg)( x)\,.
\tag10.46
$$
The two sheaves are defined on different (open)  Schubert cells, so we cannot
relate the two morphisms directly. Our choices (10.45) imply that 
$N''(t,s_\alpha x)$ fibers over  $N''(t,x)$ with fiber $\Bbb C$; for each $z\in
N''(t,x)$, the inverse image in $N''(t,s_\alpha x)$ consists of all points in relative
position $\alpha$ relative to $z$. Differently put, in the diagram (10.9), $p$ maps
the inverse image $q^{-1}(N''(t,x))$ isomorphically onto $N''(t,s_\alpha x)$.
Thus, for economy of notation, we write the fibration as $q:N''(t,s_\alpha x)\to
N''(t,x)$. 
\proclaim{10.47 Lemma} The sheaves $(\h J_\alpha\cg)( x)|_{ N''(t,x)}$ and
$Rq_!(\cg(s_\alpha  x)|_{ N''(t, s_\alpha x)})[1]$ are ca\-nonically isomorphic, so
that
$$
e^{\alpha_x}(t)\phi(\h J_\alpha\cg, x)|_{ N''(t,x)} =
Rq_!\left(\phi(\cg,s_\alpha x)|_{N''(t,s_\alpha x)}\right)[1].
$$
\endproclaim

\demo{Proof} We consider the commutative diagram
$$
\CD
\h X @<{\h p}<< \h Y_\alpha @>{\h q}>> \h X
\\
@A{\h j_{s_\alpha x}}AA @AAA @AA{\h j_x}A
\\
\h N''(t,s_\alpha x)@<{\sim}<< \h q^{-1}\h N''(t,x)@>{\tilde q}>>\h N''(t,x)\,.
\endCD
\tag10.48
$$
By base change in the right square, 
$$
\h j^*_x\h J_\alpha \cg = \h j^*_xR\h q_!\h p^* \cg [1] \simeq R\tilde q_!\h
j^*_{s_\alpha x}\cg[1]
\tag10.49
$$
in $\operatorname D_{\ttr}( N''(t,x))_{s_\alpha\lambda}$\,. In section 5, we
constructed the morphisms (10.46) from the $\ttr$-equivariant structure of
the two sheaves on the two open Schubert cells. Restricting these morphisms to
(the $\ttr$-invariant subsets)  $N''(t,s_\alpha x)$ and $N''(t,x)$ is equivalent to
applying the construction directly to the restricted sheaves $j^*_{s_\alpha
x}\cg$ and $j^*_x\h J_\alpha \cg$. Thus (10.49) relates the restricted
morphisms. The construction of the first of the two morphisms (10.46)
involves multiplication by $e^{\lambda_{s_\alpha x}-\rho_{s_\alpha x}}$, that of
the second, multiplication by $e^{(s_\alpha\lambda)_x-\rho_x}$ -- in both cases,
these expressions are considered as functions on $\fh\cong\tilde\ft$; cf. (5.7-8).
This accounts for the factor $e^{\alpha_x}(t)$ in the statement of the lemma:
the value, at $t$, of  the quotient of the two exponential expressions, 
$$
e^{(s_\alpha\lambda)_x-\rho_x-(\lambda_{s_\alpha x}-\rho_{s_\alpha x})} =
e^{\alpha_x}\,,
$$
which is well defined on all of $\ttr$.
\enddemo

Let us reformulate the remaining problem, stripping away the inessential
aspects of our situation. We consider a complex vector space $V$ -- in our case,
$N''(t, x)$ identified with its Lie algebra via $\exp$ -- and  a linear
transformation 
$$
f: \Bbb C\times V \to \Bbb C\times V\,, \ \ \ f(\Bbb C) \subset \Bbb C\,, \ \ \
f(V)\subset V\,,
\tag10.50a
$$
such that 
$$
\aligned
\text{a)}\ \ &\text{all eigenvalues $\zeta$ of $f$ on $V$ satisfy $|\zeta|>1$\,,}
\\
\text{b)}\ \ &\text{the eigenvalue  of $f$ on $\Bbb C$ has absolute value
$\geq1$\,.}
\endaligned
\tag10.50b
$$
In the situation at hand, $\Bbb C\times V$ is $N''(t,s_\alpha x)$, again identified
with its Lie algebra, and $f$ is the map induced by $t$; the hypotheses (10.50b)
are satisfied because of (10.45). Next, we consider a constructible sheaf $\Cal
E\in \db{\Bbb C\times V}$ together with the datum of a morphism
$$
\phi: f^*\Cal E \to \Cal E\,,
\tag10.51
$$
namely $\Cal E = \cg(s_\alpha  x)|_{ N''(t, s_\alpha x)}$ and $\phi=
\phi(\cg,s_\alpha x)|_{N''(t,s_\alpha x)}$\,. Let $q:\Bbb C\times V\to V$ be the
projection. Then $Rq_!\Cal E\in \db{ V}$ and $Rq_!\phi: f^*(Rq_!\Cal E)
\to Rq_!\Cal E$. 

\proclaim{10.52 Lemma} Under the hypotheses just mentioned,
$$
\tsize\sum(-1)^i\operatorname{tr}(Rq_!\phi: \!\oh^i_{\{0\}}(Rq_!\Cal E) \to
\oh^i_{\{0\}}(Rq_!\Cal E)) = \sum(-1)^i\operatorname{tr}(\phi:\!
\oh^i_{\{0\}}(\Cal E) \to
\oh^i_{\{0\}}(\Cal E)).
$$
\endproclaim

We not that this completes the induction step for the fixed point theorems
5.12, 5.24, and 5.27: the present lemma, combined with lemma 10.47 implies
the identity (10.38).

\demo{Proof of 10.52} We compactify $q$ to $\bar q : \Bbb P^1\times V
\to V$ and identify $V \cong \{\infty\}\times V$. We write $j:\Bbb C \times V
\hookrightarrow   \Bbb P^1\times V$ for the open embedding, $\bar k: V\cong
\{\infty\}\times V
\hookrightarrow \Bbb P^1\times V$ for the closed embedding,
$\bar i:\Bbb P^1
\cong \bar q^{-1}(0)
\hookrightarrow  \Bbb P^1 \times V$ for the natural inclusion, and $i$ for the
restricted map $i: \Bbb C \hookrightarrow \Bbb C\times V$.  Thus, 
$$
Rq_!\Cal E = R\bar q_!j_!\Cal E = R\bar q_*j_!\Cal E,
\tag10.53
$$
and, by applying base change to the projection $\bar q$, we get
$$
\gathered
\tsize\sum(-1)^i\operatorname{tr}(Rq_!\phi: \!\oh^i_{\{0\}}(Rq_!\Cal E) \to
\oh^i_{\{0\}}(Rq_!\Cal E)) \ = \ 
\\
\tsize\sum(-1)^i\operatorname{tr}(\phi:
\oh^i(\Bbb P^1, \bar i^!j_!\Cal E) \to
\oh^i(\Bbb P^1, \bar i^!j_!\Cal E))\,.
\endgathered
\tag10.54
$$
The Grothendieck-Lefschetz fixed point formula expresses the right hand side
of (10.54) as a sum of two terms associated to the fixed points 0 and $\infty$ of
$f$ on $\Bbb P^1 \cong \bar q^{-1}(0)$; see \cite{KSa, \S 9.6} for an exposition of
this formalism. We use the fixed point formula of
\cite{GM} to compute these contributions.  We will deal with the point 0 first. 
According to (10.50b) we can view  the map $f$ (restricted to $\Bbb P^1$) as
expanding around 0.  Theorem 4.7 in \cite{GM} gives  the fixed point 
contribution from 0 in five different ways. One of these, involving the complex
which Goresky-MacPherson call $\bold A_4$, \lq\lq shrieks" the sheaf $ i^!j_!\Cal
E$ to the contacting directions -- in our case, this means taking local
cohomology at $\{0\}$  -- then restricts the resulting sheaf to the fixed point --
which, in our case, is a vacuous operation. Hence
$$
\gathered
\tsize\sum(-1)^i\operatorname{tr}(\phi:
\oh^i_{\{0\}}(\Bbb P^1, \bar i^!j_!\Cal E) \to
\oh^i_{\{0\}}(\Bbb P^1, \bar i^!j_!\Cal E)) = 
\\
\tsize\sum(-1)^i\operatorname{tr}(\phi:
\oh^i_{\{0\}}(\Bbb C, i^!\Cal E) \to
\oh^i_{\{0\}}(\Bbb C, i^!\Cal E)) = 
\\
\tsize\sum(-1)^i\operatorname{tr}(\phi:\!
\oh^i_{\{0\}}(\Cal E) \to
\oh^i_{\{0\}}(\Cal E))\,,
\endgathered
\tag10.55
$$
i.e., the right hand side of the identity (10.52), is the contribution
of the fixed point 0 to the global trace (10.54).  We had argued initially that this
global trace coincides with the left had side in (10.52). To finish the proof of the
proposition, we need to show that the fixed point contribution at
$\infty$ in (10.54) is zero.  The map $f:\Bbb P^1\to \Bbb P^1$ contracts around
infinity, so  \lq\lq shrieking" to the contracting directions
becomes vacuous,  and restricting to the fixed point means
taking the stalk. Thus
$$
\tsize\sum(-1)^i\operatorname{tr}(\phi:
\oh^i((\bar i^!j_!\Cal E)_{\{\infty\}}) \to
\oh^i((\bar i^!j_!\Cal E)_{\{\infty\}})
\tag10.56
$$
is the fixed point contribution at infinity. This is also the fixed point
contribution at $(\infty,0)$ for the the pair $(j_!\Cal E, j_!\phi)$ and the map
$\bar f:\Bbb P^1\times V \to \Bbb P^1 \times V$ induced by $f$; to see this,
note that around
$(\infty,0)$, $\bar f$ is contracting only in the direction of $\Bbb P^1$. We now
use one of the other formulas in \cite{GM} to get an alternate expression for the
fixed point contribution at $(\infty,0)$. The complex $\bold A_5$ in \cite{GM} is
obtained by restricting the sheaf $j_!\Cal E$ to the expanding directions -- in
our case $\{\infty\}\times V$ -- then \lq\lq shrieking" to the fixed point.  But
$j_!\Cal E$ is the sheaf $\Cal E$ extended by 0 along  $\{\infty\}\times V$. In the
present setting, then, the entire complex $\bold A_5$ reduces to 0, so (10.56)
vanishes.
 
\enddemo

\vskip 2\jot

\enddocument